\documentclass{amsart}
\usepackage{amsmath, latexsym}
\usepackage{amssymb,amsthm}
\usepackage{amscd,graphics}

\newtheorem{proposition}{Proposition}[section]
\newtheorem{lemma}[proposition]{Lemma}
\newtheorem{sublemma}[proposition]{Sublemma}
\newtheorem{claim}[proposition]{Claim}
\newtheorem{theorem}[proposition]{Theorem}

\newtheorem{corollary}[proposition]{Corollary}

\theoremstyle{remark}

\newtheorem*{remark}{Remark}

\theoremstyle{definition}
\newtheorem*{de}{Definition}

\def\q{q}
\def\rmd{d}

\def\admc{{\mathcal{AC}}} 
\def\badmc{{\mathbf{AC}}}
\def\admm{{\mathcal{AM}}} 
\def\badmm{{\mathbf{AM}}} 

\def\ball{\mathbf{B}} 
\def\base{\Phi} 
\def\bchi{\mathbf{\chi}}
\def\bF{F_{\sharp}} 
\def\basin{\mathcal{B}}
\def\bound{\rho}
\def\chicl{\chi_{c}^{-}} 
\def\chicu{\chi_{c}^{+}}
\def\chiul{\chi_{u}^{-}}
\def\chiuu{\chi_{u}^{+}}
\def\chiudif{\chi_{u}^{\Delta}}
\def\chicdif{\chi_{c}^{\Delta}}

\def\cone{\mathbf{S}}
\def\Crit{\mathbf{C}}
\def\c{\mathbf{c}}
\def\criticalset{\mathcal{C}}
\def\curve{\gamma}
\def\Curve{\Gamma}
\def\cE{\mathcal{E}}

\def\D{\mathcal{D}} 
\def\disk{\mathbf D}

\def\eu{{\mathbf e}^{u}}
\def\ec{{\mathbf e}^{c}}
\def\e{{\mathbf{e}}} 
\def\E{\mathbf{E}}

\def\i{\mathbf{i}}

\def\j{\mathbf{j}}
\def\Jet{\mathbb{J}}

\def\k{\mathbf{k}}
\def\kone{\mathbf{I}}

\def\lattice{\mathbb{L}}
\def\leb{\mathbf{m}}
\def\lmax{\Lambda\subg}

\def\man{M} 
\def\mult{\mathbf{N}}

\def\mom{\mathcal{M}}
\def\mon{{\mathcal N}}

\def\one{\mathbf{1}}

\def\ph{\mathcal{PH}} 
\def\phys{\mathcal{R}} 
\def\proj{\Pi}
\def\proi{\pi}
\def\regs{\mathcal{Q}}

\def\Real{\mathbb{R}}

\def\sube{_{\epsilon}}
\def\subg{_{g}}
\def\subz{_{\chi}}

\def\ss{\mathcal{S}}
\def\s{\mathbf{s}}
\def\spl{\mathbf{Q}}

\def\tangentbundle{T\!M}

\def\t{\mathbf{t}}
\def\tmu{\tilde{\mu}}

\def\torus{\mathbb{T}}

\def\U{\mathcal{U}}

\def\vol{\mathbf{m}}
\def\v{\mathbf{v}}

\def\cX{\mathcal{X}}
\def\X{\mathbf{X}}

\def\y{{\mathbf y}}
\def\Y{{\mathcal Y}}

\def\Z{\mathbb{Z}}

\title[Partially hyperbolic endomorphisms]{Physical measures for partially hyperbolic
surface endomorphisms}

\author{Masato TSUJII}

\address{Department of Mathematics, Graduate School of Science, Hokkaido University}
\date{\today}

\begin{document}

\begin{abstract}
 We consider dynamical systems generated by partially hyperbolic
surface endomorphisms  of class $C^r$  with one-dimensional unstable subbundle.
As the main result, we prove that such a dynamical system generically admits 
finitely many  ergodic 
 physical measures whose union of basins of attraction has total Lebesgue measure,
provided that  $r\ge 19$.  
\end{abstract}

\maketitle




\section{Introduction}\label{sec-1}

In the study of smooth dynamical systems from the standpoint of ergodic theory, one of the most fundamental question is whether the following preferable picture is true for almost all of them: The asymptotic
distribution of the orbit for Lebesgue almost every initial point exists and coincides with one of the finitely many ergodic invariant measures that are given for the dynamical system. The answer is expected to be affirmative in general\cite{Palis}.  However it seems far beyond the scope of researches at present to answer the question in the general setting. The purpose of this paper is to provide an affirmative answer to the question in the case of partially hyperbolic endomorphisms on surfaces with one-dimensional unstable subbundle. 

Let $\man$ be the two-dimensional torus $\torus=\Real^{2}/\Z^{2}$ or, more generally,  a region on the torus $\torus$ whose boundary consists of finitely many simple closed $C^{2}$curves:  e.g. an annulus $(\Real/\Z) \times [-1/3,1/3]$. 
We equip $M$ with the Riemannian metric
$\|\cdot\|$ and the Lebesgue measure~$\leb$  that are induced by the standard ones on the Euclidean space $\Real^{2}$ in an obvious manner.    
We call a $C^{1}$mapping $F:\man\to \man$ a {\em partially hyperbolic endomorphism}
if there are positive constants
$\lambda$ and $c$ and  a continuous decomposition of the tangent bundle
$\tangentbundle= \E^{c}\oplus \E^{u}$ with $\dim \E^{c}=\dim \E^{u}=1$
such that
\begin{itemize}
\item[(i)] $\|DF^{n}|_{\E^{u}(z)}\|> \exp(\lambda n-c)$ and 
\item[(ii)] $\|DF^{n}|_{\E^{c}(z)}\|< \exp(-\lambda n+c)\|DF^{n}|_{\E^{u}(z)}\|$
\end{itemize}
for all $z\in \man$ and $n\ge 0$.
The subbundles $\E^c$ and $\E^u$ are called the central and unstable
subbundle, respectively. Notice that we do not require these subbundles to be invariant 
in the definition, though the central subbundle $\E^c$ turns out to be forward invariant from the condition (ii).
The totality of partially hyperbolic
$C^r$endomorphisms on $M$ is an open subset in the space $C^r(M,M)$, provided $r\ge 1$.

An invariant Borel probability measure $\mu$ for a dynamical system  $F:M\to M$ is said to be a {\em physical measure} if its basin of attraction,
\[
\basin(\mu)=\basin(\mu;F):=\left\{z\in \man\; \left|\;\;
\frac1n\sum_{i=0}^{n-1}\delta_{F^{i}(z)}\to
\mu
\mbox{ weakly as }n\to
\infty\right.\right\},
\]
has positive Lebesgue measure. 
One of the main results of this paper is 
\begin{theorem}\label{theorem-main} A partially
hyperbolic $C^{r}$endomorphism on $\man$ generically admits finitely many ergodic physical measures whose
union of basins of attraction has total Lebesgue measure, provided
that
$r\ge 19$. 
\end{theorem}
More detailed versions of this theorem will be given in the next section.
Here we intend to explain the new idea behind  the results of this paper. The readers should
notice that we do {\em not} (and  will {\em not}) claim  that the physical measures in the theorem above  are
hyperbolic.  Instead, we will show that the physical measures for  generic partially hyperbolic
endomorphisms have nice properties even if they are not hyperbolic. This is the novelty of the argument in this paper. 

Let us consider a
partially hyperbolic endomorphism  $F$ on $\man$. The  
 Lyapunov exponent of $F$ takes two distinct values at each point: The
larger is positive and the smaller   
 indefinite. The latter is called  the central Lyapunov exponent, as it is attained by the
vectors in the central subbunle. An invariant measure for $F$ is  hyperbolic if the central
Lyapunov exponent is non-zero at almost every pont with respect to it. 
In former part of this paper, we study hyperbolic invariant measures
for partially hyperbolic endomorphisms using the techniques in the Pesin theory or the smooth ergodic theory. And, as the conclusion, we show that  
the following hold under some generic conditions on~$F$: {\em  For any $\chi>0$, there are only finitely many ergodic physical measures whose
central Lyapunov exponents are larger than~$\chi$ in absolute value. Further, if  the complement $X$ of the union of the
basins of attraction of such physical measures has positive Lebesgue measure, and if a measure $\mu$ is a weak limit point  of the sequence
\begin{equation}
\label{eqn:avemu}
\frac{1}{n}\sum_{i=0}^{n-1}(\leb|_{X})\circ F^{-i},\quad n=1,2,\dots,\qquad\mbox{($\leb|_{X}$ : the restriction of $\leb$ to $X$),}
\end{equation}
then the absolute value of the central Lyapunov
exponent is not larger than~$\chi$ at almost every point with respect to $\mu$.} Though these facts
are far from trivial, the argument in the proof does not deviate far from the existing ones
in the smooth ergodic theory. 

\begin{figure}[htbp]
\begin{center}
\includegraphics{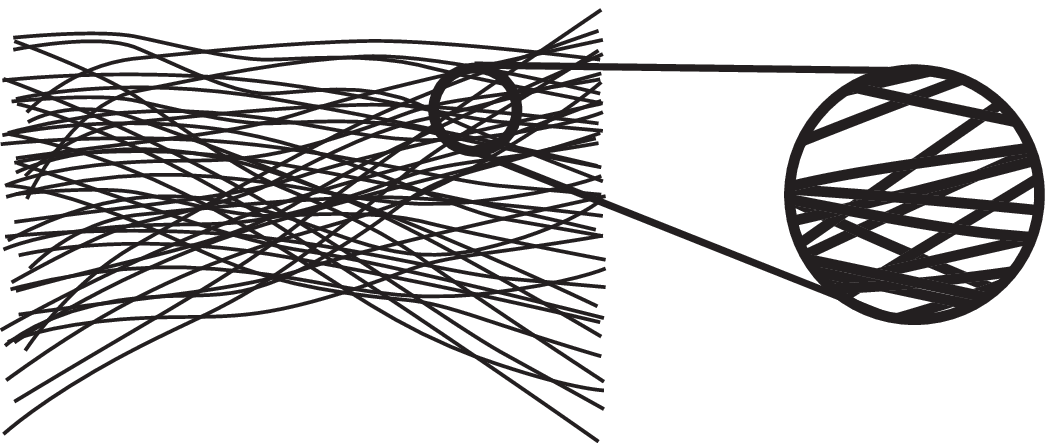}
\caption{The curve $F^n(\curve)$}
\label{fig1}
\end{center}
\end{figure}

The key claim in our argument is the following: {\em If the number $\chi$ is small enough and if $F$ satisfies some additional generic conditions, then such 
measure~$\mu$ as in the preceding paragraph is absolutely continuous
with respect to the Lebesgue measure~$\leb$. Further, the density $d\mu/d\leb$ satisfies some regularity conditions (from which  we can conclude 
theorem \ref{theorem-main}).}  This claim might appear 
unusual, since the measure $\mu$ may have neutral or even negative central Lyapunov exponent
while  we usually meet absolutely continuous invariant measures as a consequence of  expanding
property of dynamical systems in all directions. We can explain it intuitively  as
follows: As a consequence of the dominating expansion in the  unstable
directions $\E^u$, the measure~$\mu$ should have  some smoothness or uniformity in those
directions. In fact, we can show that the natural extensions of $\mu$ and its ergodic components
to the inverse limit are absolutely continuous along the (one-dimensional) unstable manifolds. So, for
each ergodic component $\mu'$ of $\mu$, we can cut a curve~$\curve$ out of an    unstable
manifold so that 
$\mu'$ is attained  as a weak limit point of the sequence 
$n^{-1}\sum_{i=0}^{n-1}\nu_{\curve}\circ F^{-i}$, $n=1,2,\dots$, where $\nu_{\curve}$ is 
a smooth measure on~$\curve$. Since $F$
 expands the curve~$\curve$ uniformly, the image $F^n(\curve)$ for large~$n$
should be a very long curve which is transversal to the central subbundle
$\E^c$.  Imagine to look into  a
small neighborhood of a point in the support of $\mu'$. The image $F^n(\curve)$ should
appear  as a bunch of short pieces of curve in that neighborhood. (Figure \ref{fig1}) 

The number of the pieces of curve should grow exponentially as $n$ gets large. And they would not 
concentrate  in the central direction strongly,  as  the central Lyapunov
exponent is nearly neutral almost everywhere with respect to~$\mu'$.  These suggest that the ergodic component~$\mu'$ should have some smoothness or uniformity in the central
direction as well as in the unstable direction and, so, the measure $\mu$  should be
absolutely continuous with respect to the Lebesgue measure~$\leb$.

On the technical side,  an important idea in the proof of the key claim is that we look at the
angles between the short pieces of curve mentioned above
rather than their positions.  As we  perturb the mapping $F$,  it
turns out that we can  control the angles between those pieces of curve to some extent, though we can not control 
their positions by the usual problem of interference. And we can show that the pieces of
curve satisfy some  transversality condition generically.  In order to show the conclusion of  the
key claim, we  relate that transversality condition to absolute continuity of  the measure $\mu$.
To this end, we make use of an idea in the  paper\cite{PeresSolo} by Peres and Solomyak with some modification.
We will illustrate
the idea in the beginning of section~\ref{sec-neu} by using a simple example.   Actually we have used
the same idea in our previous paper\cite{Tsujii}, which can be regarded as a study for this
work. Lastly, the author would like to note that the idea in \cite{PeresSolo}   can traces back
to the papers of Falconer\cite{F} and 
 Simon
and Pollicot\cite{Simon}.

\noindent{\bf Acknowledgment}:\;\; I would like to  thank  J\'er\^ome Buzzi, Viviane Baladi, Artur Avila and Mitsuhiro
Shishikura for valuable comments in writing this paper. 


\section{Statement of the main results}\label{s-state}
Let $\ph^{r}$ be the set of  partially
hyperbolic $C^{r}$endomorphisms on $\man$ and $\ph^{r}_{0}$
the subset of those without critical points. The subset $\mathcal R^{r}\subset \ph^{r} $ is the totality of mappings
$F\in \ph^{r}$
 that satisfy the following two conditions:
\begin{itemize}
\item[{\rm (a)}]  $F$
admits a finite collection of ergodic physical measures whose union of basins of attraction has total
Lebesgue measure on
$M$, and 
\item[{\rm (b)}] a physical measure for $F$ is absolutely continuous with respect
to the Lebesgue measure~$\leb$ if the sum of its Lyapunov exponents
is positive.
\end{itemize} 
In this paper, we claim that almost all partially hyperbolic
$C^r$endomorphisms on $\man$ satisfy the conditions (a) and
(b) above or, in other words, belong to the subset~$\mathcal R^r$. 
The former part of our main result is stated as follows:
\begin{theorem}\label{main_theorem1} {\rm (I)} \;The subset $\mathcal R^{r}$
is a residual subset in
$\ph^{r}$,  provided  $r\ge 19$. \hfill\ \linebreak
{\rm (II)}  The intersection $\phys^{r}\cap \ph^{r}_{0}$ is a
residual subset in 
$\ph^{r}_{0}$, provided $r\ge 2$. 
\end{theorem}
\noindent
The conclusions of this theorem mean  that the complement of the subset $\phys^r$ is  a meager subset in the
sense of Baire's category argument.  However the recent progress  in dynamical system theory has thrown serious
doubt that the notion of genericity  based on Baire's category argument may not have its literal meaning.  
In fact, it can happen that the dynamical systems in some meager subset  appear as 
 subsets with positive  Lebesgue measure in the  parameter spaces of typical families. For example, compare 
Jakobson's theorem\cite{Tsujii3} and the density of Axiom A\cite{Koz, Shen} in one-dimensional
dynamical systems.  For this reason, we dare to state our claim also in a measure-theoretical framework, though no
measure-theoretical definition that corresponds to  the notion of genericity has been firmly established yet.

Let $B$ be a Banach space. Let $\tau_v:B\to B$ be the
translation by $v\in B$, that is,  $\tau_v(x)=x+v$. A Borel finite measure $\mom$ on 
$B$ is said to be {\em quasi-invariant}
along a linear subspace $L\subset B$ if 
 $\mom\circ \tau_v^{-1}$ is equivalent to $\mom$ for any  $v\in 
L$. In the case where $B$ is finite-dimensional, a
Borel finite measure on~$B$ is equivalent to  the Lebesgue measure if and only if it is  quasi-invariant along the
whole space~$B$. But, unfortunately enough, it is  known that no Borel finite 
measures on an infinite-dimensional Banach space are quasi-invariant along the whole space\cite{Feldman}. This is
one of the  reasons why we do not have obvious  definitions for the concept like {\em Lebesgue
almost everywhere} in the cases of infinite-dimensional Banach spaces or Banach manifolds such as
the space
 $C^{r}(\man,\man)$. However, there may be  Borel finite measures
on
$B$ that are quasi-invariant along  dense subspaces. In fact,  on the
Banach space
$C^{r}(\man,\Real^{2})$, there exist Borel finite measures that are quasi-invariant along the
dense subspace
$C^{r+2}(\man,\Real^{2})$. (Lemma~\ref{lem-mom}.) For integers $s\ge r\ge 1$, we will  denote by $\mathcal Q_s^r$
the set of Borel probability  measures on $C^{r}(\man,\Real^{2})$ that are quasi-invariant along the subspace 
$C^{s}(\man,\Real^2)$ and regard the
measures in these sets as substitutions for the (non-existing) Lebesgue measure. 

Let us consider the space  $C^{r}(\man,\torus)$  of $C^{r}$mappings from $\man$ to
the torus $\torus$, which contains the space $C^{r}(\man,\man)$ of $C^r$endomorphisms on $M$.
For a mapping $G$ in 
$C^{r}(\man,\torus)$, we consider the mapping 
\begin{equation}\label{PhiG}
\base_{G}:C^{r}(\man,\Real^2)\to C^{r}(\man,\torus),\quad F\mapsto G+F.
\end{equation}
We now introduce the following notions.
\begin{de} A  subset $\cX\subset C^{r}(\man,\man)$ is {\em shy}
with respect to a measure
$\mom$ on 
$C^{r}(\man,\Real^2)$ if 
$\base_{G}^{-1}(\cX)$ is a null subset with respect to $\mom$ for any $G\in
C^{r}(\man,\torus)$. 
\end{de}
\begin{de} A subset $\cX\subset C^{r}(\man,\man)$ is {\em timid} for the class $\mathcal Q_s^r$ of measures if $\mathcal Q_s^r$ is not empty and
 if  $\cX$ is shy with respect to {\em all} measures in  $\mathcal Q_s^r$. 
\end{de}
\begin{remark}
The former of the definitions above is a slight modification of that of shyness introduced by Hunt, Sauer and Yorke\cite{HSY1, HSY2}.  By the definitions, a subset $\cX\subset C^{r}(\man,\man)$ is shy with respect to some compactly supported measure~$\mom$ in the sense above if and only if the inverse image $\base_{G}^{-1}(\cX)$ for same (and thus any) $G\in
C^{r}(\man,\torus)$ is shy in the sense of Hunt, Sauer and Yorke. 
Note that a controversy to the notion of shyness is given in the paper\cite{SC} of Stinchcombe.
\end{remark}

Put $\ss^r:=\ph^{r}\setminus\phys^{r} $. The latter part of our main result is stated as follows.
\begin{theorem}\label{main_theorem2}
{\rm (I)}\; The subset  $\ss^r$ is shy with respect to a  measure $\mom_s$ in
$\mathcal Q_{s-1}^r$ if 
the integers $r\ge 2$ and  $s\ge r+3$ satisfy 
\begin{equation}\label{rs}
(r-2)(r+1)<(r-\nu-2)\left(2(r-3)-\frac{(2s-r-\nu+1)}{\nu}\right)
\end{equation}
for some  integer $3\le \nu \le r-2$. 
Moreover, $\ss^r$ is timid for $\mathcal Q_{s-1}^r$ if
 $r\ge 2$ and  $s\ge  r+3$ satisfy
 the
condition (\ref{rs}) with $s$ replaced by $s+2$ for some integer $3\le \nu \le r-2$.
\par\noindent
{\rm (II)}
 $\ss^{r}\cap \ph_{0}^{r}$  is timid  for $\mathcal Q_s^r$ for any
 $r\ge 2$ and $s\ge r+2$.
\end{theorem}
\begin{remark}The measure $\mom_{s}$ in the claim (I) above will be constracted explicitly as a Gaussian measure.
\end{remark}
\begin{remark} The inequality (\ref{rs}) holds for the combinations 
$(r,s,\nu)=(19,22,3)$ and $(21, 26, 3)$ for example.  But it does not hold for any $s\ge r+3$ and $3\le \nu \le r-2$ unless $r\ge 19$.
\end{remark}
As an advantage of the measure-theoretical notion of timidity introduced above, we can derive  the following corollary  on the families of mappings 
in $\ph^r$, whose proof is given in the appendix.
  Let us regard the
space 
$C^{r}(\man\times [-1,1]^k, \man)$ as that of  $C^{r}$families of endomorphisms on $\man$
with parameter space $[-1,1]^k$. We can introduce the notion of shyness and timidity for the Borel subsets in this space
in the  same way as we did for those in  $C^{r}(\man, \man)$. Let $\leb_{\Real^k}$ be the
Lebesgue measure on
$\Real^k$. 
\begin{corollary}\label{lm-trans-meas} The set  of 
$C^{r}$families $F(z,\t)$ in $C^{r}(\man\times [-1,1]^k, \man)$ such that 
\[
\leb_{\Real^k}(\{\t\in [-1,1]^k\mid F(\cdot,\t)\in \ss^r\})>0
\]  is
timid with respect to the class of Borel finite measure on  $C^{r}(\man\times [-1,1]^k, \Real^2)$ that is
quasi-invariant along the subspace $C^{s-1}(\man\times [-1,1]^k,
\Real^2)$, provided that the integers $r\ge 2$ and  $s\ge r+3$ satisfy
 the
condition (\ref{rs}) with $s$ replaced by $s+2$ for some integer $3\le \nu \le r-2$. 
\end{corollary}

We give a few comments on the main results above. The restriction  that the surface $\man$  is a region on
the torus is actually not very essential. We could prove
theorem~\ref{main_theorem1} with $M$ being a general compact
surface by 
modifying the proof  slightly. The main  reason for this
restriction is  the difficulty in generalizing the notion of shyness and timidity to the spaces of endomorphisms on
general compact surfaces. Since the definitions depend heavily on the linear structure
of the space $C^r(M,\Real^2)$,   we hardly
know how we can  modify these notions naturally so that it is consistent under the non-linear
coordinate transformations. The generalization or modification of these notions should be
an important issue in the future. Besides, the restriction on
$M$ simplifies  the proof considerably and does not exclude the interesting examples such as the
so-called Viana-Alves maps\cite{A,V}.

The assumptions 
on differentiability in the main theorems are crucial in our argument especially in the part where we 
consider the influence of the critical points on the dynamics.  We do not know whether they are technical ones or
not.  

As we called attention in the introduction, the main
theorems tell nothing about hyperbolicity of  the physical measures. Of course, it is natural to expect that the
physical measures are  hyperbolic generically. The author think it not too optimistic to expect that 
$\phys^r$ contains an open dense subset of $\ph^r$ in which the
physical measures for the mappings are hyperbolic and depend on the mapping continuously.

Generalization of the main theorems to partially hyperbolic diffeomorphisms on higher dimensional
manifolds is an  interesting  subject to study. Our argument on 
physical measures with nearly  neutral central Lyapunov exponent seems to be complementary to the recent
works\cite{AVB,BV} of Alves, Bonatti and Viana. However, as far as the author understand, there exist essential
 difficulties in the case where the dimension of the central subbundle is higher  than one.

The plan of this paper is as follows: We give some preliminaries
 in section \ref{sec-pre}. We first define some basic notations and then introduce the
notions of {\em admissible curve} and {\em admissible measure}, which play central rolls in our
argument. The former is taken from the paper\cite{V} of Viana with slight modification and the
latter is a corresponding notion for measures. Next we introduce two conditions on partially
hyperbolic endomorphisms, namely, {\em the transversality condition on unstable cones} and {\em
the no flat contact condition}. At the end of section \ref{sec-pre}, we shall give a  concrete
plan of the proof of the main theorems  using the terminology introduced in this section. 
In section~\ref{sec-hyp}, we study hyperbolic physical measures using the Pesin theory. Section
\ref{sec-dist} is devoted to basic  estimates on the distortion of the iterates of partially hyperbolic
endomorphisms. Then we go into the main part of this paper, which consists of three mutually independent
sections.  In section \ref{sec-neu}, we prove that a partially hyperbolic endomorphism belongs to the subset
$\phys^r$ if  it satisfies the two conditions above.
In section \ref{sec-trans}
and \ref{sec-nonflat} respectively,  we prove that each of the two conditions holds for almost all
partially hyperbolic endomorphisms.


\section{Preliminaries}\label{sec-pre}
In this section, we prepare some notations, definitions and basic lemmas that
we shall use frequently in the following sections. 
\subsection{Notations} \label{sec-notation}
For a tangent vector
$v\in TM$, we denote by $v^\perp$ the tangent vector that
is obtained by rotating $v$ by the right angle  in the
counter-clockwise direction. For two tangent vectors $u$ and $v$, we 
denote by $\angle(u,v)$ the angle between them even if they belongs to the tangent spaces at different points.
Let $\exp_{z}:T_{z}\torus\to \torus$ be the
exponential mapping, which is defined simply by $\exp_{z}(v)=z+v$ in our case.
 For a point $z$ in the torus $\torus$ (or in some metric space, more generally) and a 
positive number $\delta$, we denote by
$\ball(z,\delta)$ the open disk with center at $z$ and radius~$\delta$. 
Likewise, for a subset $X$, we denote  its open $\delta$-neighborhood by
$\ball(X,\delta)$. For a positive number~$\delta$, we define a lattice 
$\mathbb L(\delta)$ as the subset of 
 points $(x,y)$ in $\torus$ whose components, $x$ and $y$, are  multiples of $1/([1/\delta]+1)$, so that
the disks
$\ball(z,\delta)$ for points $z\in \lattice(\delta)$ cover the torus $\torus$.

Throughout this paper, we assume $r\ge 2$. Let $F:\man\to \man$ be a $C^{r}$mapping. We denote the critical set of
$F$ by $\criticalset(F)$. For
a tangent vector
$v\in T_zM$ at a point $z\in \man$,  we define
\begin{align*}
&D_{*}F(z,v)=\|DF_{z}(v)\|/\|v\|\qquad \mbox{if  $v\neq {\mathbf 0}$}\\
\intertext{
and}
&D^{*}F(z,v)=(\det DF_z)/D_{*}F(z,v)\qquad \mbox{if  $DF(v)\neq {\mathbf 0}$}.
\end{align*}
\begin{remark}If $v\neq {\mathbf 0}$ and $DF(v)\neq {\mathbf 0}$, we can take  orthonormal bases $(v/\|v\|, v^{\perp}/\|v^{\perp}\|)$ on $T_zM$ and $(DF(v)/\|DF(v)\|, (DF(v))^{\perp}/\|(DF(v))^{\perp}\|)$ on $T_{F(z)}M$. Then the representation matrix of $DF_z:T_zM\to T_{f(z)}M$ with respect to these bases is a lower triangular matrix with $D_{*}F(z,v)$ and $D^{*}F(z,v)$ on the diagonal.
\end{remark}
Note that we have
$
|D^{*}F(z,v)|=\left\|
(DF)^{*}(v^{*})\right\|/\left\| v^{*}\right\|
$ 
for  any cotangent vector  $v^{*}\neq {\mathbf 0}$ at 
$F(z)$ that is normal to $DF(v)$. We  shall write
$D_{*}F(v)$ and
$D^{*}F(v)$ for $D_{*}F(z,v)$ and $D^{*}F(z,v)$ respectively in places where the base point $z$ is
clear from the context. 

For a $C^r$mapping $F:M\to \Real^2$, the $C^r$norm of $F$ is defined by
\[
\|F\|_{C^r}=\max_{z\in M}\;\max_{0\le  a+b\le r}\;
\left\|
\frac{\partial^{a+b}F}{\partial^{a}x\partial^{b}y}(z)
\right\|
\]
where $(x,y)$ is the coordinate on $\torus$ that is induced by the standard one on
$\Real^2$. Similarly, 
for $C^r$mappings $F$ and $ G$ in $ C^r(M,\torus)$, the
$C^r$distance  is defined by 
\[
d_{C^r}(F,G)=\max_{z\in M}\;\max\left\{d(F(z),G(z)), \max_{1\le a+b\le r}
\left\|
\frac{\partial^{a+b}F}{\partial^{a}x\partial^{b}y}(z)-\frac{\partial^{a+b}G}{\partial^{a}x\partial^{b}y}(z)
\right\|\right\}.
\]

\subsection{Some open subsets in $\ph^r$}\label{cnbl}
In this subsection, we introduce some bounded open subsets in $\ph^r$ whose elements enjoy certain estimates uniformly. As we will see,  we can restrict
ourselves to such open subsets in proving the main theorems. This simplifies the argument considerably. 

 Let
$\ss_0^r$ be the subset of mappings
$F$ in
$\ph^r$ that {\em violate} either of the  conditions:
\begin{itemize}
\item[(A1)] The image $F(\man)$ is contained in the interior of $\man$;
\item[(A2)] The function
$z\mapsto \det DF_{z}$ has $0$ as its regular value;
\item[(A3)] The restriction of $F$ to the critical set $\criticalset(F)$  is
transversal to
$\criticalset(F)$.  
\end{itemize}
Notice that the condition (A2) and (A3) are trivial if the mapping $F$ has no critical points. 
To prove the following lemma, we have only to apply
Thom's jet transversality theorem\cite{GG} and its measure-theoretical
version\cite[Theorem C]{Tsujii2}.
\begin{lemma}\label{ss0}
The subset $\ss_{0}^r$  is a closed nowhere
dense subset in  $\ph^{r}$ and  shy with respect to any measure in
$\mathcal{Q}_s^r$ for  
$s\ge r\ge 2$. 
\end{lemma}
\begin{remark}The terminology in \cite{Tsujii2} is different from that in this
paper. But we can put theorem C and other results in \cite{Tsujii2} into our terminology without difficulty. 
\end{remark}
Consider a $C^{r}$mapping $\bF$ in
$\ph^{r}$ and let $T\man=\E^{c}\oplus
\E^{u}$ be a decomposition of the tangent bundle which satisfies the  conditions
in the definition that $F=\bF$ is a partially hyperbolic endomorphism. Notice that, though the
central subbundle
$\E^{c}$ is uniquely determined by the conditions in the definition, the
 unstable subbundle
$\E^{u}$ is {\em not}. Indeed any continuous subbundle 
transversal to
$\E^{c}$ satisfies the conditions in the definition possibly with different
constants
$\lambda$ and
$c$. Making use of this arbitrariness, we can assume that  
$\E^{u}$ is a $C^\infty$subbundle. Further, by taking $\E^{u}$ 
nearly orthogonal to
$\E^{c}$ and by changing the constants $\lambda$ and $c$, 
we can assume that there exist  positive-valued
$C^{\infty}$functions
$\theta^{c}$ and
$\theta^{u}$ on $M$ such that the cone fields
\begin{align*}
&\cone^{u}(z)=\{v\in T_{z}\man\setminus \{0\}\mid 
\angle(v,\E^{u}(z))\le \theta^{u}(z)\}\quad\mbox{and}\\
&\cone^{c}(z)=\{v\in T_{z}\man\setminus \{0\}\mid 
\angle(v^{\perp},\E^{u}(z))\le \theta^{c}(z)\}
\end{align*}
satisfy the following conditions at every point $z\in \man$:
\begin{itemize}
\item[(B1)] $\cone^{c}(z)\cap
\cone^{u}(z)=\emptyset$;
\item[(B2)] $\E^{c}(z)\setminus\{0\}$ is contained in the interior of the cone $\cone^{c}(z)$;
\item[(B3)] $D\bF(\cone^{u}(z))$ is contained in the interior of $\cone^{u}(\bF(z))$;
\item[(B4)] $(D\bF)_{z}^{-1}(\cone^{c}(\bF(z)))$ is contained in
the interior of 
$\cone^{c}(z)$;
\item[(B5)] For any $v\in \cone^{u}(z)$ and $n\ge 1$, we have
\begin{itemize}
\item[(i)] $\|D_{*}\bF^{n}(z,v)\|> \exp(\lambda n-c)$ \quad and
\item[(ii)]
$\|D^{*}\bF^{n}(z,v)\|<
\exp(-\lambda n+c)\|D_{*}\bF^{n}(z,v)\|$.
\end{itemize}
\end{itemize}

Suppose that the mapping $\bF$ does {\em not} belong to $\ss_0^r$. Then we can take positive constants $\lambda$, $c$,  a
small number $\rho>0$ and a large number $\Lambda>c$  such that the following conditions hold for
any
$C^r$mapping
$F$ satisfying $d_{C^r}(F,\bF)<2\rho$:
\begin{itemize}
\item[(C1)] The conditions (B3), (B4) and (B5) with $\bF$ replaced by $F$ hold;
\item[(C2)] The parallel translation of 
$\E^{c}(\bF(z))$ to
$F(z)$ is contained in 
$\cone^{c}(F(z))\cup \{0\}$ for any $z\in M$;
\item[(C3)]  $d(F(\man),\partial \man)>\rho$;
\item[(C4)]  The
function
$z\mapsto
\det DF_z$  has no critical points on  $\ball(\criticalset(F),\rho)$ and it holds $|\det
DF_{z}|>\rho\cdot d(z, \criticalset(F))$ for $z\in \ball(\criticalset(F),\rho)$;
\item[(C5)]  If  a point $z\in M$ satisfies $d(z,w_1)<\rho$ and $d(F(z),w_2)<\rho$ for some points  $w_1,w_2\in
\criticalset(F)$ and if $ v\in \cone^{u}(z)$,  the angle between
$DF (v)$ and the tangent vector of $\criticalset(F)$ at $w_2$ is larger than $\rho$;
\item[(C6)] $\|DF_{z}\|<\Lambda$ for any $z\in M$.
\end{itemize} 
We can choose countably many pairs of a $C^{r}$mapping $\bF$ in
$\ph^r\setminus\ss_0^r$ and a positive number $\rho$ as above so that the corresponding open subsets
\[
\U=\{F\in C^{r}(M,M)\mid d_{C^{r}}(\bF,F)<\rho\} 
\] 
cover 
$\ph^r\setminus\ss_0^r$. In order to prove the main theorems, theorem
\ref{main_theorem1} and~\ref{main_theorem2}, it is enough to prove them by restricting
ourselves   to each of such open subsets~$\U$. For this reason,  we henceforth fix a 
$C^r$mapping
$\bF$ in $\ph^r\setminus \ss_0$, 
subbundles $\E^{c}$ and 
$\E^{u}$,  $C^\infty$functions $\theta^{c}$ and $\theta^{u}$, cone fields
$\cone^{c}(\cdot)$ and $\cone^{u}(\cdot)$ and positive   numbers $\rho$,  $\Lambda$, $\lambda$ and $c$ as above, and  consider the mappings in the corresponding
open subset
$\U$.

\subsection{Remarks on the notation for constants}\label{ssremc}
In this paper, we shall introduce various constants that depend only on 
\begin{itemize}
\item 
the objects that we have just fixed at the end of the last subsection, and
\item  the integer $r\ge
2$.
\end{itemize}
  In order to distinguish such kind of
constants,  
 we make it as a rule to denote them by symbols with subscript $g$. 
Obeying this rule,  we shall    denote $\lambda\subg$, $c\subg$,
$\rho\subg$ and $\Lambda\subg$ for the constants $\lambda$, $c$, $\rho$ and $\Lambda$ 
hereafter. (And we will use the symbols $\lambda$, $c$, $\rho$ and $\Lambda$ for the other purpose.) Notice that, once we denote a constant by a symbol with subscript $g$, we mean that
it is a constant of this kind. 
In order to save symbols for constants, we shall frequently use a generic symbol
$C\subg$ for large positive constants of this kind. 
Note that the value of the constants denoted by $C\subg$ may be different from place to place
even in a single expression.  For instance,  ridiculous expressions like $2C\subg<C\subg$ can be true, though we
shall not really meet such ones. Also note that we shall omit the phrases on the choice of
the constants denoted by 
$C\subg$ in most cases. 

For example, 
we can take a constant $A\subg>0$ such that it holds
\begin{equation}\label{anglerel}
A\subg^{-1}\frac{|D^{*}F^{n}(z,w)|}{D_{*}F^{n}(z,w)}
\le
\frac{\angle(DF^{n}(u),
DF^{n}(v))}{\angle(u,v)}\le A\subg\frac{|D^{*}F^{n}(z,w)|}{D_{*}F^{n}(z,w)}
\end{equation}
for any $z\in \man$, $n\ge 1$ and $u,v,w\notin \cone^{c}(z)\cup\{\mathbf{0}\}$.
 We shall use the following relations frequently: For any 
$F\in
\U$,
$z\in \man$, 
$v\in
\cone^{u}(z)$ and
$n\ge 1$, we have
\begin{align}\label{eqn-c1}
&C\subg^{-1} \cdot d(z,\criticalset(F)) \le   |\det DF_z| \le
\exp(\lmax) \|D^{*}F(z,v)\|\le C\subg \cdot d(z,\criticalset(F)) ,\\ 
&C\subg^{-1}<\|D_{*}F^{n}(z,v)\|\le \|DF^{n}_{z}\|\le C\subg \|D_{*}F^{n}(z,v)\|\label{eqn-c2}
\intertext{and, if $z\notin \criticalset(F)$,  also} \label{eqn-c3}
& C\subg^{-1}\cdot \|D^{*}F^{n}(z,v)\|\le \|(DF^{n}_{z})^{-1}\|^{-1}\le 
\|D^{*}F^{n}(z,v)\|.
\end{align}

\subsection{Admissible curves}\label{ss-admc}
In this subsection, we introduce the notion of admissible curve. From the forward invariance of the unstable cones
$\cone^u$ or the condition (B3) with $\bF$ replaced by $F$, the mappings in the set~$\U$ preserve the class of $C^1$curves  whose tangent vectors belong to $\cone^u$. We shall
investigate such class of curves  and find a subclass  which is uniformly bounded in
$C^{r-1}$sense and  essentially invariant under the iterates of mappings in~$\U$. We shall call the
curves in this subclass admissible curves.

In this paper, we 
always assume that the curves  are  regular and parameterized by length.  Let
$\curve:[0,a]\to
\man$ be a
$C^{r}$curve  such that $\curve'(t)\in \cone^{u}(\curve(t))$ for $t\in [0,a]$. As we assume
$\|\curve'(t)\|\equiv 1$,  the second differential of $\curve$ is
written in the form
\[
\frac{d^2}{dt^2}\curve(t)= d^{2}\curve(t)\cdot  (\curve'(t))^{\perp}
\]
where $d^{2}\curve:[0,a]\to \Real$ is a $C^{r-2}$function. We define $d^{k}\curve(t)$ for $3\le k\le r$
as the
$(k-2)$-th differential of the function $d^{2}\curve(t)$. For convenience, we  will denote the differential $\curve'(t)$ by $d^1\curve(t)$.

Let $F_{*}\curve:[0,a']\to \man$
be the image of the curve~$\curve$ under a mapping $F\in \U$.
Notice that $F_{*}\curve$ is not simply the composition  $F\circ \curve$, because we assume $F_{*}\curve$ to be
parameterized by length. The right relation between
$\curve
$ and
$F_{*}\curve$ is given by
\begin{equation}\label{eqn-transfer-curve}
F_{*}\curve(p(t))=F(\curve(t))
\end{equation}
where $p:[0,a]\to [0,a']$ is the unique
$C^{r}$diffeomorphism satisfying $p(0)=0$ and 
$\frac{d}{dt}p(t)=D_{*}F(\curve(t),\curve'(t))$. Differentiating the both sides of
(\ref{eqn-transfer-curve}), we get the formula
\[
D_{*}F(\curve(t), \curve'(t))\cdot
(F_{*}\curve)'(p(t))=DF_{\curve(t)}(\curve'(t))
\]
for $t\in [0,a]$.
Differentiating the both sides again and considering the components normal to
$(F_{*}\curve)'(p(t))$,
 we get 
\begin{equation}\label{induc-diff}
d^{2}F_{*}\curve(p(t))=\frac{D^{*}F(\curve(t),\curve'(t))}{D_{*}F(\curve(t),\curve'(t))^2}
\cdot
d^{2}\curve(t)+\frac{Q_{2}(\curve(t),\curve'(t);F)}{D_{*}F(\curve(t),\curve'(t))^3}
\end{equation}
where 
\[
Q_{2}(a,b;F)=(D^2 F_{a}(b,b), (DF_{a}(b))^{\perp}).
\]
Note that $Q_{2}(a,b;F)$ is a polynomial of the components of the unit vector $b$ whose coefficients are
polynomials  of the differentials of $F$ at $a$ up to the second order. 
Likewise, examining  the
 differentials of the both sides of (\ref{induc-diff}) by using the relation
\[
\frac{d}{dt}(D_*F(\curve(t),\curve'(t)))=\frac{\frac{d}{dt}(\|DF_{\curve(t)}(\curve'(t))\|^2)}{2\cdot D_*F(\curve(t),\curve'(t))},
\]
we obtain, for $3\le k\le r$,
\begin{equation}\label{eqn-inductive}
d^{k}F_{*}\curve(p(t))=\frac{D^{*}F(\curve(t),\curve'(t))}{D_{*}F(\curve(t),\curve'(t))^{k}}
\cdot  d^{k}\curve(t)+
\frac{Q_{k}(\curve(t),\curve'(t),
\{d^{i}\curve(t)\}_{i=2}^{k-1};F)}{D_{*}F(\curve(t),\curve'(t))^{3k-3}}
\end{equation}
where $Q_{k}(a,b,\{c_{i}\}_{i=2}^{k-1};F)$ is a  polynomial of the components of
the unit vector
$b$ and the scalars $c_{i}$ whose coefficients are polynomials of the differentials of $F$ 
at~$a$ up to the $k$-th order. 
\begin{remark}  In addition,  we can check that $Q_{k}(a,b,\{c_{i}\}_{i=2}^{k-1};F)$ for $2\le k\le r$ is written in  the form
\[
D_{*}F(a,b)^{2k-3}\cdot v^{*}((D^{k}F)_{a}(b,b,\cdots,b))+
{\tilde{Q}}_{k}(a,b,\{c_{i}\}_{i=2}^{k-1};F)
\]
where  $v^{*}$ is a unit cotangent vector at the point 
$F(a)$ that is normal to 
$DF_{a}(b)$ and $\tilde{Q}_{k}(a,b,\{c_{i}\}_{i=2}^{k-1};F)$ is a
polynomial of the components of 
$b$ and the scalars $c_{i}$  whose coefficients are 
polynomials of the differentials of $F$ at $a$ up to the $(k-1)$-th order.
\end{remark}

Fix an integer $n\subg>0$ such that $n\subg \lambda\subg - c\subg>0$. Then we have 
\begin{lemma} \label{lem-admc-def}
There exist constants 
$K\subg^{(k)}>1$ for  $2\le k\le r$ such that, if a curve $\curve:[0,a]\to \man$ of class $C^{r}$ satisfies
\begin{itemize}
\item[\rm{(i)}]$\curve'(t)\in
\cone^{u}(\curve(t))$ for  $t\in
[0,a]$,
\item[\rm{(ii)}] $|d^{k}\curve(t)|\le K\subg^{(k)}$ for  $2\le k\le r$ and $t\in
[0,a]$,
\end{itemize}
then   $F^{n}_{*}\curve$ for  $n\ge
n\subg$ satisfies the same conditions.
\end{lemma}
\begin{proof} Consider a $C^{r}$curve $\curve:[0,a]\to \man$ that satisfies the conditions (i) and (ii), and let $F^n_*\curve:[0,a_n]\to \man$ be its image under the iterate $F^n$. From the formulae (\ref{induc-diff}) and  (\ref{eqn-inductive}), we can see that the following holds for $n\subg \le n\le 2n\subg$:
\begin{equation}\label{eqn:clo}
|d^k F^n_*\curve(p_n(t))|\le \frac{|D^{*}F(\curve(t),\curve'(t))|}{D_{*}F(\curve(t),\curve'(t))^{k}}
\cdot  |d^{k}\curve(t)|+R(n\subg,K\subg^{(2)}, \dots, K\subg^{(k-1)})
\end{equation}
where $p_n:[0,a]\to [0,a_n]$ is the unique diffeomorphism satisfying $p(0)=0$ and $\frac{d}{dt}p(t)=DF^n_*(\curve(t),\curve'(t))$ and where $R(n\subg, K\subg^{(2)}, \dots, K\subg^{(k-1)})$ is a constant that depends only on $n\subg, K\subg^{(2)}, \dots, K\subg^{(k-1)}$ besides the objects that we have already fixed at the end of subsection ~\ref{cnbl}. The coefficient of $|d^{k}\curve(t)|$ on the right hand side of the inequality (\ref{eqn:clo}) is smaller than $\exp(-n\subg \lambda\subg + c\subg)<1$ from the condition (C1) and the choice of $n\subg$. Thus, if we take large $K\subg^{(k)}$ according to the choice of the constant $K\subg^{(2)},\dots, K\subg^{(k-1)}$ in turn for $2\le k\le r$, the conclusion of the lemma holds for $n\subg \le n\le 2n\subg$. And, employing this repeatedly,  we obtain the conclusion for all $n\ge n\subg$.
\end{proof}
Henceforth we fix the constants $K\subg^{(k)}$, $2\le k\le r$, in 
lemma~\ref{lem-admc-def}. Now we put
\begin{de}
A $C^{r-1}$curve
$\curve:[0,a]\to \man$ is called {\em an admissible curve} if it satisfies the conditions
\begin{itemize}
\item[\rm{(a)}]$\curve'(t)\in
\cone^{u}(\curve(t))$ for  $t\in
[0,a]$,
\item[\rm{(b)}] $|d^{k}\curve(t)|\le K\subg^{(k)}$ for  $2\le k\le r-1$ and $t\in
[0,a]$, and 
\item[\rm{(c)}]
the function $d^{r-1}\curve$ satisfies  Lipschitz condition with the constant $K\subg^{(r)}$:
\begin{equation*}
\left|d^{r-1}\curve(t)-d^{r-1}\curve(s)\right|\le K\subg^{(r)}|t-s|\qquad\mbox{ for any $0\le s<t\le a$.}
\end{equation*}
\end{itemize}
\end{de}
\begin{remark} When $r=2$, the condition (b) above is vacuous and the symbol $|\cdot |$ on the left hand side of the inequality in the condition (c)  should be understood as the norm on $\Real^2$. (Recall that we denote $d^1\curve(t)=\curve'(t)$.)
\end{remark}
Note that a $C^{r-1}$curve $\curve:[0,a]\to \man$ is an admissible curve if and only if  it belongs to the closure, in the space $C^{r-1}([0,a],\man)$,  of the set of $C^r$curves satisfying the conditions (i) and (ii) in lemma~\ref{lem-admc-def}. Thus we have, from lemma~\ref{lem-admc-def},
\begin{corollary}\label{cor33} If a $C^{r-1}$curve $\curve$ is admissible, so is $F^{n}_{*}\curve$ for $n\ge
n\subg$. 
\end{corollary}

\def\tadmc{\mathcal{A}}
\def\tbadmc{\mathbf{A}}

For  a positive number $a$, let $\tadmc(a)$ be  the set of $C^1$curves $\curve:[0,a]\to M$ of length~$a$ such that 
$\curve'(t)\in \cone^u(\curve(t))$ for $t\in [0,a]$. For a subset $J\subset (0,\infty)$, we define $\tadmc(J)$ as the disjoint union of $\tadmc(a)$ for $a\in J$:
\[
\tadmc(J):=\mathop{\amalg}_{a\in J}\tadmc(a).
\]
Also we define
\[
\tbadmc(J):=\mathop{\amalg}_{a\in J}(\tadmc(a)\times [0,a])\subset \tadmc(J)\times \Real
\]
We can regard  the space $\tadmc((0,\infty))$ as the totality of $C^1$curves whose length are finite and  whose tangent vectors are contained in the unstable cone $\cone^u$. 
From the condition (C1) in the choice of the open neighborhood $\U$, each mapping $F\in \U$ naturally acts on the space $\tadmc((0,\infty))$:
\[
F_*:\tadmc((0,\infty))\to \tadmc((0,\infty)),\quad \curve \in \tadmc(a)\mapsto F_*\curve\in \tadmc(p(a))
\]
and also on the space $\tbadmc((0,\infty))$:
\begin{align*}
F_*:\tbadmc((0,\infty))&\to \tbadmc((0,\infty)),\\
 (\curve,t)  \in \tadmc(a)\times [0,a] &\mapsto (F_*\curve, p(t))\in \tadmc(p(a))\times [0,p(a)]
\end{align*}
where $p:[0,a]\to [0, p(a)]$ is the unique diffeomoprhism satisfying $p(0)=0$ and $\frac{d}{dt}p(t)=D_*F(\curve(t),\curve'(t))$.

For a positive number $a$,  let $\admc(a)\subset \tadmc(a)$ be the set of admissible curves of length $a$ and, for a subset $J\in (0,\infty)$, we put \[
\admc(J):=\mathop{\amalg}_{a\in J}\admc(a) \subset \tadmc(J)
\quad \mbox{and} \quad 
\badmc(J):=\mathop{\amalg}_{a\in J}\admc(a)\times[0,a] \subset \tbadmc(J).
\]
Note that $\admc(a)$ is a compact subset of $C^{r-1}([0,a],\man)$. 

We equip the space  $\admc((0,\infty))$ with the distance $d_\admc$ defined by
\[
d_{\admc}(\curve_1,\curve_2)=\|\curve_2-\curve_1\|_{C^{r-1}}+ C\cdot |a_2-a_1|
\]
for 
$(\curve_i,t_i)\in \admc(a_i)$,
$i=1,2$, where $\|\curve_2-\curve_1\|_{C^{r-1}}$ is 
\[
\max_{0\le
\theta\le \min\{a_1,a_2\}} \left\{\; d(\curve_2(\theta),\curve_1(\theta)),\;
\angle(\curve_1'(\theta),\curve_2'(\theta)),\;
\max_{2\le k\le r-1} \left|d^{k}\curve_2(\theta)-d^{k}\curve_1(\theta)\right|\;\right\}
\]
and the constant $C$ is defined by 
\[
C=\left(\max_{2\le k\le
r} K\subg^{(k)}
\right). 
\]
Note that the constant $C$ above is chosen so that $d_{\admc}$ satisfies the axiom of distance. We equip the space $\badmc((0,\infty))$ with the distance
\[
d_{\badmc}((\curve_1,t_1),(\curve_2,t_2))=d_{\admc}(\curve_1, \curve_2)+|t_2-t_1|
\]
for 
$(\curve_i,t_i)\in
\admc(a_i)\times[0,a_i]$,
$i=1,2$. It is not difficult to check that the space $\admc(0,\infty)$ and $\badmc((0,\infty))$ with these distances are complete separable metric space and that
the subsets $\admc(J)\subset \admc((0,\infty))$ and $\badmc(J)\subset \badmc((0,\infty))$ for a subset  $J\subset (0,\infty)$ is compact if and only if $J$ is compact. 

From corollary \ref{cor33}, the iterate $F_{*}^{n}$ of the mapping $F_*:\tadmc((0,\infty))\to \tadmc((0,\infty))$ (resp. $F_*:\tbadmc((0,\infty))\to \tbadmc((0,\infty))$ for any $n\ge n\subg$ carry the subset 
$\admc((0,\infty))$ (resp. $\badmc((0,\infty))$) into itself. 
Further we have, for any $n\ge n\subg$ and $a>0$, 
\begin{equation}\label{eqn:fnaj}
F^n_{*}(\admc([a,\infty)))\subset \admc([a\exp(\lambda\subg n-c\subg),\infty))
\end{equation}
 and 
\begin{equation}\label{eqn:fnai}
F^n_{*}(\badmc([a,\infty)))\subset \badmc([a\exp(\lambda\subg n-c\subg),\infty)).
\end{equation}
We  define the
mapping
$\proj:\tbadmc((0,\infty))\to \man$ and $\proi:\tbadmc((0,\infty))\to \admc((0,\infty))$ by $\proj(\curve,t)=\curve(t)$ and $\proj(\curve,t)=\curve$. 
Obviously we have the commutative relations:
\begin{equation}\label{cdF}
\begin{CD}
\tbadmc((0,\infty)) @> F_*>> \tbadmc((0,\infty))\\
@V \proj VV @V \proj VV\\
M @>F >> M
\end{CD}
\qquad
\begin{CD}
\tbadmc((0,\infty)) @> F_*>> \tbadmc((0,\infty))\\
@V \proi VV @V \proi VV\\
\tadmc((0,\infty)) @>F_* >> \tadmc((0,\infty)).
\end{CD}
\end{equation}

\subsection{Admissible measures}\label{ss-admm}
In this subsection, we are going to introduce the notion of admissible measure. First we introduce this notion in a simple case. Let $\curve:[0,a]\to \man$ be an admissible curve and, for $n\ge 1$, let $p_n:[0,a]\to [0,a_n]$ be the unique diffeomorphism that satisfies $p_n(0)=0$ and $\frac{d}{dt}p_n(t)=D_*F^n(\curve(t), \curve'(t))$ for $t\in [0,a]$. Since mappings $F\in \U$ act on the admissible curves as uniformly expanding mappings with  uniformly bounded distortion, a standard argument on the iterations of uniformly expanding mappings gives
\begin{lemma} The mapping $p_n$ satisfies $\frac{d}{dt}p_n(t)\ge \exp(\lambda\subg n-c\subg)$
and 
\[
\log \frac{dp_n}{dt}(t)-\log \frac{dp_n}{dt}(s)\le C\subg\quad\mbox{for $t,s\in [0,a]$,}
\]
where $C\subg$ is the kind of constant that we mentioned in subsection \ref{ssremc} and, especially, does not depend on the mapping $F\in \U$, the admissible curve $\curve$ nor  $n\ge n\subg$. 
\end{lemma}
We say that a measure $\mu$ on an interval $I\subset \Real$ has {\em Lipschitz logarithmic density with constant $L$} if $\mu$ is written in the form $\mu=\varphi\cdot  \leb_{\Real}|_{I}$ where $\varphi:I\to \Real$  is a positive-valued function satisfying 
\[
|\log \varphi(t)-\log\varphi(s)|\le L|t-s|\quad\mbox{for any $t,s\in I$}
\]
and $\leb_{\Real}|_{I}$ is the restriction of the Lebesgue measure on $\Real$ to $I$.
Note that the sum (or integration) of measures on an interval $I$ having Lipschitz logarithmic density with constant $L$ again the same property. Form the lemma above, we can obtain 
\begin{corollary}\label{cor:admcd}
There is a positive constant $L\subg$ such that, if a measure $\mu$ on $[0,a]$ has Lipschitz logarithmic density with constant $L\subg$, then so does the measure $\mu\circ p_n^{-1}$ on $[0,a_n]$ for any $n\ge n\subg$,  any $F\in \U$ and any admissible curve $\curve:[0,a]\to \man$.
\end{corollary}
We henceforth fix the constant $L\subg$ for which the claim of  corollary \ref{cor:admcd} holds. And we say that a measure $\nu$ on $\man$ is an admissible measure on an admissible curve $\curve:[0,a]\to \man$ if $\nu= \mu\circ \curve^{-1}$ for a measure $\mu$ on $[0,a]$ that has Lipschitz logarithmic density with constant $L\subg$. It follows from corollary \ref{cor:admcd} that
\begin{corollary} 
If a measure $\nu$ is an admissible measure on an admissible curve $\curve:[0,a]\to \man$, then,  for $n\ge n\subg$ and $F\in \U$,  the measure $\nu\circ F^{-n}$ is an admissible measure on the admissble curve $F^n_*\curve$.
\end{corollary}
We have introduced the notion of admissible measure on a single curve and seen that the iterates of mappings $F\in \U$ preserve such class of measures. Now we are going to introduce more general definitions. 
Let $\Xi_{\badmc}$ be the measurable partition of the space  
 $\badmc((0,\infty))$  into
the subsets
$\{\curve\}\times [0,a]$ for 
$a>0$ and $\curve
\in\admc(a)$. In other words, we put $\Xi_{\badmc}=\pi^{-1}\varepsilon$ where $\varepsilon$ is the measurable partition of $\admc((0,\infty))$ into individual points and $\pi$ is the mapping defined at the end of the last subsection. On each element $\xi=\{\curve\}\times [0,a]$ of the partition $\Xi_{\badmc}$, we
consider the measure $\leb_{\xi}$ that corresponds to  the Lebesgue measure on $[0,a]$ through the bijection
$(\curve,t)\mapsto t$.  For a Borel finite measure
$\tmu$ on $ \badmc((0,\infty))$, let
$\{\tmu_{\xi}\}_{\xi\in
\Xi_{\badmc}}$ be the conditional measures with respect to the partition
$\Xi_{\badmc}$. We put the following two defintions.
\begin{de} A Borel finite measure $\tmu$ on $\badmc((0,\infty))$ is said to be an {\em  admissible measure} if the conditional measures $\{\tmu_\xi\}_{\xi\in
\Xi_{\badmc}}$ has Lipschitz logarithmic density with constant $L\subg$,  $\tmu$-almost everywhere.
\end{de}
\begin{de} A Borel finite measure $\mu$ on $\man$ is said to {\em have an admissible lift} if there exists an admissible measure $\tmu$ on $\badmc((0,\infty))$ such that $\tmu\circ \Pi^{-1}=\mu$.  The measure $\tmu$ is said to be an {\em admissible lift} of the measure $\mu$. 
\end{de}

For a subset $J\subset (0,\infty)$, we denote, by
$\badmm(J)$, be the set of admissible measures that is supported on $\badmc(J)$ and, by $\admm(J)$, the set of measures on $M$ that have admissible lifts contained in $\badmm(J)$. Then we have
\begin{lemma}\label{lem-admc-ind} 
If a measure $\tmu$ belongs to
$\badmm([a,\infty))$ for some $a\ge 0$ and if $F\in \U$, then
 $\tmu\circ F_{*}^{-n}$ belongs to  $\badmm([a',\infty))$ for  $n\ge n\subg$  where
$a'=a\exp(\lambda\subg n-c\subg)>a$.
\end{lemma}
\begin{proof} The conditional measures of the measure $\tmu\circ F_*^{-n}$ with respect to the partition $\Xi_{\badmc}$ are given as integrations of the images of the conditional measures $\{\tmu_{\xi}\}_{\xi\in \Xi_{\badmc}}$ under the mapping $F_{*}^n$. From corollary \ref{cor:admcd} and the fact noted just above it, they have Lipschitz logarithmic density with constant $L\subg$. Hence $\tmu\circ F_*^{-n}$ is an admissible measure. The claim of the lemma follows from this and (\ref{eqn:fnai}). 
\end{proof}
\begin{corollary}\label{cor-admc-ind} If $\mu\in \admm([a,\infty))$ for some $a> 0$ and if $F\in \U$, then
the measure  $\mu\circ F^{-n}$ belongs to  $\admm([a',\infty))$  for  $n\ge n\subg$ where
$a'=a\exp(\lambda\subg n-c\subg)> a$. Especially, if an invariant
measure for
$F\in
\U$ has an admissible lift, it belongs to $\admm([a,\infty))$ for any $a>0$.
\end{corollary}

\begin{lemma}\label{lik}
The subset $\badmm(J)$ for a closed subset $J\subset (0,\infty)$ is closed in the
space of Borel finite measures on
$\badmc((0,\infty))$. 
\end{lemma}
\begin{proof} For a real number $\epsilon$, we define the mapping $T_{\epsilon}
$ from $\admc((0,\infty))\times \Real$ to itself  by $T_{\epsilon}(\curve, t)=(\curve, t+\epsilon)$. Then a measure $\tmu$ on $\badmc((0,\infty))\subset \admc((0,\infty))\times \Real$ is admissible if and only if it satisfy
\[
\int_{\badmc((0,\infty))\cap T_{\epsilon}^{-1}(\badmc((0,\infty)))}  \varphi\circ T_{\epsilon}^{-1} \;d\tmu\le \exp(L\subg |\epsilon|)
\int_{\badmc((0,\infty))}  \varphi \; d\tmu
\]
for any non-negative-valued continuous function $\varphi$ on $\admc((0,\infty))\times \Real$.
For each non-negative-valued continuous function $\varphi$ on $\admc((0,\infty))\times \Real$ and $\epsilon\in \Real$, the set of Borel measures $\tmu$ that satisfy the inequality above and that are supported on $\badmc(J)$ is a closed subset  in the
space of Borel finite measures on
$\badmc((0,\infty))$. Hence so is their intersection, $\badmm(J)$.
\end{proof}
\begin{lemma}\label{lmata} $\admm([a,\infty))=\admm([a,2a])$ for $a>0$.
\end{lemma}
\begin{proof}
For $a>0$, let $\Delta_{a}:\badmc([a,\infty))\to \badmc([a,2a])$ be the mapping that brings an element 
$(\curve,t)\in
\admc(b)\times[0,b]$ to 
\begin{equation}\label{def-delta}
\Delta_{a}((\curve, t))=(\curve|_{[m(t), m(t)+(b/n)]}, t-m(t))\in \admc(b/n)\times[0,b/n]
\end{equation}
where $n=[b/a]$ and $m(t)=[t n/b](b/n)$. Then we have $\proj\circ \Delta_{a}=\proj$ and, for any
$\tmu\in\badmm([a,\infty)$, the image $\tmu\circ \Delta_a^{-1}$ belongs to 
$\badmm([a,2a])$. Thus we obtain the lemma.
\end{proof}
From the lemma above and lemma \ref{lik}, it follows
\begin{corollary}\label{lem-adm-cpt} The set $\admm([a,\infty))$ for $a>0$ is 
a compact subset in the space of Borel finite measures on $M$. Especially, for a mapping $F\in \U$, the set of
\hbox{$F$-invariant} Borel probability measures that have admissible lifts is compact.
\end{corollary}

Suppose that $P$ is  a small parallelogram on the torus $\torus$ whose center $z$ belongs to~$M$ and two
of whose sides are parallel to the unstable subspace
$\E^{u}(z)$. Then the restriction of the  Lebesgue
measure $\leb$  to $P$  has an admissible lift, provided that $P$ is sufficiently small. Moreover 
any linear combination of such measures has admissible lifts. Thus we obtain
\begin{lemma}\label{est-Leb} For any Borel finite measure $\nu$ on $M$ that is absolutely continuous with
respect to the Lebesgue measure $\leb$,
there exist a sequence $b_{n}\to +0$ and measures $\nu_{n}
\in 
\admm([b_{n},\infty))$ such that 
$|\nu-\nu_{n}|\to 0$ as $n\to \infty$. Further we can take the measures $\nu_n$ so that
the densities
$d\nu_{n}/d\leb$  are square integrable. 
\end{lemma}
The following is a consequence of the last two lemmas and
corollary~\ref{cor-admc-ind}.
\begin{lemma}\label{lem-lim} Let $F$ be a mapping in $\U$ and $\nu$  a probability measure on $\man$ that is
absolutely continuous with respect to the Lebesgue measure $\leb$. Then any limit
point of the sequence $n^{-1}\sum_{i=0}^{n-1}\nu\circ F^{-i}$ is contained in 
$\admm([a,\infty))$ for any $a>0$. Especially, physical measures for $F$ are contained in 
$\admm([a,\infty))$ for any $a>0$.
\end{lemma}
Finally we prove
\begin{lemma}\label{lem-adm2} 
Let $F$ be a mapping in $\U$. If  an $F$-invariant  Borel probability measure has an
admissible lift, so do its ergodic components. 
\end{lemma}
\begin{proof}
From corollary \ref{lem-adm-cpt}, it is enough to show the following claim: If an $F$-invariant  measure $\mu$ that
has an admissible lift splits into two non-trivial 
$F$-invariant measures $\mu_1$ and
$\mu_2$ that are totally singular with respect to each other, then
the measures $\mu_1$ and
$\mu_2$ have admissible lifts. We are going to show this claim.  From corollary
\ref{cor-admc-ind}, we can take an admissible lift $\tmu$ of $\mu$ that is supported on $\badmc([1,\infty))$.
Consider the mapping 
$G=\Delta_1\circ F^{n\subg}_{*}:\badmc([1,\infty))\to\badmc([1,2])$,  where $\Delta_1$ is the mapping defined by
(\ref{def-delta}). Then the measure $\tmu\circ G^{-1}$ is an admissible lift  of $\mu$.  Replacing $\tmu$ by  $\tmu\circ G^{-1}$, we can and do assume that
$\tmu$ is supported on $\badmc([1,2])$. 
From the assumption of the claim, we can take an $F$-invariant  Borel subset 
$X\subset M$ such that $\mu_1(M\setminus X)=\mu_{2}(X)=0$.
Then, by  the relation
 $F^{n\subg}\circ \proj=\proj\circ G$, the set  $\tilde X:=\pi^{-1}(X)$ is $G$-invariant.
Below we prove that $\tilde X$ is a $\Xi_{\badmc}$-set, that is, a union of elements of 
the partition $\Xi_{\badmc}$, modulo null subsets with respect to $\tmu$. This implies the claim above because
the restriction of the measure $\tmu$ to
$\tilde X$ is an admissible lift of
$\mu_1$. 

Put $\Xi_1=\Xi_\badmc$ and define the sequence $\Xi_n$, $n=1,2,\dots$ inductively by the relation
$\Xi_{n+1}=G^{-1}(\Xi_n)\vee \Xi_{1}$. Then $\Xi_{n}$ is increasing with respect to $n$ and the limit 
 $\bigvee_{n=1}^{\infty}\Xi_n$ is the measurable partition into individual points. Thus the conditional expectation 
$E(\tilde X|\Xi_n)$ with respect to $\tmu$ converges to the indicator function of $\tilde X$ as $n\to \infty$,
$\tmu$-a.e.  Note that the restriction of $G^n$ to each element of the partition $\Xi_n$ is a bijection onto an
element of
$\Xi_1$ and  its distortion is uniformly bounded. Hence, using the
assumption that
$\tmu$ is an admissible measure and the invariance of
$\tilde X$, we can see that the conditional expectation 
$E(\tilde X|\Xi_1)$ equals to the indicator function of $\tilde X$, or 
$\tilde X$ is a
$\Xi_{\badmc}$-set modulo null subsets with respect to $\tmu$.
\end{proof}

\subsection{The no flat contact condition}
In this subsection, we consider the influence of the  critical points on ergodic behavior of
partially hyperbolic endomorphisms. We first explain a problem
 that  the critical points may cause. And then we give a mild condition on the mappings in $\U$, {\em the no flat
contact condition}, which allows us to avoid that  problem.  In the last part of this paper, we will prove that
this condition  holds for almost all partially hyperbolic endomorphisms in $\U$.

Let us consider a mapping $F\in \U$.  We denote, by 
$\chi_{c}(z;F)<\chi_{u}(z;F)$, the Lyapunov exponents at  $z\in \man$.  For a Borel finite  measure $\mu$ on $\man$, we define
\begin{align*}
\chi_c(\mu;F)&=\frac{1}{|\mu|}\int \log \|DF|_{E^{c}(z)}\| \; d\mu(z)\quad \mbox{and} \\
\chi_u(\mu;F)&=\frac{1}{|\mu|}\int \log (|\det DF_z|/\|DF|_{E^{c}(z)}\|)\;  d\mu(z).
\end{align*} 
These are called the central and unstable  Lyapunov exponent of $\mu$, respectively. For an invariant probability measure $\mu$ for $F$, we have
\[
\chi_c(\mu;F)=\int \chi_c(z) d\mu(z)\quad \mbox{and} \quad
\chi_u(\mu;F)=\int \chi_u(z)  d\mu(z).
\]
Further, if~$\mu$ is  an ergodic
invariant measure for
$F\in
\U$, the Lyapunov exponents
$\chi_{c}(z;F)$ and
$\chi_{u}(z;F)$ take constant values $\chi_{c}(\mu;F)$ and $\chi_{u}(\mu;F)$ at $\mu$-almost every point $z$,
 respectively.

Let  
$\mu_{n}$, $n=1,2,\cdots$, a sequence of ergodic invariant probability measures for $F$ that have admissible
lifts. And suppose that $\mu_n$ converges weakly to some measure
$\mu_{\infty}$ as $n\to\infty$. Then $\mu_{\infty}$ has an admissible lift from 
corollary~\ref{lem-adm-cpt}. It is not
difficult to see that the  Lyapunov exponent
$\chi_{u}(\mu_{n};F)$  always converges to $\chi_{u}(\mu_{\infty};F)$. However, for the central
Lyapunov  exponent, we   only have the inequality
\begin{equation}\label{ine-Lyapunov}
 \limsup_{i\to \infty}\chi_{c}(\mu_{n};F)\le \chi_{c}(\mu_{\infty};F)
\end{equation}  
when $F$ has critical points, because the function $\log \|DF|_{E^{c}(z)}\|$ is not continuous at
the critical points. Though the strict inequality in (\ref{ine-Lyapunov}) is not likely to hold
often,  we can not avoid such cases in general. And, once the strict inequality holds,  the
ergodic behavior of
$F$ can be  wild by the influence of the critical point. 
\begin{remark}
It is not easy to construct  examples in which the strict inequality 
(\ref{ine-Lyapunov}) holds. For example, consider the
direct product of the quadratic mappings given in the paper\cite{HK}  and an angle
multiplying mapping
$\theta\mapsto d\cdot\theta$ on the circle. 
\end{remark}
\begin{remark} We could consider a similar but more general problem:
Suppose that a point $z\in M$ is generic for an invariant probability measure~$\mu$, that is, the sequence $n^{-1}\sum_{i=0}^{n-1}\delta_{F^i(z)}$ converges to~$\mu$ as $n\to \infty$. The problem is that the strict inequality $\chi_{c}(z;F)<\chi_{c}(\mu;F)$ can hold, though the equality $\chi_{u}(z;F)=
\chi_{u}(\mu;F)$ always holds. (If we did not assume the mapping $F$ to be partially hyperbolic, these relation would be looser.) We may call this kind of problems {\em Lyapunov irregularity}, as this is  the case where the so-called Lyapunov regularity condition\cite{O} does not hold. 
\end{remark}
In order to avoid the irregularity described above, we introduce 
 a mild condition: 
\begin{de}
We say that a mapping
$F\in
\U$ satisfies {\em the no flat contact condition} if
there exist positive constants
$C=C(F)$, $n_{0}=n_{0}(F)\ge n\subg$ and
$\beta=\beta(F)$ such that, for any admissible curve 
$\curve\in \admc(a)$ with $a>0$, $n\ge n_{0}$ and $\epsilon>0$, it holds
\begin{equation*}
\leb_{\Real}\left(\left\{t\in [0,a]\; ; \; d(F^{n}(\curve(t)),\criticalset(F))<\epsilon\right\}\right)<C\cdot 
\epsilon^{\beta}\max\{a,1\}
\end{equation*}
where $\leb_{\Real}$ is the Lebesgue measure on $\Real$. If $F$ has no critical points, we regard that
$d(z,\criticalset(F))=1$ for 
$z\in M$ and that
$F$ satisfies the no flat contact condition.
\end{de}

\begin{remark}
The definition above is motivated by the argument in a paper of Viana\cite{V}, in which the condition as above for
 $\beta=1/2$ is considered.
\end{remark}
Below we give simple consequences
of the no flat contact condition. 
For  $F\in \U$ and $z\in M$, we define 
\begin{equation}\label{def-L}
L(z;F):=\log\left(\min_{v\in \cone^{u}(z)} |D^{*}F(z,v)|\right)\in \Real\cup\{-\infty\}.
\end{equation}
This function is continuous outside the critical set $\criticalset(F)$ and satisfies
\[
L(z;F)\ge \log d(z,\criticalset(F))-C\subg
\]
from (\ref{eqn-c1}), provided that $\criticalset(F)\neq \emptyset$. Thus we can get the following lemma.
\begin{lemma}\label{lem-nfcc2} Suppose that $F\in \U$ satisfies the no flat contact condition
and  let $n_{0}=n_{0}(F)$ be that in the condition. For any 
$\delta>0$ and $a>0$,    we can choose a positive  number
$h=h(\delta,a;F)$  such that
\[
\int\min\{0, L(z;F)+h\}\; d(\mu\circ F^{-n})(z)\ge -\delta\cdot  |\mu|
\]
for any $\mu\in
\admm([a,\infty))$ and   $n\ge n_0$. 
\end{lemma}
Using the inequality $
\log
\|DF|_{E^{c}(z)}\|\ge L(z;F)-C\subg$, which follows from (\ref{eqn-c1}), together with 
 lemma \ref{lem-nfcc2},
corollary \ref{cor-admc-ind} and corollary~\ref{lem-adm-cpt}, we can obtain
\begin{corollary}\label{cor-bdd} Suppose that  a mapping $F\in \U$ satisfies the no flat contact condition.  Then the central
Lyapunov exponent 
$\chi_{c}(\mu;F)$ considered as a function on the space of \hbox{$F$-invariant} probability measures having 
 admissible lifts is continuous and,  especially, uniformly bounded away from $-\infty$.
\end{corollary}
\noindent
This corollary implies that the irregularity of the central Lyapunov exponent we mentioned  does not take place
under the no flat contact condition.


\subsection{Multiplicity of tangencies between the images of the unstable cones
} \label{sec-mult-trans}
By an iterate of a mapping $F\in \U$, the unstable cones $\cone^u(z)$
at many points $z$ will be brought to one point and some pairs of their images may tangent, that is,
have non-empty intersection.  (Recall that $\cone^u(z)$ does not contain the origin $\mathbf 0$. ) In this subsection, we introduce quantities that measure the multiplicity of such
tangencies  and then formulate a condition, {\em the transversality condition on unstable cones}, for  mappings
in $\U$. 

We introduce analogues of the so-called Pesin subsets. Let
$\bchi=\{\chicl,\chicu,\chiul,\chiuu\}$ be a quadruple of  real numbers that satisfy
\begin{equation}\label{chi-order}
\chicl<\chicu<\chiul<\chiuu.
\end{equation}
Let $\epsilon$ be a  small positive number. For a mapping $F\in \U$,  an
integer
$n>0$ and a real number
$k>0$, we define a closed subset $\Lambda(\chi, \epsilon, k,n ;F)$ of $M$  as the set of all points $z\in
\man$  that satisfy
\begin{align*}
&\chicl (j-i)-\epsilon (n-j)-k
\le \log |D^{*}F^{j-i}(v)|\le
\chicu (j-i)+\epsilon (n-j)+k
\intertext{and}
&\chiul (j-i)-\epsilon (n-j)-k
\le \log D_{*}F^{j-i}(v)
\le
\chiuu (j-i)+\epsilon (n-j)+k
\end{align*}
for any $0\le i<j\le n$ and
$v\in \cone^{u}(F^{i}(z))$. 
Applying the standard argument in the Pesin theory\cite{Pes, PS} to the inverse limit system, we can show
 
\begin{lemma}\label{lem-pesinset}
If $\mu$ is an invariant  probability measure for $F\in \U$ and if
\[
\chicl<\chi_{c}(z;F)<\chicu\quad \mbox{and}\quad \chiul<\chi_{u}(z;F)<\chiuu\quad
\mbox{
$\mu$-a.e. $z$,}
\]
we have $
\lim_{k\to \infty} \liminf_{n\to \infty}\mu(\Lambda(\chi, \epsilon, k,n;F))=1$.
\end{lemma}
The subset $\Lambda(\chi, \epsilon, k,n;F)$ is increasing with respect to $k$ and
$\epsilon$, and satisfies
\begin{align}\label{eqn-shiftlambda0}
&F^{i}(\Lambda(\chi, \epsilon, k, n;F))\subset \Lambda(\chi, \epsilon,k, n-i;F)\quad \mbox{and}\\
\label{eqn-shiftlambda}
&\Lambda(\chi, \epsilon, k, n;F)\subset \Lambda(\chi, \epsilon,k+\epsilon i, n-i;F)\quad\mbox{for 
$0\le i< n$}.
\end{align}

From (\ref{anglerel}), we can and do take a constant
$H\subg$ such that  
\begin{equation}
\angle(DF^{n}(u), DF^{n}(v))
<H\subg
\frac{|D^{*}F^{n}(z,w)|}{D_{*}F^{n}(z,w)} \le H\subg\exp((\chicu-\chiul)n+2k)
\label{hgdef}
\end{equation}
for any $z\in \Lambda(\bchi,\epsilon, k,n;F)$ and 
$u,v,w\in\cone^u(z)$.
 For  $z\in \man$,  let $\cE(z;\chi,\epsilon, k, n;F)$ be  the set of all pairs
$(w,w')$ of points  in $F^{-n}(z)\cap \Lambda(\chi,\epsilon, k,n;F)$ such that
\begin{equation}\label{cond-trans}
\angle(DF^{n}(\E^{u}(w')),
DF^{n}(\E^{u}(w)))\le 5H\subg\exp((\chicu-\chiul)n+2k).
\end{equation}
Note that, if a pair $(w,w')$ of points  in $F^{-n}(z)\cap \Lambda(\chi,\epsilon, k, n;F)$ does 
{\em not} belong to $\cE(z;\chi,\epsilon, k, n;F)$, we have
\begin{equation}\label{eqn-tran-1}
\angle(DF^{n}(u),
DF^{n}(u'))>3H\subg\exp((\chicu-\chiul)n+2k)
\end{equation}
for any $u\in \cone(w)$ and $u'\in \cone(w')$, from  (\ref{hgdef}).

As a  measure for  the multiplicity of tangencies, we consider the number 
\[
\mult(\chi,\epsilon, k, n;F)=\max_{z\in \man}\;\;\max_{w\in F^{-n}(z)\cap \Lambda(\chi,\epsilon,k,n;F)}
\#\{w'\mid (w,w')\in \cE(z;\chi,\epsilon, k,n;F)\}.
\]
This is increasing with respect to $k$ and
$\epsilon$. 

\begin{de} Let $\X=\{\chi(\ell)\}_{\ell=1}^{\ell_0}$ be a   finite collection  of
quadruples  of numbers $\chi(\ell)=\{\chicl(\ell), \chicu(\ell), \chiul(\ell), \chiuu(\ell)\}$
that satisfy
 (\ref{chi-order}).  We say that a mapping $F\in \U$ satisfies {\em the transversality condition on unstable cones
for $\X$} if
\[
\lim_{\epsilon\to +0}\lim_{k\to \infty}\liminf_{n\to\infty}\;\max
\left\{
\frac{\log(\mult(\chi(\ell),\epsilon, k, n;F))}
{n\cdot (\chicl(\ell)+\chiul(\ell)-\chicdif(\ell)-\chiudif(\ell))}\; ; \; 1\le \ell\le \ell_0\right\}<1
\]
where $
\chicdif(\ell)=\chicu(\ell)-\chicl(\ell)$ and $ \chiudif(\ell)=\chiuu(\ell)-\chiul(\ell)$.
\end{de}

\subsection{Measures on the space of mappings}\label{subsec-mom}
In this subsection, we give some additional argument concerning measures on
the space of mappings.  Recall that we denote by
$\tau_\psi: C^{r}(\man,\Real^2)\to C^{r}(\man,\Real^2)$ the translation by  $\psi\in
C^{r}(\man,\Real^2)$, that is, $\tau_\psi(\varphi)=\varphi+\psi$. For an integer $s\ge 0$ and a positive number $d>0$, we put
\begin{equation}\label{defD}
\disk^{s}(d)=\{G\in C^{s}(\man,\Real^{2})\mid \|G\|_{C^{s}}\le d\}.
\end{equation}
The following lemma gives measures on $C^{r}(\man,\Real^2)$ with nice properties.
\begin{lemma}\label{lem-mom}
For an  integer $s\ge 3$, there exists a Borel probability measure $\mom_s$ on
$C^{s-3}(\man,\Real^2)$ such that
\begin{itemize}
\item[{\rm (1)}] $\mom_s$ is quasi-invariant along
$C^{s-1}(\man,\Real^2)$ and
\item[{\rm (2)}] there exists a positive constant  $\bound=\bound_{s}(d)$ for any  $d>0$
such that
\begin{equation*}
\frac{1}{2}\le 
\frac{d(\mom_{s}\circ \tau_{\psi}^{-1})}{d\mom_s}
\le
2\quad  \mbox{$\mom_s$-almost everywhere on $\disk^{s-3}(d)$}
\end{equation*}
for any $\psi\in C^{s}(\man,\Real^{2})$ with
$\|\psi\|_{C^{s}}<\bound$.
\end{itemize}
\end{lemma}
We will give the proof of  lemma~\ref{lem-mom} in the appendix at the end of this paper. This is  one because the lemma itself
has nothing to do with  dynamical systems and one because the proof is merely a combination of
some results in probability theory.  

Henceforth, we fix the measures
$\mom_s$ for
$s\ge 3$ in lemma~\ref{lem-mom}. Note that the measure
$\mom_s$ belongs to
${\mathcal Q}^{r}_{s-1}$ when $s\ge r+3$. 
\begin{lemma}\label{lem-mea-Q} Suppose  $s\ge r+3$. If a Borel subset $X$ in
$C^{r}(\man,\man)$ is shy with respect to the measure $\mom_{s+2}$,  then $X$ is timid for the class $\mathcal Q_{s-1}^r$ of measures.
\end{lemma}
\begin{proof} Take an arbitrary  measure $\mon$ in ${\mathcal Q}^{r}_{s-1}$. The measure $\mom_{s+2}$
is supported on the space $C^{s-1}(M,\Real^2)$, along which $\mon$ is quasi-invariant. Hence  the convolution
$\mon*\mom_{s+2}$ is equivalent to
$\mon$.  From the assumption, we have
\begin{align*}
\mon*\mom_{s+2}(\base_{G}^{-1}(X))&=\int \!\mom_{s+2}\circ \tau_{\psi}^{-1}(\base_{G}^{-1}(X))  d\mon(\psi)\\
&=\int
\mom_{s+2}(\base_{G+\psi}^{-1}(X)) d\mon(\psi)=0
\end{align*}
for any $G\in C^{r}(M,\torus)$. Thus $X$ is shy with respect to $\mon$.
\end{proof}

In order to  evaluate subsets in  $C^{r}(\man,\torus)$  with respect to the measures 
$\mom_{s}$, we will use the following lemma:
\begin{lemma}\label{lem-meas-est}
Let $s\ge r+3$ and $d>0$. Suppose that mappings 
$\psi_{i}\in C^{s}(\man,\Real^{2})$ 
 and positive numbers $T_{i}$ for $1\le i\le m$
 satisfy 
\begin{equation}\label{eqn-small-norm}
\sup_{|t_{i}|\le
T_{i}}\left\|\sum_{i=1}^{m}t_{i}\psi_{i}\right\|_{C^{s}}\le \bound_{s}(d)
\end{equation}
where $\bound_{s}(d)$ is that in  lemma \ref{lem-mom}. 
If  a Borel subset $X$ in $C^{r}(M,\torus)$ satisfies, for  some $\beta>0$, that 
\begin{equation}\label{eqn-para-meas}
\leb_{\Real^m}\left(\left\{(t_{i})_{i=1}^{m}\in \prod_{i=1}^{m}[-T_{i},T_{i}]\; \left|\;   
\varphi+\sum_{i=1}^{m}t_{i}\psi_i\in X\right.\right\}\right)<\beta
\prod_{i=1}^{m} 2T_{i}
\end{equation}
for every $\varphi\in X$, then we have 
\[
\mom_s(\base^{-1}_{G}(X)\cap \disk^{s-3}(d))\le 2^{m+1}\beta\cdot \mom_s (\Phi_{G}^{-1}(Y))\le
2^{m+1}\beta
\]
 for any $G\in C^{r}(\man,\torus)$, where
\[
Y=\left\{\left.\psi+\sum_{i=1}^{m}t_{i}\psi_i\; \right|\; \psi\in X,\;\; |t_i|\le  T_{i}/2\right\}.
\]
\end{lemma}
\begin{proof} Put $Z=\base^{-1}_{G}(X)\cap \disk^{s-3}(d)$ and denote by $\one_Z$ the indicator function of $Z$.
Using Fubini theorem and the properties of $\mom_s$, we get
\begin{align*}
\int \leb_{\Real^m}&\left(\left\{ \t\in \Real^m\; \left|\;  |t_{i}|\le \frac{T_{i}}{2},\; \tilde\psi+\sum
t_i
\psi_i \in Z \right.\right\}\right) d\mom_s(\tilde\psi)
 \\
&=\int_{\left\{\t;|t_{i}|<T_{i}/2\right\}}\left( \int \one_{Z}\left(\tilde\psi+\sum
t_i
\psi_i \right) d\mom_s(\tilde\psi)\right) d\leb_{\Real^m}(\t)
 \\
 &=\int_{\left\{\t;|t_{i}|<T_{i}/2\right\} }\!\!\!\mom_s\left(Z-\sum t_i
\psi_i\right) d\leb_{\Real^m}(\t)\ge 2^{-1}\mom_s(Z)\prod_{i=1}^{m} T_{i}.
\end{align*} 
The integrand of the integral on the first line is positive only if 
$\tilde\psi$ belongs to $\Phi_{G}^{-1}(Y)$ and  bounded by $\beta\prod_{i=1}^{m} 2T_{i}$ from the assumption
(\ref{eqn-para-meas}).  Thus we obtain  the lemma.
\end{proof}

\subsection{The plan of the proof of the main theorems}\label{ss-plan}
Now we can describe the plan of  the proof of the main results, theorem~\ref{main_theorem1} and
\ref{main_theorem2},  more concretely by using the terminology introduced
in the preceding subsections.  We split  the proof into two parts. In the former part, which will be carried out 
in sections 4, 5 and~6, we study ergodic properties of partially hyperbolic endomorphisms in $\U$ that satisfy
the no flat contact condition and the transversality condition on unstable cones for some finite collection of
quadruples. The conclusion in this part is the following theorem:  For a finite or countable
collection 
$\X=\{\chi(\ell)\}_{\ell\in L}$ of quadruples $\chi(\ell)=\{\chicl(\ell),\chicu(\ell)),
\chiul(\ell),\chiuu(\ell)\}$ that satisfy the condition (\ref{chi-order}), we denote by $|\X|$ the union
of the open 
 rectangles $(\chicl(\ell),\chicu(\ell))\times (\chiul(\ell),\chiuu(\ell))$ over $\ell\in L$.
\begin{theorem}\label{th-fp} Let $\X$ be a finite collection of quadruples that satisfy  (\ref{chi-order}), 
\begin{equation}\label{chiorder1}
\chicl+\chiul>0,\qquad \chicl<0
\end{equation}
and also
\begin{equation}\label{xinclude}
\{0\}\times[\lambda\subg,\lmax]\subset |\X|\subset (-2\lmax,2\lmax)\times (0,2\lmax).
\end{equation}
Suppose that
a mapping $F$  in $\U$ satisfies 
the no flat contact condition and the transversality condition on unstable cones for
 $\X$. 
Then $F$
admits a finite collection of ergodic physical measures whose union of basins of attraction has total Lebesgue
measure on $M$.  
In addition,  if
an ergodic physical measure $\mu$ for $F$ satisfies either 
$(\chi_{c}(\mu;F),\chi_{u}(\mu;F))\in |\X|$ or\/ \hbox{$\chi_{c}(\mu;F)>0$}, then $\mu$ is
absolutely continuous with respect to the Lebesgue measure~$\leb$. 
\end{theorem}

In the latter part of the proof, which will be  carried out in section \ref{sec-trans} and
\ref{sec-nonflat}, we show that the two conditions assumed on the mapping $F$ in the theorem above hold for almost
all partially hyperbolic endomorphisms in $\U$, provided that we choose the finite collection $\X$ of quadruple
appropriately.  On the one hand, we will prove  the following theorem in section
\ref{sec-trans}:
For a finite collection $\X$ of quadruples that satisfy (\ref{chi-order}), we denote by $\ss_{1}(\X)$ 
the set of  mappings $F\in \U$ that does {\em not} satisfy the transversality condition on
unstable cones for $\X$.
\begin{theorem}\label{th-ge-trans} 
There exists a countable collection $\X=\{\bchi(\ell)\}_{\ell=1}^{\infty}$ of quadruples satisfying
(\ref{chi-order}) and (\ref{chiorder1}) such that 
\begin{itemize}
\item[{\rm (a)}] $|\X|$ contains the subset  $\{(x_{c},x_{u})\in \Real^2 \mid x_{c}+x_{u}>0,
\lambda\subg\le x_{u}\le
\lmax, x_{c}\le 0\}$, 
\item[{\rm (b)}] $|\X|$ is contained in $(-2\lmax,2\lmax)\times
(0,2\lmax)$, and 
\item[{\rm (c)}] the subset $\ss_{1}(\X')$ for any finite sub-collection $\X'\subset \X$  is shy with
respect to the measures $\mom_s$ for  $s\ge r+3$ and is
a meager subset in $\U$ in the sense of Baire's category argument.
\end{itemize}
\end{theorem}
On the other hand,  we will show the following theorem  in section \ref{sec-nonflat}. Let $\ss_{2}$ be the set of
mappings $F\in \U$ that does {\em not} satisfy the no flat contact condition.
\begin{theorem}\label{th-ge-nonflat}   If an integer $s\ge r+3$ satisfies
the condition (\ref{rs}) for some  integer
$3\le \nu \le r-2$,
then the subset  $\ss_{2}$ is shy with respect to the measure $\mom_s$. Moreover, 
$\ss_{2}$ is contained in a closed nowhere dense subset in
$\U$, provided that $r\ge 19$. 

\end{theorem}

It is easy to check that the three theorems above imply the main theorems: Consider  a
countable set of quadruples $\X=\{\bchi(\ell)\}_{\ell=1}^{\infty}$   in theorem
\ref{th-ge-trans} and put
$\X_{m}=\{\bchi(\ell)\}_{\ell=1}^{m}$ for $m>0$.  Theorem
\ref{th-fp} implies that the complement of $(\cup_{m=1}^{\infty}\ss_{1}(\X_{m}))\cup
\ss_{2}$ in $\U$ is contained in
$\mathcal R^{r}$. Thus  the main theorems, theorem~\ref{main_theorem1} and
\ref{main_theorem2}, restricted to $\U$ follow from  theorem
\ref{th-ge-trans}, 
\ref{th-ge-nonflat} and lemma \ref{lem-mea-Q}. As we noted in subsection \ref{cnbl}, this is enough for the
proof of the main theorems.

\section{Hyperbolic physical measures}
\label{sec-hyp}
In this section, we study hyperbolic physical  measures for partially hyperbolic endomorphisms. Throughout this
section,  {\em we consider a mapping
$F$ in
$\U$ that  satisfies  the no flat contact condition}.

\subsection{Physical measures with negative central  exponent}
In this subsection, we study  physical measures whose central Lyapunov exponent
 is negative. 
\begin{lemma}\label{lemma-negative-exp1}
If an  ergodic probability measure $\mu$ with negative central Lyapunov  exponent has an admissible lift, then it
is a physical measure. 
\end{lemma}
\begin{proof}The central Lyapunov exponent of the measure $\mu$ is bounded way from $-\infty$ by
 corollary
\ref{cor-bdd}. From Oceledec's theorem and the assumption that $\mu$ has an admissible lift,  we can find an
admissible curve $\curve$ such that almost all points with respect to the smooth measure on it are 
forward Lyapunov  regular for~$\mu$.  
According to the Pesin theory,  the so-called Pesin's local stable manifold exists for each of such points on
$\curve$. These local stable manifolds are  transversal to $\curve$ and contained in the basin $\basin(\mu)$ of
$\mu$. Further,  the union of them has positive Lebesgue measure from   absolute
continuity of Pesin's local stable manifolds\cite[\S 4]{PS}. Therefore $\mu$ is a physical
measure. 
\end{proof}
From this lemma and lemma \ref{lem-adm2}, we can get 
\begin{corollary}\label{corabsmu} If an $F$-invariant probability measure $\mu$ has an admissible lift, it has at
most countably many ergodic components with negative central Lyapunov exponent, each of which is a physical
measure and absolutely continuous w.r.t. $\mu$.
\end{corollary}
The basin of an ergodic physical measure with negative central Lyapunov exponent may not have
interior even though we ignore  null subsets with respect to the Lebesgue measure~$\leb$. Nevertheless, we have
\begin{lemma}\label{lem-dis-neg} For an ergodic physical measure $\mu$ with negative central Lyapunov
exponent, there is an open subset $U$  with $\mu(U)=1$ such that,  for a Borel finite measure $\nu$ that has an admissible lift, we have 
$
\nu(\basin(\mu))>0$ if and only if $\limsup_{n\to \infty} \nu\circ F^{-n}(U)>0$. 
 Especially, if we assume $\nu$ to be  $F$-invariant,  we have 
$\nu(\basin(\mu))>0$ if and only if $\nu(U)>0$.
\end{lemma}
\begin{proof}
Recall the proof of lemma \ref{lemma-negative-exp1}.  From  absolute continuity
of Pesin's local stable manifolds, 
there exists an open neighborhood $U_{z}$ for $\mu$-almost every point $z$ such that, if an admissible
curve
$\curve:[0,a]\mapsto \man$ with length $a>2$ satisfies  
$\curve([1,a-1])\cap U_{z}\neq
\emptyset$, the inverse image  $\curve^{-1}(\basin(\mu))$ has positive
Lebesgue measure. Let  $U$ be the union
of such neighborhoods~$U_{z}$. Then we have $\mu(U)=1$ obviously.  If $\limsup_{n\to
\infty}
\nu\circ F^{-n}(U)>0$ for a Borel finite measure $\nu$ that has  an admissible lift, we have
$\nu(\basin(\mu))>0$ from  the choice of $U_z$ and corollary~\ref{cor-admc-ind}. 
Conversely, if we have $\nu(\basin(\mu))>0$, it holds
$
\limsup_{n\to
\infty}
\nu\circ F^{-n}(U)\ge \nu(\basin(\mu))\cdot  \mu(U)>0$. 
\end{proof}
\begin{lemma}\label{lem-conv-adl}
Let $\mu_i$, $i=1,2,\dots$, be  a sequence of mutually distinct $F$-invariant Borel probability
measures each of which has an  admissible lift. If $\mu_i$  converges to some 
measure~$\mu_\infty$ as
$i\to
\infty$, we have 
$\chi_{c}(z;F)\ge 0$ for
$\mu_\infty$-almost every
$z\in M$. 
\end{lemma}
\begin{proof} From corollary~\ref{lem-adm-cpt},  $\mu_\infty$ has an
admissible lift. If the conclusion of the lemma were not true, there should be an ergodic physical measure
$\mu'_{\infty}\ll\mu_{\infty}$  with negative central Lyapunov exponent, from corollary
\ref{corabsmu}.   Take the open set $U$  in lemma~\ref{lem-dis-neg} for~$\mu'_{\infty}$.  On the
one hand, 
$\mu'_{\infty}(U)=1$ and hence $\mu_{\infty}(U)>0$. On the other hand, since $\mu_i\neq \mu'_{\infty}$  except
for  one $i$ at most,   we should have $\mu_i(\basin(\mu'_\infty))=0$ and hence 
$\mu_{i}(U)=0$. These contradict the fact that
$\mu_i$ converges  to
$\mu_\infty$.
\end{proof}
From this lemma and corollary \ref{cor-bdd}, it follows
\begin{corollary}\label{cor-neg-fini}
For any negative number $\chi<0$, there
exist at most finitely many ergodic physical measures for $F$ that satisfies $\chi_{c}(\mu;F)\le
\chi$.
\end{corollary}
Finally we show
\begin{lemma}\label{lem-neg-exp2} Let $\nu$ be a Borel finite measure that is absolutely continuous with respect
to  the Lebesgue measure~$\leb$ and  $\mu$  a limit point of the sequence of measures $n^{-1}\sum_{i=0}^{n-1} \nu\circ F^{-i}$,
$n=1,2,\dots$. Then  we have either
\begin{itemize}
\item[{\rm(a)}] $\chi_{c}(z;F)\ge 0$ for 
$\mu$-almost every point $z\in \man$, or
\item[{\rm(b)}]there is an ergodic physical measure $\mu'\ll\mu$ with negative central Lyapunov exponent and
$\nu(\basin(\mu'))>0$.
\end{itemize}
Especially, for a physical measure $\mu$ for $F$, we have either {\rm (a)} or
\begin{itemize}
\item[{\rm (b$'$)}] $\mu$ is ergodic and has negative
central Lyapunov exponent.
\end{itemize}
\end{lemma}
\begin{proof} Suppose that (a) does not
hold. Then, from corollary \ref{corabsmu}, there exists  an ergodic physical measure
$\mu'\ll\mu$ with negative central Lyapunov exponent. Take
the open set $U$  in lemma \ref{lem-dis-neg} for~$\mu'$.  We have
$\mu'(U)=1$ and hence $\mu(U)>0$. Thus
\[
\limsup_{n\to \infty} n^{-1}\sum_{i=0}^{n-1} \nu\circ
F^{-i}(U)\ge  \mu(U)>0.
\]
Though the measure $\nu$ may not have an admissible lift, we can use the approximation in lemma
\ref{est-Leb} to  conclude
$\nu(
\basin(\mu'))>0$ from the property of $U$.  
\end{proof}

\subsection{Physical measures with positive central exponent}

In this subsection, we investigate  physical measures with  positive central
Lyapunov exponent.  We shall prove the following three propositions. 

\begin{proposition}\label{prop-positive-exp1}
Any physical measure $\mu$ with positive central Lyapunov exponent is ergodic and  absolutely
continuous with respect to the Lebesgue measure~$\leb$. Moreover the basin
$\basin(\mu)$ is an open set modulo Lebesgue null subsets.
\end{proposition}
\begin{proposition}\label{prop-positive-exp3}
For any positive number 
$\chi>0$, there exist at most finitely many ergodic physical measures for $F$ that satisfies 
$\chi_{c}(\mu;F)\ge \chi$. 
\end{proposition}
Let $\basin^{+}(F)$ (resp. $\basin^{-}(F)$) be the union of the basins of ergodic
physical measures  with positive (resp. negative) central Lyapunov exponent. 
\begin{proposition}\label{prop-positive-exp2}
Suppose that a Borel probability measure  $\nu$ on $M$  is absolutely continuous  with respect to
the Lebesgue measure~$\leb$ and supported on the complement of 
$\basin^{-}(F)\cup\basin^{+}(F)$.  
If $\nu_\infty$ is a weak
limit point   of the sequence of measures 
$n^{-1}\sum_{j=0}^{n-1}\nu\circ F^{-j}$, $n=1,2,\dots$, then we have $\chi_{c}(z;F)= 0$ for $\nu_\infty$-almost
every point
$z$.
\end{proposition}

We derive the propositions above from the following single proposition: 
Let $X(i)$, $i=1,2,\cdots$, be Borel subsets in $M$ with positive Lebesgue measure.
We denote by
$\leb_{X(i)}$ the normalization of the restriction  of the Lebesgue
measure $\leb$ to
$X(i)$. For each $i\ge 1$, let $\mu_{i,\infty}$ be a weak limit point of the sequence 
$n^{-1}\sum_{j=0}^{n-1}\leb_{X(i)}\circ F^{-j}$, $n=1,2,\cdots$.  Assume that the
sequence 
$\mu_{i,\infty}$ converges weakly to some measure~$\mu_{\infty}$ as $i\to \infty$. Also assume that 
$\chi_{c}(\mu_\infty;F)>0$ and  that
$\chi_{c}(z;F)\ge 0$,
$\mu_{\infty}$-a.e. 
\begin{proposition} \label{prop-nonhyp}
In the situation as   above, there exist an ergodic physical measure
$\nu_{i,\infty}$  and an open disk $D_i$  in $M$ for  sufficiently
large
$i$ such that
\begin{itemize}
\item[{\rm (a)}] $\nu_{i,\infty}\ll \mu_{i,\infty}$ and  $\nu_{i,\infty}\ll \leb$,
\item[{\rm (b)}] $\chi_{c}(\nu_{i,\infty};F)>0$, 
\item[{\rm (c)}] the radius of $D_i$ is positive and independent of $i$, 
\item[{\rm (d)}] 
$\nu_{i,\infty}(D_{i})>0$ and 
$D_i\subset  \basin(\nu_{i,\infty})$ modulo Lebesgue null subsets.
\end{itemize}
\end{proposition}
Below we prove  proposition \ref{prop-positive-exp1}, \ref{prop-positive-exp3} and
\ref{prop-positive-exp2} using  proposition \ref{prop-nonhyp}.
\begin{proof}[Proof of proposition~\ref{prop-positive-exp1}]
Let $\mu$ be a physical measure such that  $\chi_{c}(\mu;F)> 0$.  From lemma
\ref{lem-neg-exp2}, we have $\chi_{c}(z;F)\ge 0$ for $\mu$-almost
every point $z$.
Apply proposition~\ref{prop-nonhyp} to the situation where $X(i)=\basin(\mu)$ and $\mu_{i,\infty}=\mu_\infty=\mu$
for all
$i\ge 1$. And  let $\nu_{i,\infty}$ and $D_i$ be those in the corresponding
conclusion, which we can
 assume to be
 independent of $i$.  Consider the open set
$V=\bigcup_{n=0}^{\infty}F^{-n}(D_i)$. Then $\basin(\nu_{i,\infty})=V$ modulo
Lebesgue null subsets. Since $\nu_{i,\infty}(V)\ge \nu_{i,\infty}( D_i)>0$ and since
$\nu_{i,\infty}\ll\mu$, we have
$\mu(V)>0$.  Hence
\[
\leb_{\basin(\mu)}(\basin(\nu_{i,\infty}))= \lim_{n\to \infty}n^{-1}\sum_{i=0}^{n-1}
\leb_{\basin(\mu)}\circ F^{-i} (\basin(\nu_{i,\infty}))\ge \mu(V)>0.
\] 
This implies $\mu=\nu_{i,\infty}$. We have proved proposition \ref{prop-positive-exp1}.
\end{proof}
\begin{proof}[Proof of proposition~\ref{prop-positive-exp3}]
Suppose that there exist infinitely many ergodic
physical measures
$\mu_{i}$, $i=1,2,\cdots$, that satisfy $\chi_{c}(\mu_i;F)\ge \chi>0$. By taking a subsequence, we assume
that
$\mu_{i}$ converges to an invariant probability measure $\mu_\infty$ as $i\to \infty$.  Then
we have
$\chi_{c}(\mu_\infty;F)\ge \chi>0$ from corollary \ref{lem-adm-cpt} and corollary
\ref{cor-bdd}.
From lemma \ref{lem-conv-adl}, we have  $\chi_{c}(z;F)\ge 0$ for
$\mu_\infty$-almost every point $z$.
Thus we can apply proposition \ref{prop-nonhyp} to the situation where  $X(i):=\basin(\mu_{i})$ and 
$\mu_{i,\infty}=\mu_i$ for
$i\ge 1$. 
 Since  $\mu_i$'s are ergodic,  the disks $D_{i}$ in the
corresponding conclusion are contained in $\basin(\mu_i)$ modulo Lebesgue
null subsets  and hence  mutually disjoint. But this is impossible because the radii of
the disks
$D_i$ are positive and  independent of $i$.
\end{proof}
\begin{proof}[Proof of proposition~\ref{prop-positive-exp2}]
Let $X=\man\setminus (\basin^{-}(F)\cup\basin^{+}(F))$. 
For the
proof of the proposition, it is enough to show the claim in the case $\leb(X)>0$ and $\nu=\leb_X$.
 Let $\nu_{\infty}$ be a weak limit point of the sequence
$n^{-1}\sum_{j=0}^{n-1}\nu\circ F^{-j}$. 
From lemma~\ref{lem-neg-exp2}, it holds $\chi_{c}(z;F)\ge 0$ for
$\nu_{\infty}$-almost every $z\in M$. 
Thus we have only to prove $\chi_{c}(\nu_{\infty};F)\le 0$. Suppose that we have 
$\chi_{c}(\nu_{\infty};F)>0$. Then we can apply  proposition
\ref{prop-nonhyp}  to the situation  where
$X(i):= X$ for all $i\ge 1$. Let $\nu_{i,\infty}\ll \nu_{\infty}$ and $D_i$ be those in the corresponding
conclusion, which we can assume to be independent of~$i$.  We should have
\[
\nu(\basin(\nu_{i,\infty}))
\ge 
\limsup_{n\to \infty}
\nu(F^{-n}(D_i))
\ge 
\nu_{\infty}(D_i)
>0.
\]
But this contradicts the definition of  $X$ because $\nu_{i,\infty}$ is an ergodic physical measure with positive
central Lyapunov exponent.  
\end{proof}

We proceed to the proof of proposition \ref{prop-nonhyp}.  For positive numbers
$\chi$, $\epsilon$,
$k$ and a positive integer $n$, we define a closed subset 
$\Gamma(\chi,\epsilon, k,n;F)$  as the set of all points $z\in M$ such that,
for  any $0\le
m< n$ and any $v\in \cone^{u}(F^{m}(z))$, 
\begin{itemize}
\item[($\Gamma 1$)] $|D^{*}F^{n-m}(v)|\ge
\exp(\chi (n-m)-k)$ and
\item[($\Gamma 2$)] $|D^{*}F(v)|\ge
\exp(-\epsilon (n-m)-k)$.
\end{itemize}
For the points in $
\Gamma(\chi,\epsilon, k,n;F)$, we have the following estimates on distortion:
\begin{lemma}\label{lemma-exp-dist}For positive numbers $\chi>0$, $0<\epsilon<\chi/10$ and $k>0$, there exists a
positive constant
$\alpha=\alpha(\chi,\epsilon,k)$, which depends only on $\chi$, $\epsilon$ and $k$ besides the objects that we
fixed at the end of  subsection \ref{cnbl},  such that, for any 
$n>0$ and 
$z\in 
\Gamma(\chi,\epsilon,k,n;F)$,  the restriction of $F^n$  to some neighborhood $V$ of $z$
is a diffeomorphism onto the disk\/ $\ball(F^{n}(z),\alpha)$ and we have
\begin{itemize}
\item[(1)] $\|(DF^{n}_{w})^{-1}\|^{-1}> C\subg^{-1}\exp(\chi n-k)$
for $w\in V$, and\\
\item[(2)] $\left|\log|\det DF^{n}_{w}|-\log |\det DF^{n}_{w'}|\right|<1$ for
$w,w'\in V$.
\end{itemize}
\end{lemma}
\begin{proof} 
Fix $v\in \cone^u(z)$ and put $\delta(i)=| D^{*}F^{n-i}(DF^{i}(v))|^{-1}$ for $0\le i<n$. 
Let  $D_n$ be the disk in the tangent space $T_{F^{n}(z)}M$ with
center at the origin and radius~$\alpha$. We define the regions $D_i\subset T_{F^{i}(z)}M$ for $0\le i<n$  so that 
$DF(D_i)$ is the $\delta(i) \alpha$-neighborhood of $D_{i+1}$.  Then we have
\[
{\mathrm{diam}}D_i\le \left\|(DF^{n-i}_{F^{i}(z)})^{-1}\right\|\cdot
\alpha+ \sum_{j=i}^{n-1}\left\|(DF^{j+1-i}_{F^{i}(z)})^{-1}\right\|\cdot 
\delta(j)\cdot \alpha
\]
for $0\le i<n$. Using the relation (\ref{eqn-c3}),  we can check 
\[
\left\|(DF^{j+1-i}_{F^{i}(z)})^{-1}\right\|\cdot 
\delta(j)
\le C\subg\cdot  | D^{*}F(DF^{j}(v))|^{-1}\cdot \delta(i).
\]
Thus, from the conditions ($\Gamma 1$) and  ($\Gamma 2$), it holds
\begin{align*}
{\mathrm{diam}}D_i &\le C\subg  (n-i+1) \cdot  \exp(\epsilon(n-i)+k)  \delta(i) \cdot \alpha\\
&\le C\subg  (n-i+1) \cdot  \exp(-(\chi-\epsilon)(n-i)+2k)\cdot \alpha.
\end{align*}
From  the condition ($\Gamma 2$) and the relation (\ref{eqn-c3}), we have
\[
\| DF_{F^{i}(z)}^{-1}\|^{-1} \ge C\subg^{-1}
\exp(-\epsilon (n-i)-k).
\]
For $v\in D_i$,  we have the estimates
\begin{align*}
&\left\|\exp_{F^{i+1}(z)}^{-1}\circ F\circ \exp_{F^i (z)}(v)-DF_{F^i(z)}(v)\right\|\\
&\qquad\qquad\le C\subg
(\mathrm{diam} D_i)^2\le C\subg n^2 \exp(-(\chi-2\epsilon)(n-i)+3k) \delta(i)\alpha^2
\end{align*}
and
\begin{align*}
&\left\|D(\exp_{F^{i+1}(z)}^{-1}\circ F\circ \exp_{F^i (z)})_{v}-DF_{F^i(z)}\right\|\\
&\qquad\qquad\le C\subg
\mathrm{diam} D_i\le C\subg  n \cdot  \exp(-(\chi-\epsilon)(n-i)+2k)\cdot \alpha.
\end{align*}
Hence, if we take  sufficiently small $\alpha$ depending only on $\chi$, $\epsilon$, $k$ and~$C\subg$,  the
restriction of  $F$ to 
$\exp_{F^{i}(z)}(D_i)$ is a diffeomorphism onto a neighborhood of the subset $\exp_{F^{i+1}(z)}(D_{i+1})$  for $0\le i<n$.
This implies  the first claim of the lemma. We can check, by straightforward
estimates, that the other claims,  (1) and (2),  holds if we take sufficiently small $\alpha$.
\end{proof} 

 From now to the end of this section, we consider the situation 
in
 proposition \ref{prop-nonhyp}. For each $i$, we take a subsequence $n(j;i)\to \infty$ ($j\to
\infty$) such that  the
sequence of measures $n(j;i)^{-1}\sum_{m=0}^{n(j;i)-1}\leb_{X(i)}\circ F^{-m}$ converges to 
$\mu_{i,\infty}$ as $j\to \infty$.
The
following is the key lemma  in  the proof of proposition \ref{prop-nonhyp}.
\begin{lemma}\label{lemma-ample-gamma} There
exist  $\chi>0$, $0<\epsilon<\chi/10$ and  $k>0$ such that
\begin{equation}\label{equ-omega}
\liminf_{j\to
\infty}\frac{1}{n(j;i)}\sum_{m=0}^{n(j;i)-1}\vol_{X(i)}(\Gamma(\chi,\epsilon, k,m;F))>0\quad \mbox{for
sufficiently large $i$}.
\end{equation} 
\end{lemma}
\noindent 
The point of this lemma is that we can take $\chi$, $\epsilon$ and  $k$ uniformly for
sufficiently  large $i$.  Before proving
this lemma, we finish the proof of  proposition
\ref{prop-nonhyp} assuming it.
\begin{proof}[Proof of proposition \ref{prop-nonhyp}] Let the constants $\chi$, $\epsilon$ and $k$ be those in
lemma \ref{lemma-ample-gamma} and 
$\alpha=\alpha(\chi,\epsilon,k,F)$ that in
lemma~\ref{lemma-exp-dist}. We consider a large integer $i$ for which (\ref{equ-omega})  holds. Then we
can take a compact subset $K\subset X(i)$ and a point
$z_0\in
\man$ such that 
\begin{equation}\label{non-tri}
\liminf_{j\to \infty}\frac{1}{n(j;i)}\sum_{m=0}^{n(j;i)-1}(\leb|_{K\cap \Gamma(\chi,\epsilon,
k,m;F)}\circ F^{-m})(\ball(z_0,\alpha/2))>0.
\end{equation}
  Let $\mathcal{D}_{m}$
be the union of the connected components of $F^{-m}(\ball(z_0,\alpha/2))$ that meet
$K\cap \Gamma(\chi,\epsilon,k,m;F)$. Then, on each of the connected components of $\mathcal{D}_{m}$, the mapping
$F^n$ is a diffeomorphism onto
$\ball(z_0,\alpha/2)$ and satisfies the estimates in lemma~\ref{lemma-exp-dist}.  Let
$\nu_i$ be a  limit point of  the sequence 
$\{n(j;i)^{-1}\sum_{m=0}^{n(j;i)-1}(\leb|_{\mathcal{D}_{m}})\circ F^{-m}\}_{j=1}^{\infty}$. Then we have
$\nu_i\le
\leb(X(i))\cdot \mu_{i,\infty}$ and $\nu_i\ll\leb$ and, further,
\[
e^{-1}\cdot
\frac{\nu_i (\ball(z_0,\alpha/2))}{\vol(\ball(z_0,\alpha/2))}
\le \frac{d\nu_i}{d\leb}\le
e\cdot
\frac{\nu_i(\ball(z_0,\alpha/2))}{\vol(\ball(z_0,\alpha/2))}.
\]
We can check that $\nu_i$ is ergodic 
and $\chi_{c}(z;F)>0$ for $\nu_i$-almost every point $z$. (See the remark
below.) Hence there is an ergodic component
$\nu_{i,\infty}$ of
$\mu_{i,\infty}$ such that
$\nu_i\ll 
\nu_{i,\infty}\ll \mu_{i,\infty}$.
The measure  $\nu_{i,\infty}$ and disk $D_i=\ball(z_0,\alpha/2)$ satisfy the conditions in
proposition
\ref{prop-nonhyp}.
\end{proof} 
\begin{remark} Actually, it is not completely simple  to prove that the measure $\nu_i$ in the proof above is ergodic and  that 
$\chi_{c}(z;F)>0$ for $\nu_i$-almost every point $z$.  But there are a few
 standard ways for it. For example, we can argue as follows:  Consider the inverse limit space
of 
$F$
\[
\tilde M_{F}=\{(z_j)_{-\infty<j\le 0}\mid  z_j\in M, \; F(z_j)=z_{j+1}\}
\]
 and the projection  $\pi:\tilde M_{F}\to M$ defined by
$\pi((z_j)_{-\infty<j\le 0})=z_0$.  Let 
$\tilde\mu_{i,\infty}$ be the natural
extension of $\mu_{i,\infty}$.  We can check  that the part $\tilde\nu_i$
of $\tilde\mu_{i,\infty}$ that corresponds to 
$\nu_i$ is supported on a union of local unstable manifolds, each of which is projected onto
the disk $\ball(z_0,\alpha/2)$ by~$\pi$. Further, the conditional measures on those local
unstable manifolds given  by $\tilde\nu_i$ are absolutely continuous with respect to the
smooth measures on them. For any continuous function
$\varphi$ on
$M$, the backward time-average of $\varphi\circ \pi$ is constant on each of the local unstable manifolds.   From
the ergodic theorem, the forward time-average coincides with the backward time-average almost everywhere with
respect to
$\tilde\nu_i\ll \tilde\mu_{i,\infty}$, and is the  pull-back of a function on $M$ by $\pi$. 
Thus it must be constant
$\tilde\nu_i$-almost everywhere.  This implies that
$\nu_i$ is ergodic.   The positivity of the central Lyapunov exponent  is obtained by
considering Lyapunov exponents with respect to  the backward iteration.
\end{remark}

In the remaining part of this subsection, we prove lemma \ref{lemma-ample-gamma}. To begin with, we fix several
constants: Fix
$\chi_{0}>0$ and
$0<s_{0}<1$ such that 
\begin{equation}\label{eqn-schi}
\mu_{\infty}(\{z\in \man\mid \chi_{c}(z)>\chi_{0}\})>s_{0}.
\end{equation}
Also fix a positive number $\epsilon_{0}$ such that $0<\epsilon_{0}<10^{-4}s_{0}\chi_{0}$.
Recall that we are considering a mapping $F\in \U$ that satisfies the no flat contact condition. From
lemma~\ref{lem-nfcc2}, we  can take and fix a large positive constant
$h_{0}>\chi_0$ such that
\[
\int \min\{0,L(z;F)+h_0\} \; d(\mu\circ F^{-n})(z)>
-10^{-1}s_0 \epsilon_{0}
\]
for any measure $\mu$ in $\admm([1,\infty))$ and $n\ge n_0(F)$, where  $L(z;F)$ is the function defined by
(\ref{def-L}) and $n_0(F)$ is the constant in the definition
of the no flat contact condition.  From (\ref{eqn-schi}) and the assumption that $\chi_{c}(z;F)\ge 0$
for $\mu_\infty$-a.e. $z$,  we can take  and fix a  constant $k_0 >h_{0}$ such
that 
\[
\mu_{\infty}(\{z\in M; |D^{*}F^{n}(v)|\ge
\exp(\chi_{0} n-k_0),\forall v\in \cone^{u}(z), \forall n\ge 0\})>s_0
\]
and
\[
\mu_{\infty}(\{z\in M; |D^{*}F^n (v)|\ge
\exp(-\epsilon_{0} n-k_0), \forall v\in \cone^{u}(z), \forall n\ge
0\})>1-\frac{s_0 \epsilon_{0}}{10 h_{0}}.
\]
Finally we fix a positive integer $m_{0}$ that satisfies $\epsilon_{0} m_{0}>10 k_0$.

Next we introduce the  following subsets of $M$:
\begin{align*}
A&=\{z\in M\;;\; |D^{*}F^{m}(v)|>
\exp(\chi_{0} m-2k_0), \forall v\in \cone^{u}(z), 0\le \forall m\le m_0\},\\
B&=\{z\in M\;;\; |D^{*}F^{m}(v)|>
\exp(-\epsilon_{0} m-2k_0), \forall v\in \cone^{u}(z), 0\le\forall m\le m_0\}\supset A,\\
C&= \man\setminus B\quad \mbox{ and} \quad D=\{z\in C\mid L(z;F)\le -h_{0}\}\subset C.
\end{align*}
Note that $A$ and $B$ are open subsets.
From the assumption that the
sequence
$\mu_{i,\infty}$ converges to $\mu_\infty$ as $i\to \infty$,  
we have 
\begin{align}
&\liminf_{j\to\infty}\frac{1}{n(j;i)}\sum_{m=0}^{n(j;i)-1}\leb_{X(i)}(F^{-m}(A))>s_0,\qquad\mbox{ and}
\label{largei1}\\
&\liminf_{j\to\infty}\frac{1}{n(j;i)}\sum_{m=0}^{n(j;i)-1}\leb_{X(i)}(F^{-m}(
B))>1-\frac{s_0\epsilon_0}{10h_0}
\label{largei2}
\end{align}
for sufficiently large $i$.

We fix a large integer $i$ for which  (\ref{largei1}) and (\ref{largei2}) hold. Using lemma
\ref{est-Leb}, we can find a small number
$b_{0}>0$ and a probability measure
$\mu_{0}$  in
$\admm([b_{0},\infty))$ such that
\begin{align}
&|\leb_{X(i)}-\mu_{0}|<s_0/10,\label{eqn-ref-0}\\
&\liminf_{j\to\infty}
\frac{1}{n(j;i)}\sum_{m=0}^{n(j;i)-1}\mu_{0}(F^{-m}(A))>s_{0}, \quad \mbox{and} \label{eqn-ref-1}\\
&\liminf_{j\to\infty} \frac{1}{n(j;i)}\sum_{m=0}^{n(j;i)-1}\mu_0 (F^{-m}(
B))>1-\frac{s_0\epsilon_0}{10h_0}.
\label{eqn-ref-2}
\end{align}
By modifying the measure $\mu_0$ slightly if necessary, we can  and do assume  
\[
\sum_{m=0}^{n_0(F)}\int 
\min\{0, L(F^{m}(z);F)+h_0\}
d\mu_{0}>-\infty
\]
 in addition. 
Then, from corollary \ref{cor-admc-ind} and the choice of 
$h_{0}$, we have also
\begin{equation}
\liminf_{j\to\infty}\frac{1}{n(j;i)}\sum_{m=0}^{n(j;i)-1}\int 
\min\{0, L(F^{m}(z);F)+h_0\}
d\mu_{0}\ge -\frac{s_0 \epsilon_{0}}{10}. \label{eqn-ref-3}
\end{equation}

For $z\in \man$ and integers $m<n$, we denote, by $A_{z}(m,n)$ $B_{z}(m,n)$, $C_{z}(m,n)$ and
$D_{z}(m,n)$,  the set of integers
$m\le q<n$ for which
$F^{q}(z)$ belongs to $A$, $B$, $C$ and $D$ respectively.
Then we have
\begin{lemma}\label{lemma-build-expansion}
A point $z\in M$ belongs to $\Gamma(s_{0}\chi_{0}/40, 4\epsilon_{0}, 6k_0,n;F)$ for $n>0$ if
\begin{itemize}
\item[{\rm (A)}] $\# A_{z}(m,n)\ge s_{0}(n-m)/10$\quad for any $0\le m<n$,\\
\item[{\rm (C)}] $\# C_{z}(m,n)\le \epsilon_{0}(n-m)/h_{0}$\quad  for any $0\le m<n$,\mbox{
and}\\
\item[{\rm (D)}] $\displaystyle \sum_{q\in
D_{z}(m,n)}\min\{0,L(F^{q}(z);F)+h_0\} \ge -\epsilon_{0}(n-m)$ for any $0\le m<n$.
\end{itemize}
\end{lemma}
\begin{proof} 
Consider a point
$z\in M$ and an integer $n$ that satisfy the conditions (A), (C) and (D). Let
$0\le m<n$ and $I=\{m,m+1,\dots, n-1\}$. We call a set of $m_0$~consecutive integers
$\{q,q+1,\dots,q+m_0-1\}$ is an
$A$-interval (resp. a
\hbox{$B$-interval}) if its smallest element $q$ belongs to $ A_{z}(m,n)$ (resp.
$B_{z}(m,n)$). 
If $\{q,q+1,\dots,q+m_0-1\}$ is an
$A$-interval, we have
\begin{equation}\label{aintest}
\sum_{j=0}^{m_0 -1}
\log|D^{*} F(DF^{j}(v))|
\ge \chi_{0} m_0
-2k_{0}> (\chi_{0}-\epsilon_0) m_0
+2k_{0}
\end{equation}
for $v\in \cone^u(F^{q}(z))$, where the second inequality follows from the choice of $m_0$.
Similarly, if $\{q,q+1,\dots,q+m_0-1\}$ is a
$B$-interval, we have
\begin{equation}\label{bintest}
\sum_{j=0}^{m_0 -1}
\log|D^{*} F(DF^{j}(v))|
\ge -\epsilon_{0} m_0
-2k_{0}> -2\epsilon_0 m_0
\end{equation}
for $v\in \cone^u(F^q(z))$.

Take mutually disjoint $A$-intervals  that
cover $A_{z}(m,n)$ and let $I_{A}$ be the union of them. Then
take mutually disjoint \hbox{$B$-intervals} that cover $B_{z}(m,n)\setminus I_{A}$ and let $I_{B}$ be the union of
them. 
  We can
and do take the $B$-intervals in $I_B$ so that their smallest elements are {\em not} contained  in~$I_{A}$. 
Note that $I_{A}$ and $I_{B}$ are not necessarily contained in $I$.

Consider an arbitrary vector $v\in \cone^u(F^{m}(z))$. Then  $DF^{q-m}(v)$ belongs to $\cone^u(F^{q}(z))$ for
$q\ge m$.  From (\ref{aintest}) and the fact that all the $A$-intervals in $I_A$ but one is
contained in $I$, we have
\begin{align*}
\sum_{q\in I_{A}\cap I}\log|D^{*} F(DF^{q-m}(v))|
&\ge (\chi_{0}-\epsilon_0) \#(I_{A}\cap I)
+2k_{0}((\#I_{A}/m_0)-1)-2k_0.
\end{align*}
Each  $A$-interval in $I_A$ meets at most one 
$B$-interval in $I_B$. Thus the number of \hbox{$B$-intervals} in $I_B$ whose intersection with $I\setminus I_A$
has  cardinality less than
$m_0$ is at most
$ (\# I_{A}/m_0)+1$. From this and (\ref{bintest}), we obtain
\[
\sum_{q\in I_{B}\cap (I\setminus I_{A})}\log|D^{*} F(DF^{q-m}(v))|
\ge 
-2\epsilon_{0} \# (I_{B}\cap (I\setminus I_{A}))-2k_0 ((\#I_{A}/m_0)+1).
\]
Since the complement of $I_{B}\cup I_{A}$ in $I$ is contained in
$C_{z}(m,n)$, the condition ($\Gamma 1$) in the definition of the set 
$\Gamma(s_{0}\chi_{0}/40,4\epsilon_{0}, 6k_0,n;F)$ follows from the two inequalities above, the assumptions
(A),(C) and (D) and the choice of $\epsilon_0$. If
$m$ belongs to
$B_{z}(m,n)$, the condition ($\Gamma 2$) holds obviously. Otherwise, 
the condition ($\Gamma 2$) follows from (D) because we have
$\epsilon_0(n-m)/h_0\ge \#C_{z}(m,n)\ge 1 $ in that case from  (C).  
\end{proof}
In order to prove lemma~\ref{lemma-ample-gamma}, we  see how often the assumptions (A), (C) and (D) in the
lemma above hold. For this purpose, we prepare the following  elementary lemma, which we shall use again in
section \ref{sec-neu}.
\begin{lemma}\label{lemma-refine}
Let $\mu$ be a  measure on a measurable space $X$ and  $\psi_{m}$, $m=0,1,\dots$,  a sequence of
non-negative-valued integrable function on
$X$. For a positive number
$\alpha>0$ and an integer $p\ge 0$, let $Y_p(\alpha)$ be the set of points $y\in X$ such that
\[
\sum_{\ell=q}^{p-1}\psi_{\ell}(y)\ge\alpha(p-q)\quad\mbox{ for some $0\le q<p$.} 
\]
Then it holds, for any $n>0$,
\[
\sum_{m=0}^{n-1}\mu(Y_{m}(\alpha))\le \sum_{m=0}^{n}\mu(Y_{m}(\alpha))\le \alpha^{-1}
\sum_{m=0}^{n-1}\int \psi_{m} d\mu.
\]
\end{lemma}
\begin{proof}
For each point $z\in \man$, we define integers  
\[
n=q_{0}(z)\ge  p_{1}(z)>q_{1}(z)\ge p_{2}(z)>q_{2}(z)\ge \cdots\ge p_{j(z)}>q_{j(z)}\ge 1
\]
in the following inductive  manner:  Suppose that  $q_{j}(z)$ has been  defined. If there exist
integers
$p\le q_{j}(z)$  such that $F^{p}(z)\in Y_{p}(\alpha)$,  let $p_{j+1}(z)$ be the maximum of
such integers  and  $q_{j+1}(z)$ the smallest
integer $q<p_{j+1}(z)$  such that
\begin{equation}\label{eqn-pq}
\sum_{\ell=q}^{p_{j+1}(z)-1}\psi_{\ell}(y)\ge \alpha(p_{j+1}(z)-q).
\end{equation}
Otherwise we put $j(z)=j$ and finish the definition. 
Consider the subsets 
\[
Z_{m}=\{z\in \man \mid q_{j}(z)\le m< p_{j}(z)\mbox{ for some $1\le j\le j(z)$}\}
\]
for $0\le m<n$. Then we have $Y_{n}(\alpha)\subset Z_{n-1}$. From (\ref{eqn-pq}), we obtain
\[
\sum_{m=0}^{n-1}\int \psi_{m} d\mu
\ge \alpha \sum_{m=0}^{n-1}\mu(Z_{m})\ge \alpha\sum_{m=1}^{n}\mu(Y_{m}(\alpha))=\alpha\sum_{m=0}^{n}\mu(Y_{m}(\alpha)). \qed
\]
\renewcommand{\qed}{}
\end{proof}
Now we can complete the proof of lemma \ref{lemma-ample-gamma}.
\begin{proof}[Proof of lemma \ref{lemma-ample-gamma}] For $n\ge 0$, let $\tilde{A}_{n}$, 
$\tilde{C}_{n}$ and $\tilde{D}_{n}$ be the set of points $z\in M$ for which
the condition (A), (C) and (D) does  NOT hold, respectively. First, apply lemma~\ref{lemma-refine} to the
case where $\alpha=1-s_0/10$, $n=n(j;i)$ and $\psi_m$ is the indicator function of the complement of
$F^{-m}(A)$. Then, from (\ref{eqn-ref-1}), we obtain
\begin{align*}
\frac{1}{n(j;i)}\sum_{m=0}^{n(j;i)-1}\mu_{0}(\tilde{A}_{m}(z))&\le
\frac{1}{1-\frac{s_0}{10}}\frac{1}{n(j;i)}\sum_{m=0}^{n(j;i)-1}\mu_{0}(M\setminus  F^{-m}(A))\\
&\le \frac{1-s_0}{1-\frac{s_0}{10}}\le 1-\frac{9}{10}s_0
\end{align*}
for sufficiently large $j$. 
Second, apply lemma~\ref{lemma-refine} to the
case where $\alpha=\epsilon_0/h_0$, $n=n(j;i)$ and $\psi_m$ is the indicator function of the set
$F^{-m}(C)=M\setminus F^{-m}(B)$. Then, from (\ref{eqn-ref-2}), we obtain
\[
\frac{1}{n(j;i)}\sum_{m=0}^{n(j;i)-1}\mu_{0}(\tilde{C}_{m}(z))\le
\frac{h_0}{\epsilon_0}\frac{1}{n(j;i)}\sum_{m=0}^{n(j;i)-1}\mu_{0}(F^{-m}(C))
\le \frac{1}{10}s_0
\]
for sufficiently large $j$. 
Third, apply lemma~\ref{lemma-refine} to the
case where $\alpha=\epsilon_0$, $n=n(j;i)$ and $\psi_m(z)=-\min\{0, L(F^{m}(z);F)+h_0\}$. Then,   from
(\ref{eqn-ref-3}), we obtain
\begin{align*}
\frac{1}{n(j;i)}\sum_{m=0}^{n(j;i)-1}\mu_{0}(\tilde{D}_{m}(z))&\le
\frac{-1}{\epsilon_0}\frac{1}{n(j;i)}\sum_{m=0}^{n(j;i)-1}\int\min\{0, L(F^{m}(z);F)+h_0\} d\mu_{0}(z)\\
&\le \frac{1}{10}s_0
\end{align*}
for sufficiently large $j$. From the three inequalities above and  (\ref{eqn-ref-0}), we conclude
\[
\frac{1}{n(j;i)}\sum_{m=0}^{n(j;i)-1}
\left(\leb_{X(i)}(\tilde{A}_{m}\cup \tilde{C}_{m}\cup \tilde{D}_{m})\right)
\le 1-\frac{6}{10}s_0
\]
for sufficiently large $j$. Since the complement of $\tilde{A}_{m}\cup \tilde{C}_{m}\cup \tilde{D}_{m}$ is
contained in the subset $\Gamma(s_{0}\chi_{0}/40, 4\epsilon_{0}, 6k_0,m;F)$ from
lemma~\ref{lemma-build-expansion}, this implies 
lemma~\ref{lemma-ample-gamma}.
\end{proof}


\section{Some  estimates on distortion}\label{sec-dist}
In this section, we give some basic estimates on distortion of the iterates of  mappings in~$\U$. The estimates
 are straightforward   and may look rather tedious.
But we need to check  that some constants in the estimates can be taken uniformly for the
mappings in~$\U$.  This is important especially in our argument in section \ref{sec-trans}, where
we consider perturbations of mappings in
$\U$.

 Let $\bchi=\{\chicl,\chicu,\chiul,\chiuu\}$ be a quadruple 
  satisfying (\ref{chi-order}) and (\ref{chiorder1}), and 
$\epsilon>0$ a small positive number satisfying
\begin{equation}\label{epsilonchoice}
\epsilon<10^{-3}\min\{\chicl+\chiul,\; |\chicl|, \;\chiuu-\chiul,\; \lambda\subg \}.
\end{equation} 
In the argument below, we will  take several constants that depend only on
$\bchi$ and
$\epsilon$ besides the integer $r$ and the objects that we have already fixed in subsection
\ref{cnbl}.  In order to distinguish such kind of constants, we  denote them by symbols
with subscript
$\epsilon$. Also we will use a generic symbol
$C\sube$ for large positive constants of this kind. The usage of these notations is the same as those introduced
in subsection~\ref{ssremc}. The following lemma is the main ingredient of this section.

\begin{lemma}\label{prop-rho} There exist positive constants $0<\rho\sube<1$, $\kappa\sube>1$ and $\kappa\subg>1$
such that the following claim holds for any $F\in\U$,  $k>0$, $n\ge 1$, $z_0\in\Lambda(\bchi,\epsilon, k,n;F)$ and
$0<\rho\le
\rho_{0}$ where 
\begin{align}\label{eqncondrho}
\rho_{0}&:=\rho\sube e^{-4\epsilon n-2k} \min_{0\le i\le j\le 
n}\;\min_{v\in \cone^{u}(F^{i}(z))}|D^{*}F^{j-i}(v)|\\
&\ge \rho\sube
\exp((\chicl  -5\epsilon) n-3k):\notag
\end{align}
For every mapping $G\in
C^{r}(M,M)$ that satisfies
$d_{C^{1}}(F,G)\le \rho$, we can take a point $z(G)$ and its neighborhood $V_{\rho}(G)\ni z(G)$ in a unique
manner so  that
\begin{itemize}
\item[{\rm (i)}] $z(G)$ depends on $G$ continuously and  $z(F)=z_0$,
\item[{\rm (ii)}] $G^{n}(z(G))\equiv F^{n}(z_0)$ and
\item[{\rm (iii)}] the restriction of $G^n$ to $V_{\rho}(G)$ is a diffeomorphism onto $\ball(F^{n}(z_0),\rho)$.
\end{itemize}
Further it holds 
\begin{itemize}
\item[{\rm (iv)}]
${\mathrm{diam}}(V_{\rho}(G))<\kappa\subg\rho
\exp(-\chicl n+k)$,
\item[{\rm (v)}]
$\ball(z(G),\kappa\subg^{-1}\rho\exp(-\chiuu n-k))\subset V_{\rho}(G)$,
\item[{\rm (vi)}]$V_{\rho}(G)\subset
\Lambda(\chi,\epsilon,k+1, n ;F)$,
\item[{\rm (vii)}] $\angle(DG^{n}(\E^{u}(w)), DF^{n}(\E^{u}(z_0)))\le
\kappa\sube e^{2k}\rho$ for any point $w\in
V_{\rho}(G)$, and
\item[{\rm (viii)}] any admissible curve in $\ball(z_0,\kappa\subg^{-1})$
meets 
$V_{\rho}(F)$ in a single curve.
\end{itemize}
\end{lemma}
\begin{proof} 
First of all, notice that the inequality in (\ref{eqncondrho}) follows from the assumption $z_0\in
\Lambda(\bchi,\epsilon, k,n;F)$. We will give the conditions on the choice of  the constants
$\rho\sube$,
$\kappa\sube$ and
$\kappa\subg$ in the course of the argument below.  For
$0\le i\le n$, we put
$\zeta(i)=F^{i}(z_0)$  and 
\[
\delta(i)=\frac{\rho\cdot \exp(\epsilon (n-i)+k)}{\min_{i\le \ell \le
n}\min_{v\in \cone^{u}(\zeta(i))}|D^{*}F^{\ell-i}(v)|}.
\]
Then it holds
\begin{equation}\label{eqn-deltaest}
\rho<\rho\exp(\epsilon(n-i)+k)\le \delta(i)\le \rho\sube\exp(-3\epsilon n-k)
\quad \mbox{for $0\le i\le n$.}
\end{equation}
Using the relation (\ref{eqn-c3}), we can see
\begin{align}\label{eqn-deltest}
\frac{\delta(j)}
{\delta(i)}&\le \exp(-\epsilon(j-i))\frac{\min_{j\le \ell \le
n}\min_{v\in \cone^{u}(\zeta(i))}|D^{*}F^{\ell-i}(v)|}{\min_{j\le \ell \le
n}\min_{v\in \cone^{u}(\zeta(j))}|D^{*}F^{\ell-j}(v)|}\\
&\le  C\subg \exp(-\epsilon(j-i))\|(DF^{j-i}_{\zeta(i)})^{-1}\|^{-1}\notag
\end{align}
for
$0\le i\le j\le n$, and
\begin{align}\label{eqn-deltest2}
\frac{\delta(i+1)}
{\delta(i)}&= \exp(-\epsilon)\frac{\min\{ 1,\;  \min_{i+1\le \ell \le
n}\min_{v\in \cone^{u}(\zeta(i))}|D^{*}F^{\ell-i}(v)|\}}{\min_{i+1\le \ell \le
n}\min_{v\in \cone^{u}(\zeta(i+1))}|D^{*}F^{\ell-i-1}(v)|}\\
&\ge  C\subg^{-1} \exp(-\epsilon)\|(DF_{\zeta(i)})^{-1}\|^{-1}\notag
\end{align}
for
$0\le i\le n$.

We put 
$D_{n}=\ball(0,\rho)\subset T_{\zeta(n)}M$ and define the region $D_{i}\subset T_{\zeta(i)}M$ for $0\le i<
n$ inductively so that 
$DF_{\zeta(i)}(D_i)$  is the $2\delta(i+1)$-neighborhood of
$D_{i+1}\subset T_{\zeta(i+1)}M$. Put
$B_{i}=\exp_{\zeta(i)}(D_{i})$. Then 
\begin{align}\label{eqn-diamD}
{\mathrm{diam}}B_{i}={\mathrm{diam}}D_i &\le 
2\rho \cdot \|(DF^{n-i}_{\zeta(i)})^{-1}\|+\sum_{j=i+1}^{n}4\delta(j)\cdot \|(DF^{j-i}_{\zeta(i)})^{-1}\|\\
&<C\sube 
\delta(i)\le C\sube \rho\sube \exp(-3\epsilon n-k) \notag
\end{align}
for $0\le i\le n$, where the second inequality follows from (\ref{eqn-deltest}) and the third
from~(\ref{eqn-deltaest}). Since $\zeta(0)=z_0\in \Lambda(\chi,\epsilon, k, n;F)$, we have
\begin{equation}\label{inff}
\|(DF_{\zeta(i)})^{-1}\|^{-1}\ge C\subg^{-1}\exp(-\epsilon (n-i)-k)\quad\mbox{for $0\le i\le n$}
\end{equation}
 by (\ref{eqn-c3}). 
Therefore, if we take the constant $\rho\sube$ sufficiently small, we can obtain
\[
\|DG_w -DF_{\zeta(i)}\|\le d_{C^1}(F,G) +C\subg \cdot {\mathrm{diam}} B_{i}<\|(DF_{\zeta(i)})^{-1}\|^{-1}
\]
and
\begin{align*}
d(G(w),\exp_{\zeta(i+1)}\circ DF_{\zeta(i)}\circ \exp^{-1}_{\zeta(i)}(w))
&\le d_{C^1}(F,G) +C\subg\cdot({\mathrm{diam}}B_{i})^{2}<2 \delta(i+1)
\end{align*}
for $0\le i<n$,  $w\in
\ball(\zeta(i),
\mathrm{diam} B_i)$ and  any mapping $G\in C^r(M,M)$ satisfying  $d_{C^{1}}(F,G)\le \rho\le \rho_0$,
where we used the relation
\[
({\mathrm{diam}}B_{i})^{2}\le C\sube \delta(i)^2\le C\sube \rho\sube \exp(-2\epsilon n)\delta(i+1)
\]
which follows from (\ref{eqn-deltaest}), (\ref{eqn-deltest2}) and (\ref{inff}). These two inequalities
imply that the mapping
$G$ restricted to  $\ball(\zeta(i), \mathrm{diam} B_i)\supset B_{i}$ is a diffeomorphism and maps $B_i$ onto a
neighborhood of 
$B_{i+1}$ for $0\le i<n$.
Put $V_{\rho}(G)=\cap_{i=0}^{n}G^{-i}(B_{i})$. Then the restriction of
$G^{n}$  to 
$V_{\rho}(G)$ is a diffeomorphism onto  
$B_n=\ball(F^{n}(z_0),\rho)$. Let $z(G)$ be the unique point in $V_{\rho}(G)$ that
$G^n$ brings to $F^{n}(z_0)$.  Clearly $z(G)$ and $V_{\rho}(G)$ satisfy the conditions (i), (ii) and (iii).

We show the conditions (iv)-(viii). Using   (\ref{eqn-c2}) and (\ref{eqn-c3}), we can check that
 (iv) and (v) follows from  (vi). We  prove (vi) and (vii). Let $G\in\U$ be a mapping that satisfies
$d_{C^{1}}(F,G)\le \rho\le \rho_0$ and 
$w$ a point in $V_{\rho}(G)$.   We put  $w(i)=G^{i}(w)$ for $0\le i\le n$.  
Consider an integer  $0\le i\le  n$ and tangent vectors $v\in
\cone^{u}(\zeta(i))$ and 
$u\in\cone^{u}(w(i))$. For $0\le m\le n-i$, we have
\begin{align*}
\angle(DG_{w(i)}^{m}(u),DF_{\zeta(i)}^{m}(v))&\le \angle(DF_{\zeta(i)}^{m}(u), DF^{m}_{\zeta(i)}(v))\\
&\quad
+\sum_{j=1}^{m}\angle(DF^{m-j+1}_{\zeta(i+j-1)}(DG_{w(i)}^{j-1}(u)),DF^{m-j}_{\zeta(i+j)}(DG_{w(i)}^{j}(u))).
\end{align*}
\begin{remark} In the expression above, we identified  tangent
vectors with their parallel translations and abused the notation slightly. In fact, 
$DF^{m-j}_{\zeta(i+j)}(DG_{w(i)}^{j}(u))$ should have been written
$DF^{m-j}_{\zeta(i+j)}(\tau(DG_{w(i)}^{j}(u)))$ where
$\tau$ is the parallel translation from
$w(i+j)$ to $\zeta(i+j)$. We continue to use such  identification below.
\end{remark}
Since $w(i+j-1)\in B_{i+j-1}$ and  $DG_{w(i)}^{j-1}(u)\in \cone^u(w(i+j-1))$, the parallel
translation of $DG_{w(i)}^{j-1}(u)$ to $\zeta(i+j-1)$ does not belongs to $\cone^c(\zeta(i+j-1))$, provided that
we take sufficiently small $\rho\sube$. Also  we have  
\[
\angle(DF_{\zeta(i+j-1)}(DG_{w(i)}^{j-1}(u)), DG_{w(i)}^{j}(u))
\le C\subg \left({\mathrm{diam}}(B_{i+j-1})+d_{C^1}(F,G)\right)
\]
for $0\le j\le n-i$. Using these and (\ref{anglerel}) in the inequality above, we obtain 
\begin{align}
&\angle(DG_{w(i)}^{m}(u),DF_{\zeta(i)}^{m}(v))\label{ineangle1}\\
& \quad\le A\subg
\frac{|D^{*}F^{m}(v)|}{D_{*}F^{m}(v)}\angle(u,v)
+C\subg
\sum_{j=1}^{m}
\frac{|D^{*}F^{m-j}(DF^{j}(v))|}{D_{*}F^{m-j}(DF^{j}(v))}
({\mathrm{diam}}B_{i+j-1}+\rho)\notag\\
&\quad \le C\subg\exp(-\lambda\subg m)\angle(u,v) +C\subg\sum_{j=1}^{m-1}\exp(-\lambda\subg(m-j)) 
{(\mathrm{diam}}B_{i+j-1}+\rho).\notag
\end{align}

In order to prove the condition (vii), we consider (\ref{ineangle1}) in the case where $i=0$,
$m=n$ and  $v$ and
$u$ are unit tangent vectors in $\E^{u}(z_0)$ and $\E^{u}(w)$, respectively. In
this case, we have 
\begin{align*}
&\frac{|D^{*}F^{n-j}(DF^{j}(v))|}{D_{*}F^{n-j}(DF^{j}(v))}\cdot ({\mathrm{diam}} B_{j-1}+\rho)
\le \frac{|D^{*}F^{n-j}(DF^{j}(v))|}{D_{*}F^{n-j}(DF^{j}(v))} \cdot C\sube
\delta(j-1)
\\
& \le C\sube \rho \exp(\epsilon(n-j)+k)\max_{j\le \ell\le
n}\frac{|D^{*}F^{n-\ell}(DF^{\ell}(v))|}{D_{*}F^{n-\ell}(DF^{\ell}(v))}
\frac{|D^{*}F(DF^{j-1}(v))|^{-1}}{D_{*}F^{\ell-j}(DF^{j}(v))}\\
& \le C\sube \rho \exp(-(\lambda\subg-2\epsilon)(n-j)+2k)
\end{align*}
for $1\le j\le n$, where we used (\ref{eqn-deltaest}) and 
(\ref{eqn-diamD}) in the first inequality, 
(\ref{eqn-c3}) in the second, and  the assumption 
$z_0\in 
\Lambda(\chi,\epsilon, k, n;F)$ in the third. Likewise,  using the estimate $\angle(v,u)\le C\subg d(z_0,w)\le
C\subg {\mathrm{diam}}B_{0}$, we can show
\begin{align*}
\frac{|D^{*}F^{n}(v)|}{D_{*}F^{n}(v)}\angle(u,v)&\le 
\frac{|D^{*}F^{n}(v)|}{D_{*}F^{n}(v)}\cdot C\subg {\mathrm{diam}}B_{0}\le C\sube \rho
\exp(-(\lambda\subg-2\epsilon)n+2k).
\end{align*}
Putting these inequalities in  (\ref{ineangle1}), we  obtain
the condition (vii).

Next we prove the condition (vi). Consider an integer $0\le i\le n$ and a vector $u\in
\cone^{u}(w(i))$. Since $w(i)$ belongs to $ B_{i}$, there is a vector
$v\in
\cone^{u}(\zeta(i))$ such that 
$\angle(u,v)< C\subg {\mathrm{diam}} B_{i}$. From this, (\ref{eqn-diamD}) and (\ref{ineangle1}), we obtain
\begin{align*}
&|D^{*}G(DG_{w(i)}^{\ell}(v))-D^{*}F(DF_{\zeta(i)}^{\ell}(u))|\\
&\le C\subg(|\det DG_{w(i+\ell)}-\det DF_{\zeta(i+\ell)}|
+|D_{*}G(DG_{w(i)}^{\ell}(v))-D_{*}F(DF_{\zeta(i)}^{\ell}(u))|)\\
&\le C\subg(d_{C^1}(F,G)+{\mathrm{diam}}B_{i+\ell}+\angle(DG_{w(i)}^{\ell}(v),DF_{\zeta(i)}^{\ell}(u)))\\
&\le C\subg \rho\sube\exp(-3\epsilon n-k)
\end{align*}
for $0\le \ell\le n-i-1$.
Thus, using  (\ref{inff}), we can obtain
\begin{equation}\label{eqn-dgdf1}
\log\left|\frac{D^{*}G^{j-i}(v)}{D^{*}F^{j-i}(u)}\right|
<\sum_{\ell=0}^{j-i-1}\log\left|\frac{D^{*}G(DG_{w(i)}^{\ell}(v))}{D^{*}F(DF_{\zeta(i)}^{\ell}(u))}\right| <1\quad
\mbox{for
$0\le i\le j\le n$, }
\end{equation}
provided that we take the constant $\rho\sube$ sufficiently small. 
Likewise,  we can get
\begin{equation*}
\log\left|\frac{D_{*}G^{j-i}(v)}{D_{*}F^{j-i}(u)}\right|<1\quad \mbox{for $0\le i\le j\le  n$.}
\end{equation*}
The condition (vi) follows from these two inequalities and the assumption that  $z_0$ belongs to 
$\Lambda(\chi,\epsilon,n,k;F)$. 

Finally we check the
condition (viii). Let $\curve$ be an admissible curve in
$\ball(z_0,\kappa\subg^{-1})$. From the argument in subsection \ref{ss-admc}, the
curvature of  $F_{*}^{i}\curve$ for $0\le i\le n$ is
 bounded by
some constant
$C\subg$, even though  $F_{*}^{i}\curve$ for $0\le i\le n\subg$ may not be admissible. Thus we can take the 
constant 
$\kappa\subg$ so large that the following holds: the intersection of any arc $\tilde \curve$ in $F_{*}^{i}\curve$ with length less than $ 4\lmax
\kappa\subg^{-1}$ with any ball with diameter not larger than $2\kappa\subg^{-1}$ is a single
sub-arc of $\tilde \curve$ with  length less than 
$4\kappa\subg^{-1}$.   The diameter of $B_i$ is bounded by $2\kappa\subg^{-1}$
provided that we take the constant
$\rho\sube$ sufficiently small.
Thus, by induction on $0\le j\le n$, we can check that  
$\curve_j:=\curve\cap (\bigcap_{\ell=0}^{j}F^{-\ell}(\ball(\zeta(\ell), {\mathrm{diam} B_\ell})))$ consists of a
single arc. We obtain the condition (viii) as the case
$j=n$.
\end{proof}

Note that the claim of lemma \ref{prop-rho} remains true even if we get the constant $\rho\sube$ smaller and
$\kappa\sube$ and $\kappa\subg$ larger.  By letting the constant $\rho\sube$ smaller and
$\kappa\sube$ larger if necessary, we can show the following claim  in addition:

\noindent
{\bf Addendum to lemma \ref{prop-rho}}
Suppose that $F\in \U$, $n\ge 1$ and $k>0$. Then 
there exists a neighborhood $W(z)$ for each point $z\in \Lambda(\chi,\epsilon, k,n;F)$ such that 
\begin{itemize}
\item[{\rm (ix)}] The restriction of $F^{n}$ to  $W(z)$ is a diffeomorphism onto the image. Further, 
if
$W(z)\cap W(w)\neq\emptyset$ for some
$w\in 
\Lambda(\chi,\epsilon, k,n;F)$,
then  $F^{n}$ is injective on the union $W(z)\cup W(w)$.
\item[{\rm (x)}] $\vol(W(z))>\kappa^{-1}\sube \exp(-(\chiuu+\max\{\chicu,0\}+7\epsilon)n-6k)$.
\end{itemize}

\begin{proof} We consider a point $z_0\in \Lambda(\chi,\epsilon, k,n;F)$ and continue to use
the notations in  lemma~\ref{prop-rho} and its proof.  Let $\curve$ be  the curve
 in $V_{\rho_0}(F)$ that  $F^{n}$ maps onto the segment
$\{\zeta(n)+t\cdot \ec(\zeta(n))\mid |t|<\rho_0\}\subset \ball(\zeta(n),\rho_0)$ where $\ec(\cdot)$ is a unit
vector in~$\E^{c}(\cdot)$.  From  backward invariance of the central cones
$\cone^c(\cdot)$, the tangent vectors of $\curve$ is contained in the central cones, provided that 
we take sufficiently small $\rho\sube$.  From (\ref{eqn-dgdf1}) and (\ref{eqn-c3}), the length of 
$F^{i}_{*}\curve$ satisfies
\[
|F_{*}^{i}\curve| 
<C\subg\rho_{0}\|(DF_{\zeta(i)}^{n-i})^{-1}\|<C\subg \rho\sube \exp(-4\epsilon n-2k)
\]
and, for the case $i=0$,  
\begin{align*}
|\curve| 
&>C\subg^{-1}\rho_{0}\|(DF_{\zeta(0)}^{n})^{-1}\|>C\subg^{-1}\min_{0\le i\le j\le n}
\rho\sube e^{-4\epsilon n-2k}
 \|(DF^{n-j}_{\zeta(j)})^{-1}\|\|(DF^{i}_{\zeta(0)})^{-1}\|\\
&\ge C\subg^{-1}
\rho\sube\exp(-\max\{\chicu,0\}n-5\epsilon n-4k).
\end{align*}
Next consider the family of parallel segments 
\[
\curve_{y}(t)=y+t\cdot \eu(z_0),\quad 
|t|<\rho\sube\exp(-(\chiuu+2\epsilon)n-2k) 
\]
parameterized by the points $y\in \curve$, where $\eu(z_0)$ is a unit vector in $\E^{u}(z_0)$. And we define
$W(z_0)$ as the region that this family of segments sweeps.  From the estimate on the length of $\curve$ above, 
 we can see that $W(z_0)$ satisfies the condition (x), provided that we take sufficiently large constant
$\kappa\sube$. Since the mapping
$F$ is uniformly expanding in the unstable directions, we can show
\[
|F^{i}_{*}\curve_{y}|\le C\subg \rho\sube\exp(-(\chiuu+2\epsilon)n-2k)  D_{*}F^{i}(\eu(z_0))<C\subg
\rho\sube \exp(-\epsilon n-k)
\]
for $0\le i\le n$. 
Hence the diameter of $F^{i}(W(z_0))$ is  bounded by 
\[
|F^{i}_{*}\curve|+2\cdot \max_{y\in \curve}|F_{*}^{i}\curve_y|\le C\subg\rho\sube\exp(-\epsilon n-k).
\]
If $W(z_0)\cap W(w)\neq\emptyset$ for some point
$w\in 
\Lambda(\chi,\epsilon, n,k;F)$, the diameter of the image $F^{i}(W(z_0)\cup W(w))$ is bounded by
$4C\subg\rho\sube\exp(-\epsilon n-k)$. On the other hand,  
the distance from $\zeta(i)$ to the critical set $\criticalset(F)$  is not less than $C\subg^{-1}\exp(-\epsilon
n-k)$ from (\ref{inff}). Thus, if we take sufficiently small constant 
$\rho\sube$,  the restrictions of $F$ to  $F^{i}(W(z_0)\cup
W(w))$ for
$0\le i<n$ are diffeomorphisms and  hence (ix) holds. 
\end{proof}

The condition (ix)  implies that,  if two points
$z$ and
$w$ in
$\Lambda(\chi,\epsilon, k,n;F)$ satisfy $F^{n}(z)=F^{n}(w)$, the neighborhoods $W(z)$ and
$W(w)$ are disjoint. Thus we obtain the following corollary from  the condition (x).
\begin{corollary}\label{prop-back-card} For any $F\in \U$, $n\ge 1$, $k>0$ and  $\zeta\in M$, we have 
\[
\#(\Lambda(\chi,\epsilon, k,n;F)\cap F^{-n}(\zeta))\le \kappa\sube
\exp((\chiuu+\max\{\chi_{c}^{+},0\}+7\epsilon)n+6k).
\]
\end{corollary}


\section{Physical measures with neutral central Lyapunov exponent}\label{sec-neu}
In this section, we study physical measures with  nearly neutral central Lyapunov exponent. The
goal is the proof of theorem \ref{th-fp}, which will be carried
out in the last three subsections. 
\subsection{An illustration of the idea of the proof} 
The argument in this section is based on a new idea that relate the transversality
condition on unstable cones to absolute continuity of physical measures with nearly neutral central
Lyapunov exponent. In this first subsection, we illustrate the idea using a simple
example, one because it is quite  new in the study of dynamical systems, as far as the author understands, and one
because the argument in the following subsections is rather involved in spite of the simplicity of the idea.

As a simplified model of partially hyperbolic endomorphism, we consider the  skew product 
 $F:[0,1)\times
\Real\to [0,1)\times
\Real$ defined by 
\[
F(x,y)=(d\cdot x, a_{i}x+b_{i}y+c_{i})\quad \mbox{ on  $[i/d,(i+1)/d)\times
\Real$, $i=0,1,2,\cdots,d-1$,}
\]
 where $d\ge 2$ is an integer and $a_i$, $b_i$ and $c_i$ are real numbers.
And we assume that 
\begin{itemize}
\item $|b_{i}|<d$ for $0\le i<d$, so that $F$ is partially hyperbolic with $\E^{c}=\langle\partial/\partial
y\rangle$,
\item $|b_{i}|>d^{-1}$ for $0\le i<d$, so that $F$ is volume-expanding, and
\item $\sum_{i=0}^{d-1} \log |b_{i}|<0$,  so that most of the orbits are bounded.
\end{itemize}
Put
$\theta=\max_{1\le i\le d}|a_{i}|/(d-|b_i|)$ and $b_{\max}=\max_{1\le i\le d} |b_i|$.  Then
$F$ brings a segment with slope less than $\theta$ in absolute value to a union of 
segments with the same property. Assume in addition that 
\begin{equation}\label{eqn-transsimp}
|a_{i}-a_{i'}|>3\theta \cdot b_{\max}\quad\mbox{for any $i\neq i'$.}
\end{equation}
This is a much simplified   analogue of the transversality condition on unstable cones.
Indeed,  if
$\ell_{\sigma}$ is a  segment in
$[i_{\sigma}/d,(i_{\sigma}+1)/d)\times
\Real$ for $\sigma=1,2$,  and if their slopes are bounded by $\theta$ in
absolute value, then (\ref{eqn-transsimp}) implies that the difference between the slopes of their images
under the mapping $F$  is larger than
$\theta \cdot b_{\max}/d$, provided  $i_{1}\neq i_{2}$.

We prove the existence of  an absolutely continuous invariant measure for $F$ with {\em  negative}
central Lyapunov exponent.  First of all,  observe the following fact:  if  Lebesgue integrable functions
$\psi_{1}$ and
$\psi_{2}$ on
$[0,1]\times \Real$ take constant values on lines with slopes
$k_{1}$ and $k_2$ respectively or, in other words, satisfy $\psi_{i}(x,y)=
\psi_{i}(0,y-k_{i}x)$ for $0\le x\le 1$ and $y\in \Real$,   then we have 
\begin{align*}
(\psi_1,\psi_2)_{L^{2}}&=\int\psi_{1}(x,y)\psi_{2}(x,y)dx dy\\
&= \int \psi_{1}(0,y')\cdot \psi_{2}(0,y'+(k_{1}-k_{2})x) dx dy'
\quad \mbox{ where $y'=y-k_{1}x$}\notag\\
&\le  |k_{1}-k_{2}|^{-1}\|\psi_{1}\|_{L^{1}}\|\psi_{2}\|_{L^{1}}\quad\mbox{provided $k_1\neq
k_2$.}
\end{align*}
Let  $\psi(x,y)$ be an $L^{2}$
function on $[0,1]\times \Real$ and suppose that it is the sum of non-negative functions
$\psi_{j}(y)$, $j=1,2,\cdots,m$, that take constant values on lines with slopes $k_j$ with
$|k_j|<\theta$ respectively. Let 
${\mathcal P}_{F}$ and
${\mathcal P}_{i}$, $0\le i<d$, be the Perron-Frobenius operator associated to 
$F$ and its restriction to 
$[i/d,(i+1)/d]\times\Real$ respectively, so that ${\mathcal P}_{F}=\sum_{i=0}^{d-1} {\mathcal P}_{i}$. By using
the transversality condition (\ref{eqn-transsimp}) and the fact we observed above, we can obtain
\begin{align}
\|{\mathcal P}_{F}\psi\|_{L^2}^{2}&=\sum_{i}\|{\mathcal P}_{i}\psi\|^{2}_{L^{2}}+\sum_{
i\neq i'}({\mathcal P}_{i}\psi, {\mathcal P}_{i'}\psi)_{L^{2}}\label{ine:ly}\\
&\le \frac{1}{d\cdot \min\{|b_{i}|\}}\|\psi\|_{L^{2}}^{2}+\frac{d}{\theta\cdot
b_{\max}}\|\psi\|_{L^{1}}^{2}.\notag
\end{align}
\begin{remark}
We can regard this inequality as an analogue of the so-called Lasota-Yorke inequality.
\end{remark}
Note that the coefficient $1/(d\cdot \min\{|b_{i}|\})$ is smaller than $1$ from the assumption. The
Perron-Frobenius operator ${\mathcal P}_F$ preserves $L^{1}$ norm of  non-negative functions
and not dissipative because of the assumption $\sum_{i=0}^{d-1} \log
|b_{i}|<0$. Since the images
${\mathcal P}_{F}^{n}\psi$ again satisfy the condition that we assumed for~$\psi$, we can apply the
inequality  (\ref{ine:ly}) repeatedly and see that 
${\mathcal P}_{F}^{n}(\psi)$, $n=1,2,\dots$, are uniformly bounded with respect to the $L^{2}$norm. Thus we can 
find a non-trivial fixed point of ${\mathcal P}_{F}$ in
$L^{2}([0,1]\times \Real)$ as a $L^2$-weak limit point of the sequence 
$n^{-1}\sum_{m=0}^{n-1}{\mathcal P}_{F}^{m}(\psi)$, $n=1,2,\dots$.  The measure $\mu$ having this fixed point as
 density
is an absolutely continuous invariant measure for $F$, whose central Lyapunov exponent is 
$d^{-1}\sum_{i=1}^{d}
\log |b_{i}|<0$.

In the argument above,  we used the assumption $\sum_{i=1}^{d}
\log |b_{i}|<0$ only to ensure that the Perron-Frobenius operator ${\mathcal P}$ is not dissipative. So, if we
consider mappings on compact surfaces, the same argument should  be valid in the case where the central
Lyapunov exponent is neutral or even slightly positive.  This is the key idea that we will develop in the
following subsections.
\subsection{Semi-norms on the space of measures}  
For a Borel finite measure
$\mu$ on
$\man$ and 
$0<\delta<1$, we define the function
\[
J_{\delta}\mu:\torus\to \Real, \quad 
J_{\delta}\mu(w):=\frac{\mu(\ball(w,\delta))}{\pi\delta^2}=\frac{1}{\pi \delta^{2}}\int \one_{\delta}(w,z) d\mu(z)
\]
where 
\[
\one_{\delta}:\torus\times \torus\to \Real,\qquad \one_{\delta}(w,z)=
\begin{cases}
1,&\mbox{if $d(w,z)<\delta$;}\\
0,&\mbox{otherwise.}
\end{cases}
\] 
And we put, for Borel finite   measures
$\mu$ and $\nu$ on $M$, 
\[
(\mu,\nu)_{\delta}=(J_{\delta}\mu,
J_{\delta}\nu)_{L^{2}(\vol)},\quad 
\|\mu\|_{\delta}=\sqrt{(\mu,\mu)_{\delta}}=\|J_{\delta}\mu\|_{L^{2}(\vol)}.
\]
Obviously $\|\cdot\|_{\delta}$ is a semi-norm and satisfies
\begin{equation}\label{snbdd}
\|\mu\|_{\delta}\le \frac{|\mu|}{\pi\delta^2}.
\end{equation}
The semi-norm $\|\mu\|_{\delta}$ for a measure $\mu$ is essentially decreasing 
with respect to the auxiliary parameter $\delta$. More precisely, we have  

\begin{lemma}\label{lm-semin-dec}
There is an absolute constant $C_0>1$ such
that 
\begin{equation}\label{normess}
\|\mu\|_{\delta}\le C_0\|\mu\|_{\rho} 
\end{equation}
for any 
$0<\rho\le
\delta<1$ and  any Borel finite measure $\mu$.
\end{lemma}
\begin{proof} There is an absolute constant $C$ such that, for any $0<\rho\le\delta<1$, we
can cover the disk
$\ball(0,\delta)$ in
$\Real^2$ by disks  $\ball(w_i,\rho)$, $1\le i\le [C\delta^2/\rho^2]$, by choosing the
points~$w_i$ appropriately. Using Schwarz inequality, we obtain
\begin{align*}
\|\mu\|_{\delta}^2&=\frac{1}{\pi^2 \delta^{4}} \int \mu(\ball(z,\delta))^{2}d\leb(z)
\le \frac{1}{\pi^2 \delta^{4}}\int
\left(
\sum_{i=1}^{[C\delta^2/\rho^2]}\mu(\ball(z+w_i,\rho))
\right)^2 d\leb(z)\\
&\le 
\frac{1}{\pi^2 \delta^{4}}
\cdot C\frac{\delta^2}{\rho^2}\cdot \sum_{i=1}^{[C\delta^2/\rho^2]}\int
\mu(\ball(z+w_i,\rho))^2 d\leb(z)
\le C^2 \|\mu\|_{\rho}^2
\end{align*}
for any Borel finite measure $\mu$ on $M$.
\end{proof}

 We will make use of  the following properties  of  the semi-norm $\|\cdot\|_{\delta}$.
\begin{lemma} \label{lemgensemi}If we have 
$\liminf_{\delta\to 0}\|\mu\|_{\delta}<\infty$ for a Borel finite measure $\mu$, then the measure $\mu$
is absolutely continuous with respect to  the Lebesgue measure $\leb$ and it holds $\lim_{\delta\to
0}\|\mu\|_{\delta}=\|d\mu/d\vol\|_{L^{2}(\vol)}$.
\end{lemma}
\begin{proof}The assumption implies that there exists a sequence $\delta(i)\to +0$ such
that $J_{\delta(i)}\mu$ is uniformly bounded in $L^{2}(\vol)$. Taking a subsequence, we can  assume that
$J_{\delta(i)}\mu$ converges weakly to some  $\psi\in L^{2}(\vol)$ as $i\to \infty$. Since 
$
(f,\psi)_{L^{2}(\leb)}=\lim_{i\to
\infty}\int f\cdot J_{\delta(i)}\mu \:d\vol=\int fd\mu
$
 for any continuous function $f$ on $M$,
we have
$\mu=\psi\cdot \vol$. Now the last equality is standard. 
\end{proof}
\begin{lemma}\label{lmsemibd}
 If a sequence of Borel finite measures $\mu_i$, $i\ge 1$, converges weakly to some Borel finite measure
$\mu_\infty$, then we have
$\|\mu_{\infty}\|_{\delta}=\lim_{i\to\infty}\|\mu_{i}\|_{\delta}$ for $\delta>0$.
\end{lemma}
\begin{proof} We have $\mu_{\infty}(\partial B(z,\delta))=0$ for Lebesgue almost every point $z$,
because 
\[
\int \mu_{\infty}(\partial B(z,\delta)) d\leb(z)=\int_{d(z,w)=\delta}d\mu_{\infty}(w) d\leb(z)=\int 
\leb(\partial B(w,\delta)) d\mu_{\infty}(w)=0.
\]
This implies that $J_{\delta}\mu_{i}$ converges to  $J_{\delta}\mu_{\infty}$ Lebesgue almost
everywhere as $i\to \infty$. Since the semi-norms $\|J_{\delta}\mu_{i}\|_{\delta}$, $i\ge 1$, are uniformly
bounded from (\ref{snbdd}), the lemma follows from  Lebesgue's convergence theorem. 
\end{proof}
\subsection{Two lemmas on the semi-norm $\|\cdot \|_{\delta}$} \label{ss-lyine}
Let  $\bchi=\{\chicl,\chicu,\chiul,\chiuu\}$ be a  quadruple
  satisfying the conditions (\ref{chi-order}) and (\ref{chiorder1}), and $\epsilon$
 a small positive constant satisfying (\ref{epsilonchoice}).  For simplicity, we put 
\[
\chicdif=\chicu-\chicl,\qquad\chiudif=\chiuu-\chiul.
\]
Let $F$ be a mapping in 
$\U$, $k$ a positive number,
$n$ a positive integer  and  $\mu$ a Borel finite measure on $M$  that is supported  on the
subset 
$\Lambda(\chi,\epsilon,k,n;F)$. The aim of this subsection is to
give two lemmas that estimate $\|\mu\circ F^{-n}\|_{\delta}$. Below we shall use the
notation in section
\ref{sec-dist}. 

Suppose that  the measure $\mu$ is absolutely continuous with respect to the  Lebesgue measure $\leb$ and that the
density 
$d\mu/d\leb$ is square integrable. Then it holds 
\[
\|d(\mu\circ F^{-n})/d\leb\|_{L^2(\leb)}^{2}\le m\cdot
\exp(-(\chicl+\chiul)n+2k)\|d\mu/d\leb\|_{L^2(\leb)}^{2}
\] 
where $m=\max\{\#(F^{-n}(w)\cap \Lambda(\chi,\epsilon,k,n;F))\mid w\in M\;\}$,
because 
\[
|\det DF^n|\ge
\exp((\chicl+\chiul)n-2k)\quad \mbox{ on
$\Lambda(\chi,\epsilon,k,n;F)$.}
\]
 The following lemma is a counterpart  of this simple 
fact for the semi-norm $\|\cdot\|_{\rho}$. Recall the constants $0<\rho\sube<1$ and
$\kappa\sube,\kappa\subg>1$  in lemma
\ref{prop-rho}.
\begin{lemma}\label{lem-dir}
Let $\rho$ be a positive number satisfying
\[
0<\rho<\rho\sube\exp((\chicl -5\epsilon) n-3(k+1))/(10 \kappa\subg^{2})
\]
and put 
\[
\delta=10\kappa\subg \rho
\exp(-\chicl n+k+1).
\]
Suppose that  a measure $\mu$  in $\admm([\delta,\infty))$ is supported
on a Borel subset $X$ in  $\Lambda(\chi,\epsilon,k,n;F)$. Then we have
\begin{equation}\label{triv}
\|\mu\circ F^{-n}\|_{\rho}^{2}
\le
I\subg\cdot m\cdot
\exp((-\chicl-\chiul+\chicdif+\chiudif)n+6k)\|\mu\|_{\delta}^{2}
\end{equation}
for some constant $I\subg>0$, where $m=\max\{\#( F^{-n}(w)\cap \ball(X,\delta))\mid w\in M\}$. 
\end{lemma}
\begin{remark} The  point of the lemma above is that the auxiliary  parameter of the semi-norm on
the right hand side of  (\ref{triv}), that is $\delta$, is larger than that on the left hand side, that is $\rho$.
If  the auxiliary parameter on the right hand side were allowed to be much smaller than that on the left
hand side, the inequality  (\ref{triv}) would hold without the assumption that
$\mu$ has an admissible lift.

\end{remark}
\begin{proof} 
For each point $y\in
\Lambda(\chi,\epsilon,k+1,n;F)$,  there is a unique neighborhood $V(y)$ such that
$F^{n}$ restricted to $V(y)$ is a diffeomorphism onto the disk $\ball(F^{n}(y),\rho)$, according to lemma
\ref{prop-rho}. 
Note that the diameter of
$V(y)$ is smaller than $\delta/10$ from lemma \ref{prop-rho}(iv) and the definition of $\delta$.  
 Let 
$U$ be the union of the neighborhoods $V(y)$ for all $y\in X$.  Then $U$ is contained in
$\ball(X,\delta/10)$ and also in $\Lambda(\chi,\epsilon,k+1,n;F)$ from lemma 
\ref{prop-rho}(vi) because $
X$ is a subset of $\Lambda(\chi,\epsilon,k,n;F)$. From the definition of $U$ and the assumption that $\mu$ is
supported on
$X$, it follows
\[
J_{\rho}(\mu\circ F^{-n})(w)=\frac{1}{\pi \rho^2} \;\mu\circ F^{-n}(\ball(w,\rho))=\frac{1}{\pi\rho^2}\sum_{z\in
F^{-n}(w)\cap U}
\mu(V(z))
\]
for $w\in M$.
Suppose that we have proved
\begin{equation}\label{claim-meas-comp}
\mu(V(z))
\le
C\subg\exp(-(\chicl+\chiul)n+2k)
\left(\frac{\rho}{\delta}\right)^{2}
\mu(\ball(z,\delta))
\end{equation}
for any $z\in \Lambda(\chi,\epsilon,k+1,n;F)$.
Then it follows
\begin{equation}\label{induc-eq}
J_{\rho}(\mu\circ F^{-n})(w)\le
C\subg\exp(-(\chicl+\chiul)n+2k) \sum_{z\in F^{-n}(w)\cap U}
J_{\delta}\mu(z)
\end{equation}
for each $w\in
\man$. 
As we have
\[
|\det DF^n| \le \exp((\chicu+\chiuu)n+2k+2)
\quad \mbox{
on $U\subset \Lambda(\chi,\epsilon,k+1,n;F)$,}
\]
we can obtain the
inequality  (\ref{triv}) from  (\ref{induc-eq}) by
integrating the squares of the both sides and using Schwarz inequality. 
Therefore, in order to prove the lemma, it is enough to show the inequality
(\ref{claim-meas-comp}). Since the both sides of  (\ref{claim-meas-comp}) is
linear with respect to
$\mu$, we may assume without loss of generality that 
$\mu$ has an admissible lift that is supported on a single element of the partition  $\Xi_\badmc$ in
$\badmc([\delta,\infty))$.

Let $\curve:[0,a]\to M$ be an admissible curve with length $a\ge \delta$ and $z$  a point in 
$\Lambda(\chi,\epsilon,k+1,n;F)$. 
Consider a connected component
$I$ of $\curve^{-1}(V(z))$ and let $J$ be the connected component of 
$\curve^{-1}(\ball(z,\delta))\supset \curve^{-1}(V(z))$ that contains $I$. As $\delta<\kappa\subg^{-1}$, 
lemma~\ref{prop-rho}(viii) tells that the interval 
$I$ is the unique connected component of $\curve^{-1}(V(z))$ in $J$. 
For the length of $I$, we have 
\[
\leb_{\Real}(I)=|\curve|_{I}|\le |F^{n}_{*}(\curve|_{I})|\exp(-\chiul n+k+2)\le
C\subg\rho\exp(-\chiul n+k+2)
\]
where the first inequality follows from the fact that $\curve|_{I}$  is an admissible curve in
$V(z)\subset\Lambda(\chi,\epsilon,k+2,n;F)
$ and the second from the fact that $F^{n}_{*}(\curve|_{I})$ is a curve in
$F^{n}(V(z))=\ball(F^{n}(z),\rho)$ whose tangent vectors are contained in the unstable cones $\cone^u$.
For the length of $J$, we have $\leb_{\Real}(J)\ge \delta/2$ because the curve $\curve|_{J}$ meets 
$V(z)\subset \ball(z,\delta/10)$  while the length of $\curve$ is not less than
$\delta$.  These estimates  hold for each connected component of
$\curve^{-1}(V(z))$. Thus we obtain
\[
\frac{\leb_{\Real}(\curve^{-1}(V(z)))}
{\leb_{\Real}(\curve^{-1}(\ball(z,\delta)))}<C\subg\frac{\rho\exp(-\chiul
n+k)}{\delta} <C\subg \frac{\rho^2}{\delta^2}\exp(-(\chicl+\chiul)n+2k),
\]
where we used the definition of $\delta$ in the second inequality.  From this and the definition of
admissible measure, we can conclude (\ref{claim-meas-comp}) for any measure $\mu$ that has an admissible lift
supported on $\{\curve\}\times[0,a]$.  
\end{proof}

The next lemma is a counterpart of the inequality
(\ref{ine:ly}). Recall the definition of $\mult(\chi,\epsilon, k,n;F)$ in subsection \ref{sec-mult-trans}.
\begin{lemma}\label{prop-ly} 
Let $\rho$ and $\delta$ be positive numbers that satisfy
\[
\rho\cdot\exp((-\chicl+\epsilon)n)\le \delta\le \exp((\chicl-2\chiuu-3\epsilon)n).
\]
Suppose that a measure  $\mu$ in  $\admm([\delta,\infty))$ is supported on
$\Lambda(\chi,\epsilon,k,n;F)$. Then we have 
\[
\left\|\mu\circ F^{-n}\right\|_{\rho}^{2}
\le 
\frac{\mult(\chi,\epsilon, k+1,n;F)\|\mu\|_{\rho}^{2}}
{\exp((\chicl+\chiul-\chicdif-\chiudif-2\epsilon) n)}
 +
\frac{\exp((-2\chicu+2\epsilon) n)}{\delta^{2}} |\mu|^{2},
\]
provided that $n$ is larger than some integer
$n_{*}=n_{*}(\bchi,\epsilon,k)$ which depends only on $\bchi$, $\epsilon$ and $k$ besides the
objects that we have fixed at the end of subsection \ref{cnbl}.
\end{lemma}
\begin{proof} In the course of the proof below, we will give some conditions on the choice of $n_{*}=n_{*}(\bchi,\epsilon,k)$. First, we require that $n_*$ is so large that it holds
\[
\exp((\chicl-\chiuu-\epsilon)n_*)<\rho\sube\exp((\chicl -5\epsilon)
n_{*}-3(k+1))/(10\kappa\subg^{2}).
\]
Consider an integer $n\ge n_*$ and put  $\rho_{1}:=\exp((\chicl-\chiuu-\epsilon)n)$. Let $\lattice(\rho_1)$ be the lattice that we defined in subsection \ref{sec-notation}. 

For $w\in
\lattice(\rho_1)$, let
$D_{3}(w,i)$,
$1\le i\le m(w)$,  be the connected components of $F^{-n}(\ball(w,3\rho_1))$
that meet
$\Lambda(\chi,\epsilon,k,n;F)$. By lemma \ref{prop-rho} and 
the choice of $n_*$ above, we can check that
 the restriction of
$F^n$ to $D_{3}(w,i)$ is a diffeomorphism onto
$\ball(w,3\rho_1)$ and  that
$D_{3}(w,i)$ is contained in
$\Lambda(\chi,
\epsilon, k+1, n;F)$. Let $D_{1}(w,i)$ and $D_{2}(w,i)$ be the part of $D_{3}(w,i)$ that $F^n$ maps
onto $\ball(w,\rho_1)$ and $\ball(w,2 \rho_1)$ respectively. For $\sigma=1,2,3$,   let 
$D_{\sigma}(w)$ be the union of
$D_{\sigma}(w,i)$ for $1\le i\le m(w)$.

Since the disks $\ball(w, \rho_1)$ for $w\in \lattice(\rho_1)$ cover the torus $\torus$, we have
\[
\mu\circ F^{-n}\le \sum_{w\in
\lattice(\rho_1)}(\mu\circ F^{-n})|_{\ball(w, \rho_1)}.
\]
The function $J_{\rho}((\mu\circ F^{-n})|_{\ball(w,\rho_1)})$ is supported on the disk
$\ball(w, 2\rho_1)$ as $\rho<\rho_1$ from the assumption on $\rho$. And the intersection
multiplicity of the disks 
$\ball(w, 2\rho_1)$ for
$w\in
\lattice(\rho_1)$ is bounded by $10^2$ at most.
Thus we obtain, by Schwarz inequality,
\begin{align*}
\left\|\mu\circ F^{-n}\right\|_{\rho}^{2}
&\le \int  \left(\sum_{w\in
\lattice(\rho_1)} J_{\rho}((\mu\circ F^{-n})|_{\ball(w, \rho_1)})(z)\right)^{2}d\leb(z)
\\
&\le 10^2\int  \sum_{w\in
\lattice(\rho_1)} \left( J_{\rho}((\mu\circ F^{-n})|_{\ball(w,
\rho_1)})(z)\right)^{2}d\leb(z)
\\
&= 10^2
\sum_{w\in \lattice(\rho_1)}\|(\mu\circ F^{-n})|_{\ball(w,\rho_{1})}\|_{\rho}^{2}.
\end{align*}

Since the intersection multiplicity of the regions $D_2(w)$ for $w\in \lattice(\rho_1)$ is also 
bounded by $10^2$,   we have $\sum_{w\in
\lattice(\rho_1)}\mu|_{D_2(w)}\le 10^2
\mu$ and hence
\begin{align*}
\sum_{w\in \lattice(\rho_1)}\|\mu|_{D_2(w)}\|_{\rho}^{2}&=
\int \sum_{w\in \lattice(\rho_1)}(J_{\rho}(\mu|_{D_2(w)})(z))^2 d\leb(z)\\
&\le \int  (10^{2}\cdot
J_{\rho}\mu(z))^2 d\leb(z)\le 10^4 \|\mu\|_{\rho}^{2}.
\end{align*}
Therefore we can deduce the inequality in the lemma from its localized version:
\begin{align}
\left\|(\mu\circ F^{-n})|_{\ball(w,\rho_{1})}\right\|_{\rho}^{2}
&\le 
\frac{\mult(\chi,\epsilon,k+1, n;F)\|\mu|_{D_{2}(w)}\|_{\rho}^{2}}
{\exp((\chicl+\chiul-\chicdif-\chiudif-\epsilon) n)}\label{claim-local}\\
&\qquad\qquad\qquad   +
\frac{\exp((-2\chicu+\epsilon)n ) }
{\delta^2}(\mu(D_{2}(w)))^{2}\notag
\end{align}
for 
$w\in
\lattice(\rho_1)$, provided that we take the constant $n_*$ so large that $\exp(\epsilon n_*)>10^6$.

Below we fix $w\in \lattice(\rho_1)$ and prove the inequality (\ref{claim-local}). 
From the definition of 
$D_{3}(w,i)$ and the assumption that $\mu$ is supported on $\Lambda(\chi,\epsilon, k,n;F)$,
we have
\begin{equation*}
\left(\mu\circ F^{-n}\right)|_{\ball(w,\rho_1)}=\sum_{i=1}^{m(w)}\left(\mu|_{D_{1}(w,i)}\right)\circ F^{-n}.
\end{equation*}
Hence  the left hand side of
the inequality (\ref{claim-local}) is written in the form
\begin{equation}\label{eqn-sum}
\sum_{1\le i,j\le m(w)}((\mu|_{D_{1}(w,i)})\circ
F^{-n}),(\mu|_{D_{1}(w,j)})\circ F^{-n})_{\rho}.
\end{equation}
For $1\le i\le m(w)$, let  
$z_i$ be the  unique point in $D_{3}(w,i)$ such that $F^{n}(z_{i})=w$, which belongs to
$\Lambda(\chi,\epsilon, k+1,n;F)$.  For
$1\le i, j\le m(w)$, we write
$i\pitchfork j$ if the pair $(z_{i},z_{j})$ does not belong to the subset
$\cE(w;\chi,\epsilon,k+1,n;F)$, that is,
\[
\angle(DF^{n}(\E^{u}(z_i)),
DF^{n}(\E^{u}(z_j)))> 5H\subg\exp((\chicu-\chiul)n+2(k+1)).
\]
(See subsection \ref{sec-mult-trans} for the definition of the set $\cE(\cdot)$.) We split the sum
 (\ref{eqn-sum}) into  two parts according to the condition $i\pitchfork j$
and reduce the inequality (\ref{claim-local}) to the following two
inequalities:
\[
\sum_{i\not\pitchfork j}((\mu|_{D_{1}(w,i)})\circ F^{-n},
(\mu|_{D_{1}(w,j)})\circ F^{-n})_{\rho} 
\le \frac{\mult(\chi,\epsilon, k+1, n;F)\|\mu|_{D_{2}(w)}\|_{\rho}^{2}}
{\exp((\chicl+\chiul-\chicdif-\chiudif-\epsilon) n)}
\]
and
\[
\sum_{i\pitchfork j}((\mu|_{D_{1}(w,i)})\circ F^{-n},
(\mu|_{D_{1}(w,j)})\circ F^{-n})_{\rho} 
\le \frac{\exp((-2\chicu+\epsilon) n) }{\delta^2}(\mu(D_{2}(w)))^{2}.
\]
We  denote the sums on the left hand sides of these two inequalities
by $\sum_{\not\pitchfork}$ and $\sum_{\pitchfork}$ respectively. 

We prove the first
 inequality. By Schwarz inequality, we have
\[
\sum_{\not\pitchfork}
\le \sum_{i\not\pitchfork j}\frac{\|(\mu|_{D_{1}(w,i)})\circ
F^{-n}\|^{2}_{\rho}+
\|(\mu|_{D_{1}(w,j)})\circ F^{-n}\|^{2}_{\rho} }{2}.
\]
Since each term $\|(\mu|_{D_{1}(w,i)})\circ F^{-n}\|_{\rho}$ appears for at most 
$2\cdot \mult(\chi,\epsilon,k+1,n;F)$ times 
  on the right
hand side, this implies
\[
\sum_{\not\pitchfork}
\le
\mult(\chi,\epsilon,k+1,n;F)\sum_{i=1}^{m(w)}
\|(\mu|_{D_{1}(w,i)})\circ
F^{-n}\|^{2}_{\rho}.
\]
Besides, we have $\sum_{i=1}^{m(w)}\|\mu|_{D_{2}(w,i)}\|^{2}_{\rho}\le
\|\mu|_{D_{2}(w)}\|^{2}_{\rho}$. Therefore it is enough to show 
\begin{equation}\label{eqn-local2}
\|(\mu|_{D_{1}(w,i)})\circ F^{-n}\|_{\rho}^{2}\le
\frac{\|\mu|_{D_{2}(w,i)}\|^{2}_{\rho}}{\exp(\chicl+\chiul-\chicdif-\chiudif-\epsilon)}.
\end{equation}
We show this inequality by using lemma \ref{lem-dir}.
Unfortunately, we can not apply lemma \ref{lem-dir} directly  to the measure $\mu|_{D_{2}(w,i)}$ because some part
of its admissible lift may be supported on the part of $\badmc((0,\infty))$ that corresponds to very  short
admissible curves, as a consequence of the restriction.  We argue as follows: Observe that 
$F^n$ brings any
$C^1$curve  with length less than $\delta$  in
$D_{3}(w,i)\subset \Lambda(\bchi, \epsilon, k+1,n;F)$ to a curve with  length  less than~$\rho_1$  from the
assumption on $\delta$ and  (\ref{eqn-c2}), provided that
$n_*$ is larger than some constant which depends only on $\epsilon$, $k$ and the constant $C\subg$ in
(\ref{eqn-c2}).  Suppose that an admissible curve
$\curve$ with length 
$a\ge
\delta$ meets $D_{2}(w,i)$ and that a connected component $I$ of 
$\curve^{-1}(D_{2}(w,i))$ has length  less than~$\delta$. Then  the curve
$\curve|_{I}$  meets the boundary of $D_{2}(w,i)$ and hence  $F_{*}^{n}(\curve|_{I})$ meets
 the boundary of $\ball(w,2\rho_1)$. From the observation above, 
$F_{*}^{n}(\curve|_{I})$ does not meet $\ball(w,\rho_1)$ and hence  $\curve|_{I}$ does not meet
$D_{1}(w,i)$. 
Using this fact, we can construct a
 measure $\tilde{\mu}$ in $\admm([\delta,\infty))$ that satisfies
$\mu|_{D_{1}(w,i)}\le \tilde{\mu}\le \mu|_{D_{2}(w,i)}$ by discarding the part of the admissible lift of
$\mu|_{D_{2}(w,i)}$ that is supported on  $\badmc((0,\delta))$. Note that the
observation above implies also that the
\hbox{$\delta$-neighborhood} of $D_{2}(w,i)$ is contained in $D_{3}(w,i)$, so that
$\max\{\#(F^{-n}(z)\cap \ball(D_{2}(w,i),\delta))\mid z\in \man\}=1$. 
Now we  apply proposition~\ref{lem-dir} to $\tilde\mu$ and $X=D_{2}(w,i)$. Then the corresponding conclusion
and (\ref{normess}) imply  (\ref{eqn-local2}), provided that 
$n_*$ is larger than some constant which depends only on $\epsilon$, $k$, $\rho\sube$, $\kappa\subg$ and $I\subg$.

Next we prove the second inequality. It is enough to show 
\begin{equation}\label{seccl}
\begin{split}
\big((\mu|_{D_{1}(w,i)})\circ F^{-n},&\;
(\mu|_{D_{1}(w,j)})\circ F^{-n}\big)_{\rho} \\
&\qquad \le 
\delta^{-2}\exp((-2\chicu+\epsilon) n)\cdot\mu(D_{2}(w,i))\cdot\mu(D_{2}(w,j))
\end{split}
\end{equation}
 for $1\le i,j\le m(w)$ such that $i\pitchfork j$. The both sides of
this inequality are linear with respect to
$\mu|_{D_{2}(w,i)}$ and $\mu|_{D_{2}(w,j)}$. Hence, without loss of generality,  we
can  assume that $\mu|_{D_{2}(w,i)}$ (resp.\ $\mu|_{D_{2}(w,j)}$) has an admissible 
  lift  supported on a single element $\{\curve_i\}\times [0,a_i]$ (resp.  $\{\curve_j\}\times [0,a_j]$) of the
partition $\Xi_\badmc$ and that the curve $\curve_i$ (resp. $\curve_j$) is a  connected component of the
intersection of an admissible curve with length $\ge\delta$ with
$D_{2}(w,i)$ (resp.\ 
$D_{2}(w,j)$). From  the argument in the proof of the first inequality above,   if the length  of the
curve 
$\curve_i$ (resp.\ $\curve_j$) is less than~$\delta$, it can  not
meet $D_{1}(w,i)$ (resp.\ $D_{1}(w,j)$) and hence the inequality (\ref{seccl})
is trivial.  Thereby, we can  assume also that the lengths of $\curve_i$ and
$\curve_j$, that is, $a_i$ and $a_j$,  are not less than~$\delta$. 

By the definition of admissible measure and that of the semi-norm
$\|\cdot\|_{\rho}$, we have
\begin{align*}
&\frac{\left((\mu|_{D_{1}(w,i)})\circ F^{-n},
(\mu|_{D_{1}(w,j)})\circ F^{-n}\right)_{\rho}}
{\mu(D_{2}(w,i))\cdot \mu(D_{2}(w,j))} \\
 &\quad  \le \frac{C\subg}{a_i a_j (\pi \rho^{2})^{2}}\int_{\torus\times [0,a_i]\times [0,a_j]}
\one_{\rho}(F^{n}\circ \curve_{i}
(t),\;y)\cdot \one_{\rho}(F^{n}\circ \curve_{j}(s),\; y)
 d\vol(y)dt ds\notag\\
 &\quad \le  C\subg\delta^{-2}\rho^{-2}\int_{[0,a_i]\times
[0,a_j]}\one_{2\rho}(F^n\circ \curve_{i}(t),\; F^n\circ \curve_{j}(s))
 dt ds.\notag
\end{align*}
We  estimate the last term by using the assumption $i\pitchfork j$.  From (\ref{hgdef}), it holds
\[
\angle(
  DF^{n}(\E^{u}(\curve_i(t))),DF^{n}(\curve'_i(t)))\le H\subg\exp((\chicu-\chiul)n+2(k+1))
\]
for $t\in [0,a_i]$. From lemma~\ref{prop-rho}(vii), it holds
\begin{align*}
\angle(DF^{n}(\E^{u}(z_{i})), DF^{n}(\E^{u}(\curve_i(t))))&\le 
\kappa\sube e^{2(k+1)} \cdot 2\rho_{1}\\
&\le H\subg\exp((\chicu-\chiul)n+2(k+1))
\end{align*}
$t\in [0,a_i]$, where the second inequality follows from the definition of $\rho_1$
provided that
  $n_*$ is larger than some constant which depends only on $\epsilon$, $\kappa\sube$ and $H\subg$.
Thus we have
\[
\angle(DF^{n}(\E^{u}(z_i)), 
 DF^{n}(\curve'_i(t)))\le 2H\subg\exp((\chicu-\chiul)n+2(k+1))\quad \mbox{for $t\in [0,a_i]$}
\]
and the same estimate with the index $i$ replaced by $j$. 
Therefore   the condition  $i\pitchfork j$ implies that, for any $t\in [0,a_i]$ and $s\in
[0,a_j]$,
\[
\angle(
 DF^{n}(\curve'_i(t)),  DF^{n}(\curve'_j(s)))>H\subg\exp((\chicu-\chiul)n+2(k+1)).
\]
By simple geometric consideration using this fact, we can see that the
part of the curve  $F^{n}_{*}\curve_{i}$ that is within distance $2\rho$ from the  curve
$F^{n}_{*}\curve_{j}$ has length less than 
$C\subg\rho \exp(-(\chicu-\chiul)n-2(k+1))$. Since $\curve_i$ and $\curve_j$ are admissible curves in
$\Lambda(\bchi,
\epsilon, k+1,n;F)$, we obtain
\begin{align*}
\leb_{\Real} \{t\in [0,a_i]\mid d(F^{n}(\curve_i(t)), F_{*}^{n}\curve_j)\le 2\rho\}&\le 
\frac{C\subg\rho  \exp(-(\chicu-\chiul)n-2(k+1))}
{\exp(\chiul
n-(k+1))}
\\
&= C\subg \rho \exp(-\chicu n-(k+1))
\end{align*}
and the same inequality with the indices $i$ and $j$ exchanged. These imply
\[
\int_{[0,a_1]\times [0,a_2]}\one_{2\rho}(F(\curve_{i}(t)),F(\curve_{j}(s)))
 dt ds\le C\subg\rho^2\exp(-2\chicu n- 2(k+1)).
\]
Therefore we can conclude (\ref{seccl})  by taking the constant $n_*$ larger if necessary. 
\end{proof}

\subsection{The proof of theorem \ref{th-fp}: Part I}
We give the proof of theorem \ref{th-fp} in the following three subsection.
From this point to the end of this section, we consider the situation assumed in the theorem: Let
$\X$ be a finite collection of quadruples 
$\bchi(\ell)=\{\chicl(\ell),\chicu(\ell),\chiul(\ell),\chiuu(\ell)\}$,
$1\le
\ell\le
\ell_0$, satisfying (\ref{chi-order}), (\ref{chiorder1}) and (\ref{xinclude});
Let $F$ be  a mapping in
$\U$ that satisfy   the no flat contact condition and the transversality condition on  unstable
cones for~$\X$. 
The aim of this subsection is to  derive the conclusions of  theorem~\ref{th-fp}  from the
following proposition:
\begin{proposition}\label{prop-abs} Under the assumptions as above, the following claim holds: Let $\mu_{i}$,
$i\ge 1$, be
 a sequence of Borel probability measures on $M$. We assume either 
\begin{itemize}
\item[{\rm(A)}] every $\mu_{i}$ is invariant and has an admissible lift, or 
\item[{\rm(B)}] $\mu_{i}=n(i)^{-1}\sum_{j=0}^{n(i)-1}\leb_{X}\circ
F^{-j}$ for  some subsequence
$n(i)\to
\infty$, where $\leb_{X}$ is the normalization of the restriction of  the Lebesgue measure~$\leb$ to some Borel subset
$X\subset M$ with positive Lebesgue measure.
\end{itemize}
Further, we assume that $\mu_i$ converges weakly to a
Borel  probability measure
$\mu_{\infty}$ as $i\to \infty$ and that the pair   of Lyapunov exponents $(\chi_{c}(z;F),\chi_{u}(z;F))$  is
contained in the region
$|\X|$ for $\mu_\infty$-almost every point $z$. 
Then, for sufficiently
large $i$, 
there exists a measure $\nu_{i}\le \mu_{i}$  such that 
\begin{itemize}
\item[{\rm (a)}] $|\nu_{i}|>1/3$ and 
\item[{\rm (b)}] $\nu_i$ is absolutely continuous with respect to the Lebesgue measure~$\leb$
 and  the 
   \hbox{$L^{2}$-norm} of the density $d\nu_i/d\leb$ 
 is bounded by a constant  independent of $i$. 
\end{itemize}
\end{proposition}
\noindent
We  assume proposition~\ref{prop-abs} and prove theorem \ref{th-fp}.
\begin{proof}[Proof of theorem \ref{th-fp}]
First, note that, if an ergodic invariant  measure $\mu$  has an  admissible lift, and if the pair of
Lyapunov exponent  
$(\chi_{c}(\mu;F),\chi_{u}(\mu;F))$ of $\mu$ is contained in $ |\X|$,  then $\mu$ is absolutely continuous with
respect to the Lebesgue measure~$\leb$ and, hence, is a physical measure.   This follows immediately from proposition~\ref{prop-abs} if we set
$\mu_{i}=\mu_{\infty}=\mu$ in the assumption~(A). 

We show that there exist at most finitely many ergodic physical measures. Suppose that there exist
infinitely many  mutually distinct ergodic physical measures
$\mu_{i}$,
$i=1,2,\cdots$. By taking a subsequence, we can assume that 
$\mu_{i}$ converges weakly to some measure $\mu_\infty$ as $i\to \infty$.  We have
$\chi_{c}(\mu_{\infty};F)=0$ from corollary
\ref{cor-neg-fini}, proposition \ref{prop-positive-exp3} and  corollary \ref{cor-bdd}. Moreover, we have
$\chi_{c}(z;F)=0$ for \hbox{$\mu_{\infty}$-almost} every point $z$.  In fact, otherwise, there should be an
ergodic physical measure $\mu'_{\infty}\ll
\mu_{\infty}$ with negative central  Lyapunov exponent from lemma \ref{lem-neg-exp2} and
hence   
 $\mu_{i}=\mu'_{\infty}$ for sufficiently large $i$ from lemma~\ref{lem-dis-neg},  which
contradicts the assumption that $\mu_i$ are mutually distinct. Since $\lambda\subg\le \chi_{u}(z;F)\le \lmax$ for
any point $z\in M$ from the choice of the constants $\lambda\subg$ and $\lmax$, the assumption (\ref{xinclude})
implies that the pair of Lyapunov exponents
$(\chi_{c}(z;F), \chi_{u}(z;F))$ is contained in $|\X|$ for $\mu_\infty$-almost every point~$z$. 
Therefore we can apply proposition
\ref{prop-abs} with assumption (A) to the sequence $\mu_i$ and conclude that there is a measure
$\nu_{i}\le
\mu_{i}$ for sufficiently large~$i$ such that
$|\nu_i|>1/3$ and $\|d\nu_{i}/d\leb\|_{L^{2}(\vol)}<C$ for a constant $C$ that is independent of $i$. 
For these measures $\nu_i$, Schwarz inequality gives
\[
(1/3)^{2}<|\nu_i|^{2}\le \vol(\basin(\mu_{i}))\|d\nu_{i}/d\vol\|_{L^{2}(\vol)}^{2}<
C^{2}\vol(\basin(\mu_{i})).
\]
Obviously this contradicts the fact that the basins $\basin(\mu_{i})$ are mutually disjoint.

Let $\basin^{0}$ be the union of the basins of ergodic physical measures whose central Lyapunov exponent
is neutral. Below we prove that
the  Lebesgue measure of the subset $X:=\man\setminus(\basin^{-}\cup\basin^{0}\cup
\basin^{+})$ is zero. Again the proof is  by contradiction. Suppose that the subset $X$ has positive
Lebesgue measure. Then, by choosing a subsequence $n(i)\to \infty$ appropriately, we can assume that the sequence
of measures
$
\mu_{i}=n(i)^{-1}\sum_{j=0}^{n(i)-1}\vol_{X}\circ F^{-j}
$ 
converges  to some measure $\mu_{\infty}$ as $i\to \infty$.
Note that the measures $\mu_i$ are supported on
$X$, for
$F(X)\subset X$.  
From proposition \ref{prop-positive-exp2}, we have $\chi_{c}(z;F)=0$ for 
$\mu_{\infty}$-almost every point $z$. Thus the assumption (\ref{xinclude})
imply that the pair of Lyapunov exponents
$(\chi_{c}(z;F), \chi_{u}(z;F))$ is contained in $|\X|$ for $\mu_\infty$-almost every point~$z$. Each ergodic
component of  
$\mu_{\infty}$ has an admissible lift from lemma \ref{lem-adm2} and  
hence it is a physical measure with neutral central Lyapunov exponent
from the fact we noted in the beginning. Especially $\mu_\infty$ is supported on~$\basin^0$.
Now apply proposition
\ref{prop-abs} with assumption (B) to the sequence $\mu_{i}$ and then let~$\nu_i$ be those in the
corresponding conclusion.  Since  the
density  $\psi_i:=d\nu_{i}/d\vol$ has uniformly bounded $L^{2}$norm for sufficiently large $i$, 
we can assume that  $\psi_{i}$ converges weakly to some 
$\psi_{\infty}\in L^{2}(\vol)$,
by taking a subsequence of
$n(i)$. Note that $\psi_{\infty}$ is not trivial because
\[
(\psi_{\infty},1)_{L^{2}(\vol)}=\lim_{i\to\infty}(\psi_{i},1)_{L^{2}(\vol)}=\lim_{i\to\infty}|\nu_{i}|\ge 1/3.
\]
On the one hand, we have $\int \psi_{i}d\mu_{\infty}=0$ since $\nu_{i}\le \mu_{i}$ is supported on
$X\subset
\man\setminus\basin^{0}$. On the other hand, 
we should have
\[
\lim_{i\to\infty}\int \psi_{i}d\mu_\infty \ge \lim_{i\to\infty}\int
\psi_{i}\psi_{\infty}d\leb
=\lim_{i\to\infty}(\psi_{i}, \psi_{\infty})_{L^2}=\|\psi_{\infty}\|^{2}_{L^{2}(\vol)}>0
\]
because  $\psi_{\infty}\cdot \leb\le \mu_\infty$. We have arrived at  a contradiction.

We have proved that there exists only finitely many ergodic physical measures for $F$ and that the union of
basins of them has total Lebesgue measure. The last statement of theorem \ref{th-fp} follows from proposition 
\ref{prop-positive-exp1} and the fact we noted in the beginning of this proof.
\end{proof}

\subsection{The proof of theorem \ref{th-fp}: Part II}
In this subsection, we give the proof of  proposition \ref{prop-abs},
assuming a lemma, lemma~\ref{prop-abs2}, whose proof is left to the next  subsection. Let
$\mu_i$ and $\mu_\infty$ be those in  proposition \ref{prop-abs}. We denote 
\[
\chicdif(\ell)=\chicu(\ell)-\chicl(\ell)\quad\mbox{and}\quad \chiudif(\ell)=\chiuu(\ell)-\chiul(\ell)
\quad \mbox{ for $1\le \ell\le \ell_0$}.
\]
To begin with, we take and fix several constants in the following order:
\begin{itemize}
\item[{\rm (K1)}] Take $0<\epsilon<1$ so small that (\ref{epsilonchoice}) hold for all the quadruples $\bchi\in
\X$ and that
\[
\lim_{k\to \infty}\liminf_{n\to\infty}\max_{1\le \ell\le \ell_0}\frac{\log
(\mult(\chi(\ell),\epsilon, k,n;F))}{n\cdot (\chicl(\ell)+\chiul(\ell)-\chicdif(\ell)-\chiudif(\ell)-100\epsilon)}
<1.
\]
This is possible from the  transversality condition on unstable cones for 
$\X$.
\item[{\rm (K2)}]  Take positive constants $\rho\sube$ so small and  $\kappa\sube$ so large  that
lemma
\ref{prop-rho} and lemma~\ref{lem-dir} hold for all the quadruples  $\bchi\in \X$ and $\epsilon$ above.
\item[{\rm (K3)}] Take a positive constant $\eta$ so small that
\[
10\lmax \eta<\epsilon  \quad \mbox{and}\quad 
\eta<10^{-3}\epsilon<10^{-3}.
\]
\item[{\rm (K4)}] Take positive constants
$h_{0}$ and
$m_{0}$ so large that $h_{0}>\lmax>1$, $m_0\ge n\subg$ and 
\[
\int \min\{0,\; L(F^{n}(z);F)+h_0\}\; d\mu(z)>-\frac{\eta}{100}\cdot |\mu|
\]
 for any  
$\mu\in \admm([1,\infty))$ and $n\ge m_{0}$, where $L(\cdot)$ is the function defined in (\ref{def-L}). This is
possible from lemma
\ref{lem-nfcc2}.
\item[{\rm (K5)}] 
Take a positive constant $k_{0}$ such that  $k_{0}>h_0$ and that
\[
\mu_{\infty}\left(\bigcup_{\ell=1}^{\ell_0}\Lambda(\chi(\ell),\epsilon,k_{0}-1,
n;F)\right)>1-\frac{\eta}{200 h_{0}}\quad\mbox{for any $n>0$.}
\]
This is possible from lemma \ref{lem-pesinset}, and the assumption on $\mu_\infty$.
\item[{\rm (K6)}] Take a large positive integer $p_{0}$
such that 
\begin{itemize}
\item[(a)]  
$\displaystyle
\mult(\chi(\ell),\epsilon, k_{0}+2, p_{0} ;F)\le
\exp((\chicl(\ell)+\chiul(\ell)-\chicdif(\ell)-\chiudif(\ell)-100\epsilon)p_0) 
$, 
\item[(b)]   $p_0>n_{*}(\bchi(\ell),\epsilon,k_0+1)$ 
\end{itemize}
for  $1\le \ell\le \ell_0$ where
$n_{*}(\cdot)$ is that in lemma \ref{prop-ly}.
This is possible from the choice of $\epsilon$  and the fact that
$\mult(\chi(\ell),\epsilon, k, p_{0} ;F)$ is increasing with respect to~$k$. 
\end{itemize}
Hereafter we will never change the constants taken in (K1)-(K5). Note that we can choose
the integer $p_0$ arbitrarily large in the condition (K6) above. In some places below,  we shall
put additional conditions that
$p_0$ is  larger than some numbers that depend only on $\X$, $c\subg$, $\lambda\subg$,  $\lmax$,
$\kappa\subg$, 
$\ell_0$ and the constants taken in (K1)-(K5). 

For a point $z\in \man$, we define
\[
\k(z)=\min \left\{ k\in \Z\;\left|\; k\ge k_0\mbox{ and } z\in \bigcup_{1\le \ell\le
\ell_0}\Lambda(\chi(\ell),\epsilon, k, p_{0};F)\right.\right\}\ge k_0
\] and 
$\k(z)=\infty$ if the set $\{\cdot\}$ above is empty. 
Also we define
\[
\kone(z)=
\begin{cases}
0,&\mbox{if $\k(z)= k_{0}$;}\\
1,&\mbox{if $\k(z)> k_{0}$.}
\end{cases}
\]
This is the indicator function of the complement of $\bigcup_{1\le \ell\le
\ell_0}\Lambda(\chi(\ell),\epsilon, k_0,
p_{0};F)$. Let $m$ be a positive integer and write it in the form $m=q(m)\cdot  p_{0}+\rmd(m)$ where
$q(m)=[m/p_0]$, so that 
$0\le
\rmd(m)<p_0$. We define the subset $\phys(m)$ as the set of points 
$z\in \man$ that satisfy
\begin{itemize}
\setlength{\topsep}{3pt}
\setlength{\itemsep}{3pt}
\item[{\rm(R1)}]  $\displaystyle \#\{1\le j\le q\mid  \kone(F^{m-j p_{0}}(z))=1\}  
<\frac{\eta\cdot  \q}{10 h_{0}} $\quad
for 
$1\le \q\le q(m)$;
\item[{\rm(R2)}]  $\displaystyle\sum_{j=1}^{\q} (\k(F^{m-j p_{0}}(z))-k_{0}) <\eta \cdot \q 
p_0$\quad  for
$1\le \q\le q(m)$; and 
\item[{\rm(R3)}]  $\k(z) -k_{0} <\eta m$.
\end{itemize}
The following lemma gives a sufficient condition in order that $\phys(m)$,
$m=1,2,\cdots$, are not very small with respect to a measure $\mu$.    
\begin{lemma}\label{lemma-xn} 
Let $\mu$ be a Borel probability measure $\mu$ on $\man$ and $n$ a positive integer such that $n\ge 10p_0$. Assume
that
\begin{equation}\label{eqn-ass-exp}
\sum_{j=0}^{n-1}\int|L(F^{j}(z);F)|\cdot \kone(F^{j}(z))\; d\mu(z)<\frac{\eta n}{10}
\end{equation}
 and that 
\begin{equation}\label{eqn-ass-exp2}
\sum_{j=0}^{n-1}\int \kone(F^{j}(z))\; d\mu(z)<\frac{\eta n}{100 h_{0}}.
\end{equation}
Then we have
$n^{-1}\sum_{m=0}^{n-1}\mu(\phys(m))\ge 1/2$.
\end{lemma}
\begin{proof} For $0\le m<n$, let $\regs_{1}(m)$, $\regs_{2}(m)$ and $\regs_{3}(m)$ be the sets of points $z$ that
{\em violate} the condition (R1), (R2) and (R3)  respectively. We are going to estimate the measures of these
subsets by using lemma \ref{lemma-refine}.  First we give the estimate on the subset $\regs_{1}(m)$ for $0\le
m<n$.  If $z\in \regs_{1}(m)$, we have
\[
\sum_{j=1}^{\q}\kone(F^{m-j p_0}(z))
\ge \frac{\eta \cdot \q}{10 h_{0}}
\]
for some $1\le \q<q(m)$. Using lemma \ref{lemma-refine} with the assumption (\ref{eqn-ass-exp2}), we obtain
\begin{align*}
\sum_{m=0}^{n-1}\mu(\regs_{1}(m))&= \sum_{\rmd=1}^{p_0}\sum_{j=0}^{[(n-\rmd)/p_{0}]}\mu(\regs_{1}((n-\rmd)- j
p_{0}))
\\
&\le \sum_{\rmd=1}^{p_0}\left(\frac{10 h_{0}}{\eta }\sum_{j=0}^{[(n-\rmd)/p_{0}]}\int\kone\left(F^{(n-\rmd)-j
p_0}(z)\right) d\mu(z)\right)\le  \frac{n}{10}.
\end{align*}

Next we give the estimate on the union $\regs_{2}(m)\cup \regs_{3}(m)$. Let us put 
\[
\psi(z)=(|L(z;F)|+5\lmax)\cdot \kone(z)\quad\mbox{ for $j\ge 1$.}
\]
We claim that
\begin{equation}\label{eqn-kz}
\k(z)-k_{0}\le \sum_{j=0}^{p_0-1}\psi(F^{j}(z))\qquad \mbox{for $z\in M$}.
\end{equation}
For a   point $z$, take the smallest integer $0\le p<p_0$ such that $\k(F^{p}(z))= k_{0}$,
and  set $p=p_0$ if there are no such integers. If $p=0$, the inequality (\ref{eqn-kz})
is trivial. So we assume $p>0$. In the case $0<p<p_0$, we choose an integer
$1\le \ell\le
\ell_0$ so that
$\Lambda(\chi(\ell),\epsilon, k_0, p_0;F )$ contains $F^{p}(z)$. In the case
$p=p_0$, we choose $1\le \ell\le
\ell_0$ arbitrarily.  For $0\le i< i'\le p$ and $v\in
\cone^{u}(F^i(z))$, we have the following obvious estimates:
\begin{align*}
&\sum_{j=i}^{i'-1}L(F^{j}(z);F)\le \log|D^{*}F^{i'-i}(v)|\le \lmax(i'-i),\quad {and}\\
&-\lmax\le -c\subg\le  \log|D_{*}F^{i'-i}(v)|\le \lmax(i'-i).
\end{align*}
Using these estimates and the fact that  $F^p(z)\in \Lambda(\chi(\ell),\epsilon, k_0, p_0;F )$ in the case
$p<p_0$, we can check that
$z$ belongs to $\Lambda(\chi(\ell),\epsilon, k , p_0;F )$ for
\[
k=k_{0}+\left[\sum_{j=0}^{p-1}(|L(F^{j}(z);F)|+3\lmax+\epsilon )\right]+1.
\]
This implies (\ref{eqn-kz}). 

If a point $z$ belongs to 
$\regs_{2}(m)$ or $
\regs_{3}(m)$  for $p_0\le m<n$,  we have, from (\ref{eqn-kz}), 
\[
\sum_{j=m'}^{m-1}\psi(F^j(z))
\ge \eta (m-m')\quad \mbox{for some
$0\le m'<m$}. 
\]
As we took $h_0$ so that $h_0>\lmax$, the assumptions (\ref{eqn-ass-exp}) and
(\ref{eqn-ass-exp2}) imply
\[
\sum_{j=0}^{n-1}\int \psi(F^j(z))d\mu(z) 
\le \frac{\eta n}{5}.
\]
Therefore, by using lemma \ref{lemma-refine}, we can obtain
\[
\sum_{m=p_0}^{n-1}\mu(\regs_{2}(m)\cup \regs_{3}(m))\le \frac{n}{5}.
\]
Note that we have $\sum_{m=0}^{p_0}\mu(\regs_{2}(m)\cup \regs_{3}(m))\le p_0\le n/10$, as we assume $n\ge 10 p_0$ in the lemma. Since  $\phys(m)$ is  the complement of $\regs_{1}(m)\cup \regs_{2}(m)\cup
\regs_{3}(m)$,   we can obtain the lemma from the estimates above. 
\end{proof}
The following lemma is the key step in the proof of proposition \ref{prop-abs}. 
\begin{lemma}\label{prop-abs2}
Let $\mu$ be a Borel finite measure on $M$ and $n$ a non-negative integer. If
$\mu$ has an admissible lift $\tmu$ such that  $\tmu\circ F_{*}^{-i}$ belongs to
$\badmm([\exp(-\eta n),\infty))$
 for  $0\le i<n$,
then we have
\[
\|(\mu|_{\phys(n)})\circ F^{-n}\|_{\rho}<C|\mu|+C\exp(-\epsilon n)\|\mu\|_{\rho\exp(-10\eta
n)}
\]
for $0<\rho\le \exp(-10 \lmax p_0)$, where  $C>0$ is a constant that does not depend on the measure $\mu$ nor
the integer $n$.
\end{lemma}
\begin{remark}Actually, the constant $C>0$ above depends only on $\epsilon$, $p_0$, $c\subg$ and $\lmax$.
\end{remark}
We give the proof of this lemma in the next subsection. Below we assume this lemma and
complete the proof of  proposition
\ref{prop-abs}. 
\begin{proof}[Proof of proposition
\ref{prop-abs}]  First consider the case where the assumption (A) holds. 
From the choice of
$k_{0}$, we have 
\[
\mu_{i}\left(\bigcup_{\ell=1}^{\ell_0}\Lambda(\chi(\ell),\epsilon, k_{0}, p_{0};F)\right)>1-
\frac{\eta}{100 h_{0}}
\]
or, in other words, 
\[
\int \kone(z) d\mu_i  <\frac{\eta}{100 h_{0}}
\]
for sufficiently large $i$,
because 
$\Lambda(\chi(\ell),\epsilon, k_0,p_0 ;F)$ contains an open neighborhood of the compact subset
$\Lambda(\chi(\ell),\epsilon, k_0-1,p_0;F)$. The measure 
 $\mu_i$ belongs to $\admm([1,\infty))$  from corollary \ref{cor-admc-ind}.  Thus, it follows from the choice of
$h_{0}$ that
\[
\int \min\{0, L(z;F)+h_0\} d\mu_{i}(z)>-\frac{\eta}{ 100}.
\] 
Hence
\[
\int |L(z;F)|\cdot\kone(z)  d\mu_{i}(z)<h_0 \cdot \frac{\eta}{ 100 h_0} +\frac{\eta}{100} <\frac{\eta}{10}.
\] 
Now we can apply lemma~\ref{lemma-xn} to the invariant measure $\mu_i$ for sufficiently large $i$,  and obtain
\[
n^{-1}\sum_{j=0}^{n-1}\mu_i(\phys(j))\ge \frac12\qquad\mbox{ for 
  $n\ge 10p_0$.}
\]
We put
\[
\nu_{i,n}=\frac{1}{n}\sum_{j=0}^{n-1}(\mu_{i}|_{\phys(j)})\circ F^{-j}\;\; \le \mu_i\qquad \mbox{for $n\ge1$},
\]
so that $|\nu_{i,n}|\ge 1/2$ for $n\ge 10p_0$. Obviously the measure $\mu_i$ has an admissible lift that satisfies
the assumption of  lemma
\ref{prop-abs2} for any $n\ge 0$. Thus it holds
\[
\|\nu_{i,n}\|_{\rho}\le
\frac{1}{n}\sum_{j=0}^{n-1}\left\|(\mu_{i}|_{\phys(j)})\circ F^{-j}\right\|_{\rho}\le 
C+\frac{C}{n}\sum_{j=0}^{n-1}\exp(-\epsilon j)\|\mu_{i}\|_{\rho\exp(-10\eta  j)}
\]
for $0<\rho\le \exp(-10 \lmax p_0)$. This, together with  (\ref{snbdd}) and the choice of $\eta$, 
implies
$\limsup_{n\to
\infty}\|\nu_{i,n}\|_{\rho}\le C$. 
Let~$\nu_i$ be a weak limit point of the sequence $\nu_{i,n}$, $n=1,2,\dots$. Then it holds \hbox{$\nu_i\le
\mu_i$} and
\hbox{$|\nu_i|\ge 1/2$}. Also we have 
$\|\nu_i\|_{\rho}\le C$ for $0<\rho\le \exp(-10
\lmax p_0)$ from lemma~\ref{lmsemibd}.  From lemma \ref{lemgensemi}, this implies that $\nu_i$
is absolutely continuous with respect to the Lebesgue measure~$\leb$  and the density  satisfies
$\|d\nu_{i}/d\leb\|_{L^{2}(\leb)}\le C$. Thus the measures $\nu_i$ satisfy the conditions
in  proposition
\ref{prop-abs}. 

Next we consider the case where the assumption (B) holds. 
Let $n_0=n_{0}(F)>n\subg$ be that in the definition
of the no flat contact condition. Let $X$ and $\leb_{X}$ be those in the assumption (B). Using lemma
\ref{est-Leb},  we can find a small positive number $b>0$ and  
a probability measure $\omega'\in \admm([b,\infty))$ such that
\begin{itemize}
\item $\left|\leb_{X}-\omega'\right|<10^{-3}\eta /h_{0}$,
\item $\omega'\circ F^{-n_0}$ is
absolute continuous with respect to the Lebesgue measure~$\leb$,  
\item the density of the measure $\omega'\circ F^{-n_0}$, $d(\omega'\circ F^{-n_0})/d\leb$, is square integrable.
\end{itemize}
\begin{remark}
In the third condition above, we do not care how large the $L^2$norm is.
\end{remark}
We put $\omega=\omega'\circ F^{-n_0}$ and  
\[
\mu'_{i}=n(i)^{-1}\sum_{j=0}^{n(i)-1}\omega\circ F^{-j},\qquad \mbox{ for $i=1,2,\dots$.}
\]
Then, for sufficiently large $i$, we have $|\mu_i-\mu'_i|<10^{-3}\eta /h_{0}$ and hence
\[
\mu'_{i}\left(\bigcup_{\ell=1}^{\ell_0}\Lambda(\chi(\ell), \epsilon, k_{0}, p_0;F)\right)>1-\frac{\eta }{100
h_{0}},\quad \mbox{that is,}\quad \int \kone(z) d\mu'_i  <\frac{\eta}{100 h_{0}}
\]
from the choice of $k_0$. 
From corollary \ref{cor-admc-ind},  $\omega\circ F^{-j}$ belongs to $\admm([1,\infty))$
for
sufficiently large $j$. 
Thus we have
\[
\int |L(z;F)|\cdot \kone(z)d\mu'_{i}(z)<h_{0}\cdot \frac{\eta }{100 h_{0}}
+\frac{\eta }{100}<\frac{\eta }{10}
\]
for sufficiently large $i$, from the choice of $h_0$. Now we can apply lemma~\ref{lemma-xn} to  $\mu=\omega$
and $n=n(i)$ in order to obtain
\[
n(i)^{-1}\sum_{m=0}^{n-1}\omega(\phys(m))\ge 1/2
\]
for  sufficiently large $i$.
Let  $\tilde \omega'$ be an admissible lift of $\omega'$ that belongs to $\badmm([b,\infty))$ and put
$\tilde\omega=\tilde \omega'\circ F_{*}^{-n_0}$. Then $\tilde\omega$ is an admissible lift of $\omega$. 
Take  a large positive integer $n_1$ that satisfies 
$
\exp(-\eta\cdot n_1)<b\exp(-c\subg)$.
From lemma \ref{lem-admc-ind}, the measures $\tilde\omega\circ F_{*}^{-i}=\tilde\omega'\circ F_{*}^{-i-n_0}$
for $i\ge 0$ belongs to
$\badmm([\exp(-\eta n),\infty))$, provided $n\ge n_1$.  Thus we can apply lemma~\ref{prop-abs2} to $\omega$, and 
obtain
\[
\|(\omega|_{\phys(n)})\circ F^{-n}\|_{\rho}<C|\omega|+C\exp(-\epsilon n)\|\omega\|_{\rho\exp(-10\eta
n)}
\]
for $0<\rho\le \exp(-10 \lmax p_0)$ and $n\ge n_1$.  We put
\[
\nu'_{i}=\frac{1}{n(i)}\sum_{j=n_1}^{n(i)-1}(\omega|_{\phys(j)})\circ F^{-j}\le \mu'_{i}\quad
i=1,2,\dots.
\]
Then, for sufficiently large $i$,  we have $|\nu'_{i}|\ge (2/5)$ and 
\begin{align*}
\|\nu'_{i}\|_{\rho}&
\le C+\frac{C}{n(i)}\sum_{j=n_1}^{n(i)-1}\exp(-\epsilon j)\|\omega\|_{\rho\exp(-10\eta  j)}
\end{align*}
for $0<\rho\le \exp(-10 \lmax p_0)$. Letting $\rho\to +0$ in the last inequality, we obtain  
\[
\left\|\frac{d\nu'_{i}}{d\vol}\right\|_{L^{2}} 
\le C+\left(\frac{C}{n(i)}\sum_{j=n_1}^{n(i)-1}\exp(-\epsilon
j)\right)\left\|\frac{d\omega}{d\vol}\right\|_{L^{2}}
\]
by lemma~\ref{lemgensemi}.
 Since we have 
$|\mu'_{i}-\mu_{i}|<10^{-2}$ and $\nu'_i\le \mu'_i$, we can find a Borel measure $\nu_{i}$ such that
$\nu_{i}\le \nu'_{i}$, $\nu_{i}\le
\mu_{i}$, $|\nu_{i}|>1/3$ and $\|d\nu_i/d\leb\|_{L^2}\le 2C$ for sufficiently large $i$. The measures 
$\nu_{i}$ satisfy the conditions in proposition~\ref{prop-abs}. 
\end{proof}
\subsection{The proof of theorem \ref{th-fp}: Part III}
In this subsection, we give the proof of  lemma \ref{prop-abs2} and complete the proof of 
theorem \ref{th-fp}. 
Let $n$, $\mu$ and $\tmu$ be those in lemma \ref{prop-abs2}.  
Recall the mapping  $\proj:\tbadmc((0,\infty))\to \man$ and the commutative relation (\ref{cdF}) 
in  subsection 
\ref{ss-admc}. Below we divide the measure $\tmu$
into many parts so that we can evaluate the semi-norms  of their images under the mapping 
$\Pi\circ F_{*}^n$ by the two  inequalities we gave in subsection \ref{ss-lyine}. 

We write the integer $n$  in the form
$n=q(n) p_{0}+\rmd(n)$ where $q(n)=[n/p_0]$, so $0\le d(n)<p_0$. 
For integers $-1\le \q\le \q(n)$, 
we put
\[
\tau(\q)=
\begin{cases}
\q p_{0}+\rmd(n), &\mbox{in the case $0\le \q\le  q(n)$};\\
0 , &\mbox{in the case $\q=-1$,}
\end{cases}
\]
so that $\tau(q(n))=n$, and also put
\[
\delta(\q)=
\begin{cases}
\exp\left(-4\eta (n-\tau(\q))-7\lmax p_0-c\subg\right), &\mbox{in the case $0\le \q\le  q(n)$};\\
\exp\left(-4\eta n-7\lmax p_0\right),&\mbox{in the case $\q=-1$}.
\end{cases}
\]
Take and fix a  number  $0<\rho\le\exp(-10 \lmax p_0)$ arbitrarily and  put
\[
\rho(\q)=
\rho\exp(-10\eta (n-\tau(\q)))\quad \mbox{for  $-1\le \q\le \q(n)$}.
\]
We put $W=\badmc([\exp(-\eta n),\infty))$, so $\tmu\circ F_{*}^{-i}$  for $0\le i\le n$ are supported on~$W$,
from the assumption.

We begin with 
constructing measurable partitions
$\xi(\q)$, $-1\le \q\le \q(n)$,  of  the space $W$ such that
\begin{itemize}
\item[($\Xi$1)] $\xi(\q)$ subdivides the partition $\Xi_{\badmc}$ on~$W$, which is defined
in subsection
\ref{ss-admm}.  And $\xi(q)$ is increasing with respect to $q$, that is, $\xi(\q+1)$ subdivides
$\xi(\q)$. 
\item[($\Xi$2)] Each element of the partition $\xi(\q)$ is of the form $\{\curve\}\times J$ where $\curve$ is
an admissible curve in $
\admc(a)$ with $a\ge\exp(-\eta n)$ and $J$ is an interval in $[0,a]$ such that
 $\delta(\q)\le |F^{\tau(\q)}_{*}(\curve|_{J})|\le 2 \delta(\q)$.

\end{itemize} 
The construction is done by induction on $q$ easily. Since $\delta(-1)<\exp(-\eta
n)$, we can construct a partition  $\xi(-1)$ that satisfies ($\Xi$1) and ($\Xi$2) by
subdividing the partition
$\Xi_{\badmc}$ on $W$. Let $0\le q\le q(n)$ and suppose that we have constructed the partitions
$\xi(j)$ for $-1\le j< \q$. For each  element
$\{\curve\}\times J$ of 
$\xi(\q-1)$, the length of the curve $F^{\tau(\q)}_{*}(\curve|_{J})$ is not less
than 
\[
\delta(\q-1)\cdot \exp(\lambda\subg (\tau(q)-\tau(q-1)) -c\subg)>\delta(\q),
\] provided that we take the constant $p_0$
so large that 
$(\lambda\subg-4\eta)p_0 >c\subg$. (Recall the remark on the choice of the constant $p_0$ in the last
subsection.) Hence we can construct the partition
$\xi(\q)$ satisfying ($\Xi$1) and ($\Xi$2) by subdividing
$\xi(\q-1)$. 

A Borel measurable subset in  $W$ is said to be 
a $\xi(\q)$-subset if it is a union of
elements of
$\xi(\q)$. 
Note that, if $Y$ is a $\xi(\q)$-subset, the measure $(\tilde{\mu}|_{Y})\circ
F_{*}^{-\tau(\q)}\circ \Pi^{-1}$ is contained in $\admm([\delta(\q),2\delta(q)])$ from the
condition ($\Xi$2). 

For  $-1\le q\le q(n)$ and an element $P=\{\curve\}\times J$ of the partition $\xi(\q)$, we define 
\[
\k_{\q}(P):=
\min\{\k(F^{\tau(\q)}(\curve(t)))\mid t\in J\;\}\ge k_0,
\]
where $\k(\cdot)$ is that defined in the last subsection. 
For simplicity, we denote
\[
\|\tilde\nu\|_{\rho}:=\|\tilde\nu\circ
\Pi^{-1}\|_{\rho}\quad \mbox{ for a measure $\tilde\nu$ on $W$.}
\]
The following is a consequence of the two inequalities in subsection
\ref{ss-lyine}.
\begin{sublemma} \label{lemma-bilin} Let $Y$ be  a  $\xi(\q)$-subset in $W$  for some  $-1\le \q\le \q(n)$
and  $k$ an
integer such that
\begin{equation}\label{cond-k}
k_{0}\le k\le k_0+\eta (n-\tau(\q)).
\end{equation}
If\/
 $\k_{\q}(P)\le k$ for all elements \/$P\in \xi(\q)$ that is contained in $Y$, 
  we  have
\[
\left\|(\tilde{\mu}|_{Y})\circ F_{*}^{-\tau(\q+1)}\right\|_{\rho(\q+1)}
\le
 \exp(10\lmax p_{0}+6 (k-k_0)) \left\|(\tilde{\mu}|_{Y})\circ F_{*}^{-\tau(\q)}\right\|_{\rho(\q)}.
\]
Moreover, if $k=k_0$ and $\q\ge 0$ in addition, we have either  
\begin{align*}
\left\|(\tilde{\mu}|_{Y})\circ F_{*}^{-\tau(\q+1)}\right\|_{\rho(\q+1)}&\le 
\exp(-48\epsilon p_{0})\left\|(\tilde{\mu}|_{Y})\circ F_{*}^{-\tau(\q)}\right\|_{\rho(\q)}
\intertext{or}
\left\|(\tilde{\mu}|_{Y})\circ F_{*}^{-\tau(\q+1)}\right\|_{\rho(\q+1)}&\le 
\delta(\q)^{-1}\exp(3\lmax p_0) \cdot \tilde{\mu}(Y).
\end{align*}
\end{sublemma}
\begin{proof}
We put $p=\tau(\q+1)-\tau(\q)\le
p_{0}$. 
So $p$ is smaller than $p_0$ only if $\q=-1$. From the assumption, we can divide the subset $Y$ into $\xi(\q)$-subsets $Y(\ell)$, $1\le\ell\le
\ell_{0}$, such that
 $\Pi\circ F^{\tau(\q)}_{*}(P)\cap \Lambda(\chi(\ell),\epsilon,k,p_{0};F)\neq \emptyset$ for
each $P\in \xi(\q)$ that is contained in
$Y(\ell)$. The measures  $(\tmu|_{Y(\ell)})\circ F^{-\tau(q)}\circ \Pi^{-1}$ belong to
$\admm([\delta(q),\infty))$, as we noted above.

We prove the first claim.  By using (\ref{cond-k}) and (\ref{xinclude}), we can check
\[
2\delta(\q)\le \kappa\subg^{-1}\rho\sube
\exp((\chicl(\ell)-\chiuu(\ell)-5\epsilon)p_0-4k),
\]
provided that $p_0$ is larger than some constant that depends only
on
$k_0$, $\rho\sube$,
$\kappa\subg$ and $\lmax$. This and the claim (v) and (vi) of lemma~\ref{prop-rho} imply that
the subset
$\proj\circ F_{*}^{\tau(\q)}(Y(\ell))$ is contained in  
$\Lambda(\chi(\ell),\epsilon,k+1,p_{0};F)$ and, hence, is  contained   in 
$\Lambda(\chi(\ell),\epsilon,k+1+\epsilon p_0,p;F)$ even in the case $p<p_0$, by
(\ref{eqn-shiftlambda}).

For simplicity, we put 
\[
\delta:=10\kappa\subg\rho(\q+1)\exp(-\chicl(\ell)p +k+\epsilon
p_0+2).
\]
We have 
\begin{align}
&\delta<\kappa\subg^{-1}\rho\sube \exp((\chicl-\chiuu-5\epsilon)p-4(k+1+\epsilon p_0)),\label{est:delta}\\
&\rho(q)<\delta<\delta(q) \quad \mbox{and}\\
&0<\rho(\q+1)\le \rho\sube\exp((\chicl(\ell)-5\epsilon)p -3(k+2+\epsilon p_0))/(10\kappa\subg^2)
\end{align}
from (\ref{cond-k}) and (\ref{xinclude}), provided that  $p_0$ is larger than some constant which depends only on
$k_0$,  $\rho\sube$, $\kappa\subg$, $c\subg$ and $\lmax$. 
Since the subset  $\proj\circ F_{*}^{\tau(\q)}(Y(\ell))$ is contained in $\Lambda(\chi(\ell),\epsilon,k+1+\epsilon p_0,p;F)$ as we noted, the claims (v) and (vi) of lemma~\ref{prop-rho} and the inequality (\ref{est:delta}) imply that the $\delta$-neighborhood of $\proj\circ F_{*}^{\tau(\q)}(Y(\ell))$ is contained in $\Lambda(\chi(\ell),\epsilon,k+2+\epsilon p_0,p;F)$.  From corollary \ref{prop-back-card}, it follows
\[
\max_{w\in \man}\#(F^{-p}(w)\cap
\ball(\proj\circ F_{*}^{\tau(\q)}(Y(\ell)),\delta))<\exp(6\lmax p_{0}+6k)
\]
provided that $p_0$ is larger than some constant that depends only on $\kappa\sube$ and $\lmax$.
Now we can apply lemma
\ref{lem-dir} and obtain
\begin{align*}
&\left\|(\tilde{\mu}|_{Y(\ell)})\circ F_{*}^{-\tau(\q+1)}\right\|_{\rho(\q+1)}^2\\
&\qquad \le I\subg\cdot \exp(16\lmax p_{0}+6(k+\epsilon p_0 +1)+6k)
\left\|(\tilde{\mu}|_{Y(\ell)})\circ F_{*}^{-\tau(\q)}\right\|_{\delta}^2\\ 
&\qquad \le \ell_0^{-2}
\exp(20\lmax p_{0}+12 (k-k_0)) \left\|(\tilde{\mu}|_{Y})\circ F_{*}^{-\tau(\q)}\right\|_{\rho(\q)}^2
\end{align*}
using (\ref{normess}), provided that  $p_0$ is larger than some constant which depends only on
 $I\subg$
$k_0$, $\ell_0$ and
$\lmax$. Summing up the square root of the both sides over $1\le \ell\le \ell_0$, we obtain the first claim.

We prove the second claim by using lemma \ref{prop-ly}. 
Note that $\proj\circ F_{*}^{\tau(\q)}(Y(\ell))$ is contained in
$\Lambda(\chi(\ell),\epsilon,k_0+1,p_{0};F)$ in this case, from the argument above.
We can check 
\[
\rho(\q+1)\exp((-\chicl(\ell)+\epsilon)p_{0})<\delta(\q)<\exp((\chicl(\ell)-2\chiuu(\ell)-3\epsilon)p_0),
\]
provided that $p_0$ is larger than some constant which depends only on  $c\subg$ and $\lmax$.
Recall that we took $p_0$ so large that $p_0\ge n_{*}(\chi(\ell),\epsilon,k_0+1)$ in the condition (K6).  Hence we
can apply lemma
\ref{prop-ly} and obtain
\begin{align*}
\left\|(\tilde{\mu}|_{Y(\ell)})\circ F_{*}^{-\tau(\q+1)}\right\|_{\rho(\q+1)}^{2}&\le 
\exp(-98\epsilon p_{0})\left\|(\tilde{\mu}|_{Y(\ell)})\circ
F_{*}^{-\tau(\q)}\right\|_{\rho(\q+1)}^{2}\\ &\qquad\qquad   +
\delta(\q)^{-2}\exp((-2\chicu(\ell)+2\epsilon) p_0) \cdot\tilde{\mu}(Y(\ell))^2,
\end{align*}
where we used the condition (K6)(a) in the choice of $p_0$.
This implies
\begin{align*}
\left\|(\tilde{\mu}|_{Y(\ell)})\circ F_{*}^{-\tau(\q+1)}\right\|_{\rho(\q+1)}&\le 
\exp(-49\epsilon p_{0})\left\|(\tilde{\mu}|_{Y})\circ
F_{*}^{-\tau(\q)}\right\|_{\rho(\q+1)}\\ &\qquad\qquad   +
\delta(\q)^{-1}\exp((-\chicu(\ell)+\epsilon) p_0)  \cdot\tilde{\mu}(Y).
\end{align*}
Summing up the both sides for $1\le \ell\le \ell_0$ and using (\ref{normess}), we conclude 
\begin{align*}
\left\|(\tilde{\mu}|_{Y})\circ F_{*}^{-\tau(\q+1)}\right\|_{\rho(\q+1)}&\le 
C_{0}\ell_0\cdot \exp(-49\epsilon p_{0})\left\|(\tilde{\mu}|_{Y})\circ
F_{*}^{-\tau(\q)}\right\|_{\rho(\q)}\\ &\qquad\qquad   +
\ell_{0}\cdot \delta(\q)^{-1}\exp((2\lmax +\epsilon) p_0) \cdot \tilde{\mu}(Y).
\end{align*}
The second claim follows from this inequality,
 provided that  $p_0$ is larger than some constant that depends only on $\ell_0$, $\lmax$ and
$\epsilon$. 
\end{proof}

For integers $-1\le \q'\le \q\le  \q(n)$, let
$\mathcal{K}(\q',\q)$ be the set of sequences
$\sigma=(\sigma_{j})_{j=\q'}^{\q-1}$ of $(q-q')$ integers that satisfy
\begin{equation}\label{sigeq}
0\le \sigma_j \le \eta (n-\tau(j)) \qquad \mbox{for $\q'\le j< q$.}
\end{equation}
In the case $\q'=\q$, we regard that $\mathcal{K}(\q',\q)=\mathcal{K}(\q,\q)$ consists of one empty
sequence, which we   denote by $\emptyset_{\q}$. We put
\[
\mathcal{K}(\q)=\bigcup\left\{\;\mathcal{K}(\q',\q)\mid {-1\le \q'\le \q}\;\right\}
\]
 for $0\le \q\le \q(n)$. 
Below we construct subsets $\D(\sigma)$ in $W$ for $\sigma\in \bigcup_{q=-1}^{q(n)}\mathcal{K}(\q)$ 
so  that the following conditions hold:
\begin{itemize}
\item[(D1)] $\D(\sigma)$ for $\sigma\in \mathcal{K}(\q)$ are  mutually disjoint $\xi(\q-1)$-subsets. 
\item[(D2)] The union of $\D(\sigma)$ for all $\sigma\in
\mathcal{K}(\q)$ contains the subset $\Pi^{-1}(\phys(n))\cap W$.
\item[(D3)] For  $-1\le \q'<\q\le q(n)$ and  $\sigma=(\sigma_{j})_{j=\q'}^{\q-1}\in
\mathcal{K}(\q',\q)$,  we have
\[
\left\|(\tilde{\mu}|_{\D(\sigma)})\circ  F_{*}^{-\tau(\q)}\right\|_{\rho(\q)}
\le
 \exp(10\lmax p_{0}+6 \sigma_{\q-1}) \left\|(\tilde{\mu}|_{\D(\sigma')}) \circ
F_{*}^{-\tau(\q-1)} \right\|_{\rho(\q-1)}
\]
 where $\sigma'=(\sigma_{j})_{j=\q'}^{\q-2}\in \mathcal{K}(\q',\q-1)$ (so $\sigma'=\emptyset_{q'}$ if $q'=q-1$).
Further, it holds
\[
\left\|(\tilde{\mu}|_{\D(\sigma)})\circ F_{*}^{-\tau(\q)}\right\|_{\rho(\q)}\le \exp(-48 \epsilon
p_{0})\left\|(\tilde{\mu}|_{\D(\sigma')})\circ F_{*}^{-\tau(\q-1)}\right\|_{\rho(\q-1)}
\]
in the case where $\q\ge 1$ and $\sigma_{\q-1}=0$.
\item[(D4)] For the  empty sequence $\emptyset_{\q}\in \mathcal{K}(\q,\q)$ for $q\ge 0$, we have
\[
\left\|(\tilde{\mu}|_{\D(\emptyset_{\q})})\circ F_{*}^{-\tau(\q)}\right\|_{\rho(\q)}\le 
\delta(\q-1)^{-1}\exp(3\lmax p_0)\cdot \tilde{\mu}(\D(\emptyset_{\q})).
\]
\end{itemize}
The construction is done by induction on $q$. 
For the case $\q=-1$, we just define $\D(\emptyset_{-1})=W$. For the case  $\q=0$, we have to
define
$\D(\sigma)$ for
$\sigma=\emptyset_{0}\in
\mathcal{K}(0,0)$ and $\sigma=(\sigma_{-1})\in \mathcal{K}(-1,0)$, where $0\le \sigma_{-1}\le \eta n$ from
(\ref{sigeq}). We let
$\D(\emptyset_{0})$ be the empty set  and put
\[
 \D((\sigma_{-1}))= \bigcup\{P\in \xi(-1)\mid \k_{(-1)}(P)=k_0+\sigma_{-1}\}\quad \mbox{ for 
$0\le \sigma_{-1} \le \eta n$}.
\]
Then the conditions (D1) and (D4) hold obviously. The condition (D2) follows from the condition (R3) in the
definition of the subset
$\phys(n)$.  The
first claim of  sublemma~\ref{lemma-bilin} implies that the condition  (D3) holds also.

Next, let $\q\ge 1$ and suppose that we have defined $D(\sigma)$ for $\sigma\in \mathcal{K}(\q-1)$
so that the conditions (D1)-(D4) hold for them. Consider  an element  $\sigma=(\sigma_{j})_{j=\q'}^{\q-1}$ 
in $\mathcal{K}(\q',\q)$ with
$\q'<\q$ and put
$\sigma'= (\sigma_{j})_{j=\q'}^{\q-2}\in \mathcal{K}(\q',\q-1)$. Let us set
\begin{equation}\label{dtiled}
\D_{*}(\sigma)=\bigcup\left\{P\in \xi(\q-1)\mid P\subset \D(\sigma')\mbox{ and
}\k_{q-1}(P)=k_0+\sigma_{\q-1}\right\}.
\end{equation}
In the case 
$\sigma_{\q-1}>0$, we  define $\D(\sigma)=\D_{*}(\sigma)$. In the case $\sigma_{\q-1}=0$,
we define $\D(\sigma)$ in the following manner:
From the second claim of sublemma~\ref{lemma-bilin}, we have either
\begin{align}
\left\|({\tilde \mu}|_{\D_{*}(\sigma)})\circ F_{*}^{-\tau(\q)}\right\|_{\rho(\q)}
&\le 
\exp(-48 \epsilon p_{0})
\left\|({\tilde \mu}|_{\D_{*}(\sigma)})\circ F_{*}^{-\tau(\q-1)}\right\|_{\rho(\q-1)}\label{eqn-alt1}\\
\intertext{or}
\left\|({\tilde \mu}|_{\D_{*}(\sigma)})\circ F_{*}^{-\tau(\q)}\right\|_{\rho(\q)}
&\le 
\delta(\q-1)^{-1}\exp(3\lmax p_0) \cdot {\tilde \mu}(\D_{*}(\sigma)).
\end{align}
We define  
\[
\D(\sigma)=\begin{cases}
\D_{*}(\sigma),&\mbox{ in the case (\ref{eqn-alt1}) holds;}\\
\mbox{emtpy set},&\mbox{ otherwise.}
\end{cases}
\]
Finally we define $\D(\emptyset_{\q})$ as the union of $\D_{*}(\sigma)$ for the sequences 
$\sigma=(\sigma_{j})_{j=\q'}^{\q-1}$ in $ \bigcup_{-1<q'<q}\mathcal{K}(\q',\q)=\mathcal{K}(\q)\setminus
\{\emptyset_q\}$ such that 
$\sigma_{\q-1}=0$ and that (\ref{eqn-alt1})  does {\em not} hold. As a consequence of  this definition, the
condition (D4) holds for the empty sequence~$\emptyset_q$. The condition (D1) holds obviously. We can check the
condition (D2) by using the condition
(R2) in the definition of the subset $\phys(n)$.
 The condition (D3) follows from the first claim of sublemma \ref{lemma-bilin} and  the
construction above. We have finished the  definition of the subsets
$\D(\sigma)$.

For $-1\le \q'\le \q(n)$, let $\mathcal{K}_{*}(\q')$
be the set of sequences 
$\sigma=(\sigma_{j})_{j=\q'}^{\q(n)-1}$ in
$\mathcal{K}(\q', \q(n))$ that satisfy the conditions
\[
|\sigma|_{0}\le \eta (\q(n)-\q')\quad\mbox{ and}\qquad 
|\sigma|_{1}\le 2\eta (\q(n)-\q') p_0
\]
where 
\[
|\sigma|_{0}:=\#\{\q'\le j< \q(n)\mid j\ge 0\mbox{ and }\sigma_{j}>0\},\qquad 
|\sigma|_{1}:=\sum_{j=\q'}^{\q(n)-1}\sigma_{j}.
\]
Then, from the definition of the subsets $\phys(n)$ and $\D(\sigma)$, we have 
\[
\Pi^{-1}(\phys(n))\cap W\subset \bigcup_{\q'=-1}^{\q(n)}\bigcup_{\sigma\in \mathcal{K}_{*}(\q')}\D(\sigma) 
\]
and hence
\[
(\tilde\mu|_{\Pi^{-1}(\phys(n))})\circ F_{*}^{-n}\le \sum_{\q'=-1}^{\q(n)}\sum_{\sigma\in \mathcal{K}_{*}(\q')}
(\tilde\mu|_{\D(\sigma)})\circ  F_{*}^{-n}.
\]
 For each  $\sigma=(\sigma_{j})_{j=\q'}^{\q(n)-1}$ in $\mathcal{K}_{*}(\q')$ with $q'\ge 0$, we can obtain
\begin{align*}
&\left\|(\tilde\mu|_{\D(\sigma)})\circ F_{*}^{-n}\right\|_{\rho}=
\left\|(\tilde\mu|_{\D(\sigma)})\circ F_{*}^{-\tau(q(n))}\right\|_{\rho(q(n))}\\
&\quad
\le \exp\left(10\lmax p_0 |\sigma|_{0}+6 |\sigma|_{1}-48\epsilon(\q(n)-\q'-|\sigma|_{0})p_0\right)\cdot
\left\|(\tilde\mu|_{\D(\emptyset_{\q'})})\circ  F_{*}^{-\tau(q')}\right\|_{\rho(\q')}
\end{align*}
 from  the conditions (D3) and, hence, 
\[ 
\left\|(\tilde\mu|_{\D(\sigma)})\circ F_{*}^{-n}\right\|_{\rho}\le
\exp(-45\epsilon(\q(n)-\q')p_0+11\lmax p_0+c\subg)\cdot |\tilde\mu|
\]
from the condition (D4) and the choice of $\eta$. Similarly,  for   $\sigma=(\sigma_{j})_{j=-1}^{\q(n)-1}$ in
$\mathcal{K}_{*}(-1)$,  we can obtain
\[
\left\|(\tilde\mu|_{\D(\sigma)})\circ F_{*}^{-n}\right\|_{\rho}
\le \exp\left(6\lmax p_0 (|\sigma|_{0}+1)+3 |\sigma|_{1}-48\epsilon(\q(n)-|\sigma|_{0})p_0\right)\cdot
\left\|\tilde\mu\right\|_{\rho(-1)}
\]
and hence
\[
\left\|(\tilde\mu|_{\D(\sigma)})\circ F_{*}^{-n}\right\|_{\rho}\le
 \exp(-45\epsilon n+8\lmax p_0)\cdot
\left\|\tilde\mu\right\|_{\rho(-1)}.
\]

For the cardinality of the set $\mathcal{K}_{*}(\q')$, we have 
\[
\#\mathcal{K}_{*}(\q')\le 
\begin{pmatrix}
\q(n)-\q'\\
{[}\eta (\q(n)-\q')]
\end{pmatrix} \cdot 
\begin{pmatrix}[2\eta p_{0}(\q(n)-\q')]+[\eta (\q(n)-\q')]\\ 
{[}\eta (\q(n)-\q')]\end{pmatrix}
\]
where the first factor on  the right hand side is an upper bound for the number of
possible arrangements of integers $j\ge 0$ for which $\sigma_{j}$ may be positive and the second
factor is an upper bound for the cardinality of $\sigma\in \mathcal{K}_{*}(\q')$ when one of such
arrangements is given. For positive numbers $\alpha,\beta>0$ and an integer $m\ge 1$ such that $\alpha m\ge 1$ and
$\beta m\ge 1$,  we have 
\[
\log \begin{pmatrix}
\alpha m+\beta m\\
\beta m
\end{pmatrix}
\le  \alpha m\log \left(1+\frac{\beta}{\alpha}\right)+\beta m \log \left(1+\frac{\alpha}{\beta}\right)
+A_{0}
\]
from 
Stirling's formula, where $A_{0}$ is  an absolute constant. Hence  we
can obtain 
\begin{align*}
\frac{\log \#\mathcal{K}_{*}(\q')}
{q(n)-q'}
&\le -(1-\eta)\log(1-\eta)-\eta\log \eta\\
&\qquad\qquad  +2\eta p_0\log\left(1+\frac{1}{2p_0}\right)
+\eta \log(1+2p_0)+2A_{0}
\end{align*}
for $-1\le q'<q(n)$. This implies 
\[
\#\mathcal{K}_{*}(\q')\le \exp(\epsilon p_0 (q(n)-q'))\quad\mbox{for $-1\le q'< q(n)$,}
\]
provided that $p_0$ is larger than some constant which depends only on~$\epsilon$ and $\eta$.
Now we can conclude
\begin{align*}
&\left\|(\mu|_{\phys(n)})\circ F^{-n}\right\|_{\rho} =\left\|(\tilde\mu|_{\Pi^{-1}(\phys(n))})\circ
F_{*}^{-n}\right\|_{\rho}\\
 & \le
\left(
\sum_{\q'=0}^{\q(n)}\sum_{\sigma\in
\mathcal{K}_{*}(\q')}\left\|(\tilde\mu|_{\D(\sigma)})\circ F_{*}^{-n}\right\|_{\rho}
\right)
+
\left(
\sum_{\sigma\in
\mathcal{K}_{*}(-1)}\left\|(\tilde\mu|_{\D(\sigma)})\circ F_{*}^{-n}\right\|_{\rho}
\right)
\\ 
& \le \sum_{\q'=0}^{\q(n)} \exp(-44\epsilon
(\q(n)-\q')p_{0}+11\lmax p_0+c\subg ) |\mu|+ \exp(-44\epsilon n+8\lmax p_0)\|\mu\|_{\rho(-1)}.
\end{align*}
This implies the inequality in lemma~\ref{prop-abs2}.  


\section{Genericity of the transversality condition on unstable cones}\label{sec-trans}

In this section, we consider multiplicity of  
tangencies between the images of the unstable cones under iterates of mappings in $\U$, and investigate to what
extent we can resolve the tangencies by perturbation.  The goal is  the proof of  
 theorem~\ref{th-ge-trans}. The point of our argument in this section is that the
dominating expansion in the unstable direction
 acts as  uniform contraction on the angles between subspaces in the
unstable cones. This enables us to control the images of the unstable cones  in perturbations of 
 mappings in
$\U$. Notice that the content and the notation in this section is independent  of those in the
last section.

\subsection{Reduction of   theorem
\ref{th-ge-trans}: The first step}\label{ssthge}
In this subsection and the next, we reduce theorem \ref{th-ge-trans} to more tractable
propositions in two steps. For a quadruple 
$\bchi=(\chicl,\chicu,\chiul,\chiuu)$, we put
\[
 \chi_{c}^{++}:=\max\{\chicu,0\}, \quad \chicdif:=\chicu-\chicl\quad\mbox{ and }
\quad \chiudif:=\chiuu-
\chiul.
\]
For a quadruple $\chi$  satisfying 
(\ref{chi-order}) and a positive number $\epsilon$, let
$\ss_{1}(\bchi,\epsilon)$ be the set of mappings
$F\in
\U$ that satisfy
\begin{equation}\label{inecl}
\limsup_{n\to\infty}n^{-1}\log
\mult(\chi,\epsilon,  \epsilon n,n;F)\ge  \chicl+\chiul-\chicdif-\chiudif-\epsilon.
\end{equation}
The first step of the reduction is simple. We show that we can deduce theorem~\ref{th-ge-trans} from the
following proposition:
\begin{proposition}
\label{prop-ge-trans} 
Suppose that $s\ge r+3$ and let $\mom_s$ be the measure on $C^r(M,\Real^2)$ introduced  in 
lemma~\ref{lem-mom}. The subset  $\ss_{1}(\bchi,\epsilon)$  is shy
with respect to 
the measure $\mom_s$ for $s\ge r+3$, if the quadruple
$\bchi=(\chicl,\chicu,\chiul,\chiuu)$   satisfies the conditions 
(\ref{chiorder1}) and
\begin{equation}\label{chiorder3}
-2\lmax<\chicl<\chicu<\chiul<\chiuu<2\lmax,
\end{equation}
\begin{equation}\label{chidiforder3}
\chicdif+\chiudif<\chi_{c}^{-}+\chi_{u}^{-},
\end{equation} 
\begin{equation}
\label{chiorder2}
\chi_{u}^{-}+\chi_{c}^{-}-\chi_{c}^{++}>\left(\frac{\chi_{c}^{++}+\chi_{u}^{+}}
{\chi_{c}^{-}+\chi_{u}^{-}-\chicdif-\chiudif}+1\right)\cdot (\chicdif
+\chiudif)
\end{equation}
and if $\epsilon>0$ is smaller than some constant which depends only on $\bchi$  and $s$ besides the integer
$r\ge 2$ and the objects that we  fixed in subsection
\ref{cnbl}.
\end{proposition}
Below we prove theorem \ref{th-ge-trans} assuming this proposition.
\begin{proof}[Proof of theorem \ref{th-ge-trans}]
  For any point $(\chi_c,\chi_u)$
in the subset  given in the claim (a):
\[
\left\{\;(x_{c},x_{u})\in \Real^2 \mid \; x_{c}+x_{u}>0,\;
\lambda\subg\le x_{u}\le
\lmax,\; x_{c}\le 0\;\right\},
\] 
we can take a quadruple $\bchi=(\chicl,\chicu,\chiul,\chiuu)$  satisfying the
conditions (\ref{chiorder1}),  (\ref{chiorder3}), (\ref{chidiforder3}) and
(\ref{chiorder2}) such that the rectangle $(\chicl,\chicu)\times (\chiul,\chiuu)$
contains the point $(\chi_c,\chi_u)$. 
  Thus we  can  choose a countable collection $\X$ of
quadruples that satisfy (\ref{chiorder1}),  (\ref{chiorder3}), (\ref{chidiforder3}) and
(\ref{chiorder2}) such that
 the conditions (a) and (b) in theorem \ref{th-ge-trans} hold. We are going to show  the condition (c) in
theorem \ref{th-ge-trans}. We fix $s\ge r+3$. Let $\X'$ be an 
 arbitrary finite subset  of 
$\X$.  Then we can take positive number $\epsilon>0$ so small that the conclusion of proposition
\ref{prop-ge-trans} holds for all the quadruples in 
$\X'$.  For each $\bchi\in \X'$ and $n\ge 1$,  let $\ss_{1}^{*}(\bchi,
\epsilon, n)$ be the closed subset of mappings
$F\in
\U$ that satisfy
\[
\mult(\chi,\epsilon, \epsilon n,n;F)\ge \exp\left(\left(\chicl+\chiul-\chicdif-\chiudif-\epsilon\right)n\right).
\]
If a mapping $F\in \U$ belongs to $\ss_1(\X')$, or $F$ does not satisfy the transversality condition on
unstable cones for
$\X'$, it holds
\begin{equation*}
 \liminf_{n\to\infty}\max
\left\{
\frac{\log(\mult(\chi,\epsilon, \epsilon n, n;F))}
{n\cdot (\chicl+\chiul-\chicdif-\chiudif)}\; ;\; \bchi=(\chicl,\chicu,\chiul,\chiuu)\in \X'\right\}\ge 1
\end{equation*}
because $\mult(\chi,\epsilon, k,n;F)$ is increasing with
respect to
$\epsilon$ and $k$. Hence we have
\[
\ss_{1}(\X')\subset \bigcup_{m>0}\bigcap_{n>m}\bigcup_{\bchi\in \X'}\ss_{1}^{*}(\bchi,
\epsilon, n)\subset
\bigcup_{\bchi\in
\X'}\ss_{1}(\bchi,\epsilon). \]
From proposition
\ref{prop-ge-trans}, the subset $\bigcup_{\bchi\in
\X'}\ss_{1}(\bchi,\epsilon)$  is shy 
with respect to the measure $\mom_s$ and, hence, so is $\ss_{1}(\X')$. Further, 
the closed subset
$\bigcap_{n>m}\bigcup_{\bchi\in
\X'}\ss_{1}^{*}(\bchi,
\epsilon, n)$ is nowhere dense, because it is shy with respect to
the measure $\mom_s$. Thus  $
\ss_1(\X')$ is a meager subset in $\U$ in the sense of Baire's category argument. 
\end{proof}

\subsection{Reduction of   theorem
\ref{th-ge-trans}: The second step}\label{ssthge2}
The second step of the reduction is rather involved. We reduce proposition
\ref{prop-ge-trans} to yet another proposition, proposition \ref{ss12}, which will be  proved in the remaining
part of this section. Below we consider an integer $s\ge r+3$, a  quadruple 
$\chi=\{\chicl,\chicu,\chiul,\chiuu\}$ and a positive number $\epsilon$. We assume that the
quadruple $\bchi$  satisfies
the assumptions in proposition~\ref{prop-ge-trans}, that is, the conditions (\ref{chiorder1}), 
(\ref{chiorder3}), (\ref{chidiforder3}) and (\ref{chiorder2}). 

In  this section, we will introduce 
several constants that depends only on the quadruple 
$\bchi$ and the integers $s\ge r\ge2$ besides  the objects that we  fixed in subsection
\ref{cnbl}. In order to distinguish such kind of constants, we will denote them by symbols with subscript
$\chi$. Also we will use a generic symbol $C\subz$ for large positive constants of this kind. 
The usage of these notations is the same as those introduced
in subsection~\ref{ssremc} and section \ref{sec-dist}.

The choice of the number
$\epsilon>0$ is important for our argument not only in this subsection but also in the remaining
part of  this section.
 We claim that our argument in this section  is true if
 $\epsilon$ is smaller than some constant~$
\epsilon\subz$.  Below we will assume $0<\epsilon\le \epsilon\subz$ and give the conditions
on the choice of $\epsilon\subz$ in the course of the argument. 

From  the condition (\ref{chiorder2}), we can take and fix a positive constant
$h\subz$ such  that
\[
 h\subz+1 > \frac{\chi_{c}^{++}+\chiuu}
{\chicl+\chiul-\chicdif-\chiudif}
\]
and that
\[
\chiul+\chicl-\chi_{c}^{++}>(h\subz+2)(\chicdif+\chiudif).
\]
Then we  fix a positive integer $q\subz$ such that
\[
q\subz
>\frac{2(\chiul-\chicl)+\chi_{c}^{++}-\chicl+\chicdif+2\chiudif}{\chiul+\chicl-\chi_{c}^{++}
-(h\subz+2)(\chicdif+\chiudif)}.
\]
Also we  put 
\begin{equation}\label{defrx}
r\subz =100 (h\subz+1)^2\lmax^{2}/\lambda\subg \ge 100.
\end{equation}
\begin{de}
For  integers
$0<p<n$ and  a point 
$z\in
\man$,  let $\ss_{1}(\bchi,\epsilon, n, p, z)$ be the set of mappings $F\in \U$ such that 
there exist  a subset 
$\{w_{i}\}_{i=0}^{q\subz}$ in
$F^{-p}(z)$ and  subsets
$E_{i}$,
$0\le i\le q\subz$, in $F^{-n+p}(w_i)\subset F^{-n}(z)$ that satisfy the following three
conditions:
\begin{itemize}
\item[{\rm (S1)}] The subsets $E_i$ for $0\le i\le q\subz$ are contained in 
$\Lambda(\chi,\epsilon,2(h\subz+1)\epsilon n, n;F)$, and
\[
\#E_{i}=[ \exp(( \chicl+\chiul-\chicdif-\chiudif-r\subz \epsilon)n)]+1.
\]
\item[{\rm (S2)}] For any points
$y$ and $y'$ in the union $ \bigcup_{i=0}^{q\subz} E_{i}$, we have 
\[
\angle(DF^{n}(\E^{u}(y)), DF^{n}(\E^{u}(y')))
\le
\exp((\chicu-\chiul+6\epsilon+h\subz(\chicdif+\chiudif+4\epsilon))n).
\]
\item[{\rm (S3)}]  For $0\le
j\le  p$ and $0\le i, i'\le q\subz$, we have
\[
F^{j}(\ball(w_{i},10\exp(-r\subz \epsilon n))\cap \ball(w_{i'},10\exp(-r\subz
\epsilon n)) =\emptyset
\]
but for the case where both $i= i'$ and $j=0$ hold.
\end{itemize}
\end{de} 
For an integer $n\ge 1$, we consider the lattice 
\[
\lattice_{n}=\lattice(\exp((\chicl-\chiul)n))
\]
where $\lattice(\cdot)$ is that defined in subsection
\ref{sec-notation}.  
The following lemma is the main ingredient of this subsection:
\begin{lemma}\label{ss1} 
We have
\begin{equation}\label{borelc}
\ss_1(\chi,\epsilon)\subset \limsup_{n\to \infty}\left(\bigcup_{p}\bigcup_{z\in \lattice_n}
\ss_{1}(\bchi,\epsilon, n, p, z)\right),
\end{equation}
where $\cup_p$ indicates the union over integers $p$ satisfying
\begin{equation}\label{eqnp}
3h\subz(\lmax/\lambda\subg) \epsilon n\le p\le
3 h\subz(h\subz+1) (\lmax/\lambda\subg)  \epsilon n +1.
\end{equation}
\end{lemma}

\begin{proof}
Let  $F$ be a mapping in $\ss_{1}(\chi,\epsilon)$. We show that there are an arbitrarily
large integer $n$, an integer $p$ satisfying (\ref{eqnp}) and a point $z\in \lattice_n$ such that
$F$ belongs to 
$\ss_{1}(\bchi,\epsilon, n, p, z)$. From the
definition of $\ss_{1}(\chi,\epsilon)$,  there are 
 infinitely many integers~$m$ that satisfy
\begin{equation}\label{eqn:much}
\mult(\chi,\epsilon, \epsilon m,m;F)>
\exp((\chicl+\chiul-\chicdif-\chiudif-2\epsilon)m). 
\end{equation}
In the argument below, we consider a large integer $m$ satisfying the condition (\ref{eqn:much}). Note that, since we can take the integer $m$ as large as we like, we may and will 
replace $m$ by larger one if it is necessary. From the definition of $\mult(\cdot)$,  there exist a point
  $\zeta\in \man$ and a subset $P$ in $
\Lambda(\chi,\epsilon,\epsilon m,m;F)$ with cardinality
\[
\#P>\exp((\chicl+\chiul-\chicdif-\chiudif-2\epsilon)m)
\]
 such that $F^{m}(P)= \{\zeta\}$ and that
\[
\angle(DF^{m}(\E^{u}(w)),
DF^{m}(\E^{u}(w')))\le 10  H\subg\exp((\chicu-\chiul+2\epsilon)m)
\]
for  $w,w'\in P$.
We put
$p:=\left[3h\subz(\lmax/\lambda\subg)
\epsilon m\right]+1$ and consider the subsets of $P$,
\[
P_\ell(w)=\left\{w'\in P\mid F^{m-\ell p}(w')=F^{m-\ell p}(w)\right\}
\]
for $0\le \ell\le[m/p]$ and  $w\in P$.
Since the subset $P_\ell(w)$ is contained in the subset
$\Lambda(\bchi,\epsilon, (m+\ell p)\epsilon, m-\ell p;F)$ from (\ref{eqn-shiftlambda}),
we have
\begin{align}\label{eqnecard}
\#P_{\ell}(w)&\le \kappa\sube \exp((\chiuu+\chi_{c}^{++}+7\epsilon)(m-\ell
p)+6(m+\ell p)\epsilon))\\
&\le \exp((\chiuu+\chi_{c}^{++}+7\epsilon)(m-\ell
p)+7(m+\ell p)\epsilon))\notag
\end{align}
from corollary~\ref{prop-back-card}, where the second inequality holds when $m$ is sufficiently large. 
Especially,  for the case
$\ell=[m/p]$, we have
\[
\#P_{[m/p]}(w)\le 
\exp((\chiuu+\chi_{c}^{++}+7\epsilon)p+14\epsilon m))<\exp(-  [m/p] \cdot \epsilon
p) \cdot \# P
\]
where the second inequality holds if  $\epsilon\subz$ is smaller than some constant
that depends only on $\bchi$, $h\subz$, $\lmax$ and $\lambda\subg$ and if we consider sufficiently large $m$ according to the choice of $\epsilon\subz$.  Thus there exist  integers 
$0\le
\ell
< [m/p]$ such that 
\begin{equation}\label{eqn-p}
\max_{w\in P}\#P_{\ell+1}(w)
<\exp(-\epsilon p)\max_{w\in P}\#P_{\ell}(w).
\end{equation}
Let $\ell_0$ be the smallest integer $0\le
\ell
< [m/p]$ such that  (\ref{eqn-p}) holds. Then we have
\[
\max_{w\in P}\#P_{\ell_0}(w)\ge \exp(-\epsilon \ell_0 p)\cdot\#P.
\]
Take a point
$w_{0}\in P$ such that
$\#P_{\ell_0}(w_0)=\max_{w\in P}\#P_{\ell_0}(w)$, and put 
$n=m-\ell_0 p$, $z=F^{n}(w_0)$,  
$E=P_{\ell_0}(w_{0})$. 
It holds
\[
\#E=\# P_{\ell_0}(w_{0})\ge \exp(-\epsilon(m-n))\# P\ge
\exp((\chicl+\chiul-\chicdif-\chiudif-3\epsilon)m).
\]
Comparing  this with (\ref{eqnecard}) for $\ell=\ell_0$, we obtain 
\[
m<\frac{\chi_{c}^{++}+\chiuu}
{\chicl+\chiul-\chicdif-\chiudif-17\epsilon}\cdot n<(h\subz +1)n
\]
where the second inequality follows from the choice of $h\subz$ 
provided that $\epsilon\subz$ is smaller
than some constant that depends only on $\bchi$ and $h\subz$. Hence 
$n$ and
$p$ satisfy the condition (\ref{eqnp}) and it holds
\[
E\subset  \Lambda(\bchi,\epsilon, \epsilon m, m;F)\subset \Lambda(\bchi,\epsilon, (m+\ell_0 p)\epsilon,
m-\ell_0 p;F)\subset \Lambda(\chi,\epsilon,(2h\subz+1)\epsilon n,n;F).
\]
From  (\ref{anglerel}), we can obtain,
for any points $w$ and $w'$ in
$E$, 
\begin{align*}
\angle(DF^{n}&(\E^{u}(w)),\;
DF^{n}(\E^{u}(w')))\\
& \le
A\subg 
\frac
{D_{*}F^{m-n}(\eu(F^{n}(w)))}
{|D^{*}F^{m-n}(\eu(F^{n}(w)))|}
\angle(DF^{m}(\E^{u}(w)),
DF^{m}(\E^{u}(w')))\\
& \le A\subg \exp((-\chicl+\chiuu)(m-n)+2\epsilon m)\cdot 10 
H\subg\exp((\chicu-\chiul+2\epsilon)m)\\ &
= 10 H\subg A\subg \exp((\chicu-\chiul+4\epsilon)n+(\Delta\chi_{c}+\Delta\chi_{u}+4\epsilon)
(m-n) )\\
& \le
\exp((\chicu-\chiul+5\epsilon +h\subz(\chicdif+\chiudif+4\epsilon))n)
\end{align*}
provided that $m$ is sufficiently large. 

Let us consider  the subset $\{w_i\}_{i=1}^{i_0}\subset F^{-p}(z)$ of all points $w\in  F^{-p}(z)$ such that $F^{n-p}(w)\cap E\neq \emptyset$ . From (\ref{eqn-shiftlambda0}) and  (\ref{eqn-shiftlambda}), it is contained in 
$\Lambda
(\chi, \epsilon, \epsilon (m+p \ell_0), p;F)$. Corollary \ref{prop-back-card} gives the following estimate for its cardinarity $i_0$:
\[
i_0\le \kappa\sube\exp(5\lmax p +6\epsilon  (m+p \ell_0))\le
\kappa\sube\exp(5\lmax p+12\epsilon  m).
\]
 We denote $E_i=\{y\in E\mid
F^{n-p}(y)=w_{i}\}$ for $1\le i\le
i_{0}$, so $E=\cup_{i=1}^{i_0}E_i$.

By changing the index~$i$, we assume that the cardinality of the subset $
E_{i}$ is decreasing with respect to
$i$.  Let $i_1$ be the smallest positive integer such that
\[
\sum_{i=1}^{i_{1}}\#E_{i}>\frac{1}{2}\sum_{i=1}^{i_{0}}\#E_{i}=\frac{\# E}{2}.
\]
Then we have
$
 \# E_{i_1}\cdot (i_0-i_1+1)\ge \sum_{i=i_1}^{i_0}\# E_{i} \ge \# E/2
$
and hence
\[
\#E_{i}\ge \# E_{i_1}\ge  \frac{\# E}{2i_0} >
\exp((\chicl+\chiul-\chicdif-\chiudif-r\subz\epsilon)n)\quad \mbox{for $1\le i\le i_{1}$,}
\]
where the last inequality follows from the definitions of $p$ and $r\subz$ provided that
 $m$ is sufficiently large. Also we have 
\[
i_{1}\ge\frac{\sum_{i=1}^{i_1}\#E_i}{\#E_{1}}\ge  
\frac{\# E}{2 \# E_1}\ge \frac{\exp(\epsilon
p)}{2}
\]
  from the condition (\ref{eqn-p}) for $\ell=\ell_0$.

Notice that the point $z$ that we took above may not be contained in $\lattice_{n}$, while 
it have to be. So we want to shift it to  the closest point in 
$\lattice_{n}$. The distance from the point $z$ to the closest point in $\lattice_{n}$ is bounded by
$\exp((\chicl-\chiul)n)$ and hence by
$\rho\sube\exp((\chicl -5\epsilon -3(2h\subz+1)\epsilon))n)$,
provided that  $\epsilon\subz$ is smaller than some constant which depends only on $\bchi$  and
that we took sufficiently large
$m$. Thereby, by virtue of lemma \ref{prop-rho}, we can move the points
$w_{i}$ and those in
$E_i$ accordingly so that the relations $F^p(w_i)=z$ and $F^{n-p}(E_i)=\{w_i\}$ are
preserved.
Henceforth, we consider the points $z\in \lattice_n$, $w_i$ and the subsets $E_i$ thus obtained.  Lemma
\ref{prop-rho} guarantees that the subsets $E_{i}$ are contained in $\Lambda(\bchi,\epsilon,
2(h\subz+1)\epsilon n,n;F)$ and that
\begin{align*}
\angle(DF^{n}(\E^{u}(w)),\;
DF^{n}(\E^{u}(w'))) &\le
\exp((\chicu-\chiul+5\epsilon +h\subz(\chicdif+\chiudif+4\epsilon))n)\\
&\qquad +2\kappa\sube \exp((\chicl-\chiul+(4h\subz+2)\epsilon)n)\\
&<
\exp((\chicu-\chiul+6\epsilon +h\subz(\chicdif+\chiudif+4\epsilon))n)
\end{align*}
for any points $w,w'\in \bigcup_{i=1}^{i_1}E_i$, provided that $\epsilon\subz$ is smaller than some constant which depends only on $\bchi$ and that $m$ is sufficiently large.    Up to this
point, we have found arbitrarily large integer $n$, an integer $p$,  points $z$, $w_{i}$,
$1\le i\le i_1$, and subsets $E_{i}$,
$1\le i\le i_1$, that satisfy the conditions (\ref{eqnp}), (S1) and (S2).  It remains to choose
$(q\subz+1)$ points among  
$w_{i}$, $1\le i\le i_1$, so that the condition (S3) holds.

Put $W=\{w_{i};1\le i\le i_1\}$ and $\delta=40 p\cdot \exp(2\lmax p-r\subz
\epsilon n)$.  
Note that the points~$w_i$ belong to $\Lambda(\bchi,\epsilon,2(h\subz+1)\epsilon n,p;F)$ by (\ref{eqn-shiftlambda0}). 
We can check
\[
2\delta<\kappa\subg^{-1}\rho\sube \exp((-\chiuu+\chicl-5\epsilon)p-8(h\subz+1)\epsilon n)
\]
by using the definition of $p$ and $r\subz$ and the condition (\ref{eqnp}), provided that  $m$ is
large enough.  Thus 
$F^{p}$ is a diffeomorphism on the 
$2\delta$-neighborhood of each point in $W$ from  lemma~\ref{prop-rho}(v). This implies that  the distances
between the points in
$W\subset F^{-p}(z)$ are not less than 
$2\delta$. Let $L\subset W$ be the set of points in $W$ that are within distance $\delta$ to
either of the points $F^{j}(z)$, $ 0\le j< p$. Then we have
$\# L\le p$ obviously.

Consider a sequence $J=(j_{\nu})_{\nu=0}^{\nu_{0}}$  of integers such
that
$1\le j_\nu\le p$ for $0\le \nu\le \nu_0$.  We denote  the sum of the integers in $J$ by
$|J|:=\sum_{\nu=0}^{\nu_{0}}j_{\nu}$. For $x,x'\in W\setminus L$, we denote
$x\succ_{J} x'$ if  there is  a
sequence of points $x_{0}=x,x_1,\cdots , x_{\nu_{0}+1}=x'$ in
$W\setminus L$ such that
\[
F^{j_{\nu}}(\ball(x_{\nu},10\exp(-r\subz \epsilon n))\cap \ball(x_{\nu+1},10\exp(-r\subz
\epsilon n))
\neq \emptyset\qquad \mbox{for $0\le \nu\le \nu_{0}$. }
\]
 From the definition of $\delta$ above, it is easy to see that we have 
 $
d(F^{|J|}(x),x')<\delta$ if
$x\succ_{J}x'$  for some $J$ with $|J|\le 2p$.
Hence, given a point $x\in W\setminus L$ and an integer $1\le
i\le 2p$, there is at most one point $x'$ in
$W\setminus L$ that satisfies $x\succ_{J} x'$ for some sequence $J$ with $|J|=i$. 

Actually, the relation $x\succ_{J}x'$ holds for some points $x, x'$ in
$W\setminus L$ only if
$|J|<p$. In fact, otherwise, there should be  a sequence 
$J$ with $p\le |J|< 2p$ and points $x,x'$ in $W\setminus L$ such that $x\succ_J x'$ and hence
$d(F^{|J|-p}(z), x')=d(F^{|J|}(x), x')<\delta$. But, since 
  $0\le |J|-p< p$, this contradicts
the definition of $L$.

The relation $x\succ_{J}x'$ never holds if $x=x'$. In fact, if $x\succ_{J} x$ for
some  
$J$, the relation
$x\succ_{J^i} x$ should hold for any $i\ge 1$ where $J^i$ is the iteration of 
$J$ for
$i$~times. But this obviously contradicts the fact proved in the preceding paragraph.

 Let us
denote
$x\succ x'$ for
$x,x'\in W\setminus L$ if either
$x=x'$ or 
$x\succ_{J} x'$ for some sequence $J=(j_\nu)_{\nu=0}^{\nu_{0}}$ satisfying $1\le j_\nu \le
p$. From the argument above, this
relation
 is a  partial order on the set
$W\setminus L$ and, for  each $x\in W\setminus L$, there exist at most $p$ points  $x'$ in
$ W\setminus L$ such that
$x\succ x'$. Let $W_{\max}$ be the set of the maximal elements in $W\setminus L$ with
respect to the partial order
$\succ$. Then we have
\[
\#W_{\max}\ge
\frac{\#(W\setminus L)}{p}\ge  \frac{([\exp(\epsilon p)/2]-p)}{p}\ge (q\subz+1)
\]
provided that  $m$ is large enough.
Take $(q\subz+1)$ points $\{w_i\}_{i=0}^{q\subz}$ from  $W_{\max}$, then the
condition (S3) holds for them. We have completed the proof of  lemma
\ref{ss1}. 
\end{proof}

Using lemma \ref{ss1},   we can deduce proposition \ref{prop-ge-trans} from the following proposition:
\begin{proposition} \label{ss12}
Let $s\ge r+3$. Suppose that a quadruple $\bchi$ satisfies the conditions (\ref{chiorder1}),  (\ref{chiorder3}),
(\ref{chidiforder3}) and (\ref{chiorder2}) and that a positive number $\epsilon$ satisfies $0<\epsilon\le
\epsilon\subz$.  Then, for any
$d>0$ and any mapping
$G$ in
$C^{r}(M,\torus)$, there exists an integer $n_0$ such that 
\begin{equation}\label{claim-eq}
\mom_s(\base_{G}^{-1}(\ss_{1}(\bchi,\epsilon, n, p, z))\cap
\disk^{s-3}(d))<\exp((2\chicl-2\chiul-\epsilon) n)
\end{equation}
for $n\ge n_0$,
$z\in
\lattice_n$ and $0<p<n$ that satisfies the condition (\ref{eqnp}). 
\end{proposition}
\begin{remark}
$\base_{G}$ and $\disk^{s-3}(d)$ above are those defined by (\ref{PhiG}) and (\ref{defD})  respectively.
\end{remark}
In fact, since we have $\#L_n= ([\exp((-\chicl+\chiul)n)]+1)^2$,  it follows from proposition
\ref{ss12} and (\ref{borelc}) that
\[
\mom_s\left(\base_{G}^{-1}(\ss_{1}(\bchi,\epsilon))\cap
\disk^{s-3}(d)\right)=0\quad \mbox{ for any  $d>0$ and $G\in C^{r}(M,\torus)$}.
\]
Since the measure $\mom_s$ is supported on 
$C^{s-3}(M,\Real^2)=\bigcup_{d>0}\disk^{s-3}(d)$, this implies that the subset 
$\ss_{1}(\bchi,\epsilon)$  is shy with respect to the measure
$\mom_s$.

\subsection{Perturbations}\label{secpert}
In this subsection, we introduce some families of mappings and give  estimates
on the variations of the images of the unstable subspaces $\E^u(z)$ under iterates of the mappings in the
families. Henceforth, in this subsection and the next, we consider the situation in proposition~\ref{ss12}: 
Let  $s\ge r+3$;  Let  $\bchi$ a quadruple that satisfies the conditions (\ref{chiorder1}),  (\ref{chiorder3}),
(\ref{chidiforder3}) and (\ref{chiorder2}) and $\epsilon$ a positive number $\epsilon$ that satisfies
$0<\epsilon\le
\epsilon\subz$.

 Take and fix a
$C^{\infty}$ function
$\psi:\Real^{2}\to \Real$ such that $\|\psi\|_{C^1}\le 1$ and that 
\[
\psi(w)=\begin{cases} x,&\mbox{if $\|w\|\le 1/10$;}\\
0, &\mbox{if $\|w\|\ge 1$}
\end{cases}\qquad \mbox{for $w=(x,y)\in \Real^2$.}
\]
For each point $z\in \man$, we consider an isometric embedding
\[
\varphi_{z}:\{w\in
\Real^2\mid
\|w\|<1\}\to
\torus
\]
that carries the origin to
$z$ and the $x$-axis $\Real\times \{0\}$ to $\E^{u}(z)$. For $n\ge 1$, we put 
\[
\delta_n=\exp(-r\subz \epsilon n).
\]

Recall that we took the
subset $\U$ of mappings  as a neighborhood of a
$C^r$mapping $\bF$ in subsection \ref{cnbl}. 
For an integer $n\ge 1$ and a point $z\in M$, we define the
$C^\infty$mapping 
$\psi_{n,z}:\man\to \Real^{2}$  by
\[
\psi_{n,z}(w):=
\begin{cases}
\delta_n^{s+3}\cdot \psi\left(\varphi^{-1}_{z}(w)/\delta_n\right)\cdot
\e^{c}(\bF(z)),&\mbox{if $d(w,z)<\delta_n$;}\\
0,&\mbox{otherwise}
\end{cases}
\]
where $\e^{c}(\cdot)$ is either of the two unit vectors in the central
subspace $\E^{c}(\cdot)$.  Note that, for any mapping $F\in \U$, the parallel translation of the vector
$\e^{c}(\bF(z))$ to $F(z)$ is contained in 
$\cone^c(F(z))$  from the choice of 
the constant $\rho\subg$  in subsection \ref{cnbl}.
\begin{remark} Notice that the definition of  $\psi_{n,z}(w)$  does not depend on $F\in \U$.
\end{remark}

Let $n$ and $p$ be  positive integers that satisfy the
condition (\ref{eqnp}),
$S=\{x_i\}_{i=0}^{q\subz}$ an  ordered subset of the lattice
 $\lattice(\delta_{n}/40)$ and $F$ a mapping in $\U$. 
The family of mappings that we are going to consider is
\begin{equation*}\label{pert}
F_{\t}(w)=F(w)+\sum_{i=1}^{q\subz}t_{i}\cdot\psi_{n,x_i}(w)
\;: \;M\to \torus
\end{equation*}
where $\t=(t_{i})\in \Real^{q\subz}$ is the parameter that ranges over the region
\[
R=\{\t=(t_{i})\in \Real^{q\subz}\mid |t_{i}|\le \exp(\chicl n)\}.
\]
For this family, we have
\begin{equation}\label{magn}
d_{C^{\ell}}(F_{\t},F) \le C\subg \cdot q\subz
\delta_{n}^{s-\ell+3}\|\t\|\cdot \|\psi\|_{C^{\ell}}\qquad \mbox{for $\t\in R$ and $0\le \ell\le s$.}
\end{equation}
From this inequality in the case $\ell=0$,  we obtain
\begin{equation}\label{magn3}
d_{C^{0}}(F_{\t}^j,F^j) \le p\exp(\lmax p) \cdot 
C\subg \cdot q\subz \delta_n^{s+3}\exp(\chicl n)<\delta_n
\end{equation}
for $0\le j\le p$ and $\t\in
R$, where the second inequality follows from the condition  (\ref{eqnp}) and the
definition of 
$r\subz$ 
provided that $n$ is larger than some constant $N\sube$.
(Recall the notation introduced in section \ref{sec-dist}.)

Denote by
$\partial_{i}$  the partial differentiation  with respect to the parameter~$t_{i}$. Then
\begin{equation}\label{eqn-delta-bd0}
\|\partial_{i} F_{\t}(w)\|
\le  \delta_{n}^{s+3}
\end{equation}
and
\begin{equation}\label{eqn-delta-bd1}
\|\partial_{i} (DF_{\t})({\v})\|
\le 
C\subg \cdot \delta_{n}^{s+2}\|\v\|
\end{equation}
for  any  $w\in \man$, 
$\v\in
\cone^{u}(w)$ and $\t\in R$. If $d(w,x_i)<\delta_n/10$ in addition, we also have 
\begin{equation}\label{eqn-delta-bd2}
|\v^{*}(\partial_{i} (DF_{\t}(\v)))|\ge C\subg^{-1}\cdot\delta_{n}^{s+2}\|\v\|
\end{equation}
for  any  $w\in \man$, 
$\v\in
\cone^{u}(w)$ and $\t\in R$, where $\v^{*}$ is the unit cotangent vector at $F_{\t}(w)$
that is normal to $DF_{\t}(\v)$.

In the following argument, we assume that 
\begin{equation}\label{bdisj}
F^{j}(\ball(x_i, 2\delta_{n}))\cap
\ball(x_{i'}, 2
\delta_{n})=\emptyset
\end{equation}
for $0\le i,i'\le q\subz$ and $0\le 
j\le p$ but for the case where both $i= i'$ and $j=0$ hold. 
Note that (\ref{bdisj})  and  the estimate (\ref{magn3}) imply
\begin{equation}\label{disjbi}
F^{j}_{\t}(\ball(x_i, \delta_{n}))\cap
\ball(x_{i'},
\delta_{n})=\emptyset\qquad \mbox{for $\t\in R$.}
\end{equation}

Consider 
 a point $z\in \man$  and  families
of points
$y_{i}(\t)\in \man$,
$0\le i\le q\subz$, parameterized by
$\t\in R$ continuously. Suppose that it holds
\begin{itemize}
\item[(Y1)] $F_{\t}^{n}(y_i(\t))=z$,
\item[(Y2)] $y_i(\t)\in \Lambda(\chi,\epsilon,
(2h\subz+3)\epsilon n,n;F_{\t})$, and
\item[(Y3)] $d(F_{\t}^{n-p}(y_i(\t)), x_i)<\delta_{n}/10$
\end{itemize}
for $0\le i\le q\subz$ and
$\t\in R$. 
Let us put
\[
A_{i}(\t)=\delta_{n}^{s+2} \frac{|D^{*}F_{\t}^{p-1}(DF_{\t}^{n-p+1}(\eu(y_{i}(\t))))|}
{D_{*}F_{\t}^{p-1}(DF_{\t}^{n-p+1}(\eu(y_{i}(\t))))}
\]
for $1\le i\le q\subz$, where $\eu(z)$ is either of the two unit tangent vectors in $\E^u(z)$. Then we can show the
following estimates on the motion of the subspace $DF_{\t}^{n}(\E^{u}(y_i(\t)))$ as the parameter $\t$ moves. 
\begin{lemma}\label{lemma-diag} 
By letting the constant $N\sube$ larger if necessary, we have the following:
If $n\ge N\sube$, we have 
\[
C\subg^{-1}A_{i}(\t)\le |\partial_{i}\angle(DF_{\t}^{n}(\E^{u}(y_i(\t))),\E^{u}(z))|
\le C\subg A_{i}(\t)
\]
for $1\le i\le q\subz$, and also 
\[
|\partial_{j}\angle(DF_{\t}^{n}(\E^{u}(y_{i}(\t))),\E^{u}(z))|
\le C\subg \exp(-\lambda\subg  p)A_{i}(\t)
\]
for $0\le i\le q\subz$ and $1\le j\le q\subz$ provided $i\neq j$.
\end{lemma}
\begin{proof} Let $1\le i\le q\subz$ and $0\le j\le q\subz$.
For $0\le
m\le n$, we denote by
$\e_{m}$ the unit tangent vector in the direction of 
$DF_{\t}^{m}(\eu(y_{i}(\t)))$ and by
$\e^{*}_{m}$ the unit cotangent vector that is normal to $\e_{m}$. We can and do choose the orientation of 
the cotangent vectors $\e^{*}_{m}$ so that  $
(DF^{n-m})^{*}(\e^{*}_{n})=D^{*}F^{n-m}(\e_{m})\cdot \e^{*}_{m}$. Also we denote
$z_{m}=F_{\t}^{m}(y_{i}(\t))$ for simplicity.
Notice that $\e_m$, $\e^{*}_{m}$ and $z_m$  depend on the
parameter~$\t$. 

We first give simple consequences of the conditions (Y1) and (Y3). From 
 (\ref{disjbi}) and the condition (Y3), the point
 $z_{m}$ is not contained in $\ball(x_j,\delta_{n})$ for
$n-2p<m<n$ but for
the case where both $m=n-p$ and $j=i$ hold. Especially, the points $z_{m}$ for $n-p<m<n$ are not contained in $
\bigcup_{\ell=0}^{q\subz}\ball(x_\ell,\delta_{n})$. So  the condition (Y1)
implies that the point 
$z_m$ does not depend on the parameter $\t$. For   $0\le m\le n-p$, differentiation of
the both sides of the identity $F_{\t}^{n-p+1-m}(z_{m})\equiv z_{n-p+1}$  gives
\[
\left(DF_{\t}^{n-p+1-m}\right)_{z_m}(\partial_{j}z_m)
+\sum_{\ell=m+1}^{n-p+1}\left(DF_{\t}^{n-p+1-\ell}\right)_{z_{\ell}}
\left((\partial_{j}F_{\t})(z_{\ell-1})\right)=0.
\] 
Applying $(DF_{\t}^{n-p+1-m})_{z_m}^{-1}$ to the both sides of this identity and using
(\ref{eqn-delta-bd0}) and (\ref{eqn-c3}), we obtain
\begin{equation}\label{est-zj}
\left|\partial_{j}z_{m}\right|\le \sum_{\ell=m+1}^{n-p+1}
C\subg {\left\|((DF_{\t}^{\ell-m})_{z_m})^{-1}\right\|} \cdot \delta_{n}^{s+3}
\le \sum_{\ell=m+1}^{n-p+1}
\frac{C\subg\delta_{n}^{s+3}}{\left|D^{*}F_{\t}^{\ell-m}(\e_{m})\right|}.
\end{equation}

Now  we are going to  estimate
\[
\partial_{j}\angle(DF_{\t}^{n}(\E^{u}(y_{i}(\t))),\E^{u}(z))
=\partial_{j}\angle(\e_{n}, \E^{u}(z))=\frac{\e^{*}_{n}(
\partial_{j}(DF_{\t}^{n}(\e_{0})))}{D_{*}F_{\t}^{n}(\e_{0})}.
\]
Differentiating the both sides of 
\[
DF_{\t}^{n}(\e_{0})=(DF_{\t})_{z_{n-1}}\circ (DF_{\t})_{z_{n-2}}\circ
\cdots\circ (DF_{\t})_{z_{0}}(\e_0)
\]
 and using the relation $DF_{\t}^{m}(\e_{0})=D_{*}F_{\t}^{m}(\e_0)\cdot \e_{m}$, we can
obtain
\begin{align*}
\partial_{j}(DF_{\t}^{n}(\e_{0}))
&=\sum_{m=0}^{n-1}(DF_{\t}^{n-m-1})_{z_{m+1}}
((\partial_{j}(DF_{\t})_{z_m})(\e_{m}))\cdot  D_{*}F_{\t}^{m}(\e_{0})\\
&\qquad+
\sum_{m=0}^{n-1}(DF_{\t}^{n-m-1})_{z_{m+1}}(
D^{2}F_{\t}(\e_{m},\partial_{j}z_{m}))\cdot D_{*}F_{\t}^{m}(\e_{0})\\
&\qquad+
(DF_{\t}^{n})_{z_0}( D\eu(\partial_{j} z_{0})).
\end{align*}
From this and the relation $ (DF^{n-m})^{*}(\e^{*}_{n})=D^{*}F^{n-m}(\e_{m})\e^{*}_{m}$, 
it follows
\begin{align*}
\e^{*}_{n}(
\partial_{j}(DF_{\t}^{n}(\e_{0})))
&=
\sum_{m=0}^{n-1}
D^{*}F_{\t}^{n-m-1}(\e_{m+1})\cdot
\e^{*}_{m+1}((\partial_{j}(DF_{\t})_{z_m})(\e_{m}))  D_{*}F_{\t}^{m}(\e_{0})\\
&\quad+
\sum_{m=0}^{n-1}
D^{*}F_{\t}^{n-m-1}(\e_{m+1})\cdot
\e^{*}_{m+1}(D^{2}F_{\t}(\e_{m},\partial_{j}z_{m}))  D_{*}F_{\t}^{m}(\e_{0})\\
&\quad +D^{*}F_{\t}^{n}(\e_{0})
\cdot  \e^{*}_{0}(D\eu(\partial_{j} z_{0})).
\end{align*}
Note that  $(\partial_{j}(DF_{\t}))_{z_m}=0$ for $n-2p<m<n$ but for the case
$m=n-p$, and that $\partial_{j}z_{m}=0$ for $n-p<m\le n$ as we noted above.
Thus we obtain
\begin{align}
\frac{\e^{*}_{n}(
\partial_{j}(DF_{\t}^{n}(\e_{0})))}{D_{*}F_{\t}^{n}(\e_{0}) }&-
\frac{D^{*}F_{\t}^{p-1}(\e_{n-p+1})}
{D_{*}F_{\t}^{p-1}(\e_{n-p+1})}
\cdot
\frac{
\e^{*}_{n-p+1}((\partial_{j}(DF_{\t})_{z_{n-p}})(\e_{n-p}))}
{D_{*}F_{\t}(\e_{n-p})}\label{eqn-dq1}\\ &= 
\sum_{m=0}^{n-2p}
\frac{D^{*}F_{\t}^{n-m-1}(\e_{m+1})}
{D_{*}F_{\t}^{n-m-1}(\e_{m+1})}\cdot
\frac{
\e^{*}_{m+1}((\partial_{j}(DF_{\t})_{z_m})(\e_{m})))}
{D_{*}F_{\t}(\e_{m})}\notag\\
&\qquad +
\sum_{m=0}^{n-p}\frac{D^{*}F_{\t}^{n-m-1}(\e_{m+1})}
{D_{*}F_{\t}^{n-m-1}(\e_{m+1})}\cdot
\frac{\e^{*}_{m+1}(D^{2}F_{\t}(\e_{m},\partial_{j}z_{m}))}
{D_{*}F_{\t}(\e_{m})}\notag\\
& \qquad+\frac{D^{*}F_{\t}^{n}(\e_{0})}{D_{*}F_{\t}^{n}(\e_{0})}
\cdot \e^{*}_{0}(D\eu(\partial_{j} z_{0})).\notag
\end{align}
From (\ref{eqn-delta-bd1}), the first sum on the right hand side  
is bounded in absolute value by 
\begin{align*}
&C\subg \delta_{n}^{s+2}\cdot \frac{|D^{*}F_{\t}^{p-1}(\eu_{n-p+1})|}
{D_{*}F_{\t}^{p-1}(\eu_{n-p+1})}\sum_{m=0}^{n-2p}\exp(-\lambda\subg(
n-p+1-m+2c\subg ))\\
&\qquad \le
C\subg  \cdot  A_{i}(\t)\exp(-\lambda\subg
p).
\end{align*}
By the estimate (\ref{est-zj}) on $\partial_j z_m$  and the condition (Y2),  the second sum on the right hand
side is bounded in absolute value  by
\begin{align*}
C\subg \sum_{m=0}^{n-p}&\sum_{\ell=m+1}^{n-p+1}
\frac{|D^{*}F_{\t}^{n-m-1}(\e_{m+1})|}
{D_{*}F_{\t}^{n-m}(\e_{m})}
\frac{\delta_{n}^{s+3}}
{|D^{*}F_{\t}^{\ell-m}(\e_{m})|}\\
&= C\subg \sum_{m=0}^{n-p}\sum_{\ell=m+1}^{n-p+1}
\frac{|D^{*}F_{\t}^{n-\ell}(\e_{\ell})|}
{D_{*}F_{\t}^{n-\ell}(\e_{\ell})}
\frac{\delta_{n}^{s+3}}
{D_{*}F_{\t}^{\ell-m}(\e_{m})|D^{*}F_{\t}(\e_{m})|}\\
&<C\subg \delta_{n}^{s+3}\frac{|D^{*}F_{\t}^{p-1}(\e_{n-p+1})|}
{D_{*}F_{\t}^{p-1}(\e_{n-p+1})} \sum_{m=0}^{n-p}\sum_{\ell=m+1}^{n-p+1}
\frac{\exp(-\lambda\subg(n-p+1-m+2c\subg))}
{\exp(-(2h\subz+4)\epsilon n)}\\
&<C\subg A_{i}(\t) \cdot \delta_{n}\exp((2h\subz+4)\epsilon n)
<C\subg A_{i}(\t)\exp(-\lambda\subg p),
\end{align*}
where the last inequality follows from the definition of the constant
$r\subz$ and the condition (\ref{eqnp}) on $p$. 
Similarly, we can show that the last term on the right hand side is bounded by
\[
C\subg \sum_{\ell=1}^{n-p+1}
\frac{|D^{*}F_{\t}^{n}(\e_{0})|}
{D_{*}F_{\t}^{n}(\e_{0})}
\frac{\delta_{n}^{s+3}}
{|D^{*}F_{\t}^{\ell}(\e_{0})|}
<C\subg A_{i}(\t)\exp(-\lambda\subg p).
\]
From
(\ref{eqn-delta-bd1}) and (\ref{eqn-delta-bd2}), we have 
\[
C\subg^{-1} \delta_{n}^{s+2}<|\e^{*}_{n-p+1}(\partial_{j}DF_{\t}(\e_{n-p}))|
< C\subg\delta_{n}^{s+2}\quad \mbox{ if $j=i$}
\]
 and 
$\partial_{i}DF_{\t}(\e_{n-p})\equiv 0$ otherwise. Using these 
estimates in (\ref{eqn-dq1}), we can conclude the
lemma, by taking the constant $N\sube$ larger if necessary.
\end{proof}

Consider the mapping $
\Psi:R\to
\Real^{q\subz}$
defined by
\begin{equation}\label{defpsi}
\Psi(\t)=
\Big(\angle(DF_{\t}^{n}(\E^{u}(y_{i}(\t))),
DF_{\t}^{n}(\E^{u}(y_{0}(\t)))\Big)_{i=1}^{q\subz}.
\end{equation}
As a consequence of the lemma \ref{lemma-diag}, we have the following corollary, where we take the
constant $N\sube$ still larger if necessary.
\begin{corollary}\label{lem-jac}
The
mapping
$\Psi$ is injective and there is a constant $B\subz$ such that  
\[
|\det D\Psi(\t)|>\exp(-B\subz  \epsilon n)\quad \mbox{ for\/ $\t\in
R$,}
\] 
provided that $n\ge N\sube$.
\end{corollary}
\begin{proof} Let us denote by $D\Psi(\t)_{ij}$ the $(i,j)$-entry of the
representation matrix of
$D\Psi(\t)$ with respect to the standard coordinate on $\Real^{q\subz}$. 
Lemma
\ref{lemma-diag} tells that the diagonal entries satisfy
\[
C\subg^{-1} A_{i}(\t)<|D\Psi(\t)_{ii}|<C\subg A_{i}(\t)
\]
while the off-diagonal entries satisfy
\[
|D\Psi(\t)_{ij}|<C\subg\exp(-\lambda\subg p) A_{j}(\t),\quad \mbox{  $j\neq i$.}
\]
These imply that $\Psi$ is  injective on $R$ and  $|\det D
\Psi(\t)|$ is  bounded  from below by $\prod_{i=1}^{q\subz}C\subg
A_{i}(\t)$, provided that $n$ is larger than some constant $C\subz$. Therefore we have 
\[
|\det D
\Psi(\t)|> \Big(C\subg
\exp((\chicl-\chiuu)p-(4h\subz+6+(s+2)r\subz)\epsilon n)\Big)^{q\subz}
\]
 from the condition (Y2).
Using the condition (\ref{eqnp}), we obtain the corollary.
\end{proof}

\subsection{The proof of proposition \ref{ss12}}
In this subsection, we complete the proof
of theorem
\ref{th-ge-trans} by proving proposition \ref{ss12}. Let  
$G$ be a mapping in
$ C^{r}(M,\torus)$ and $d>0$ a positive number.
We consider a large integer  $n>N\sube$, an integer $p$ satisfying the condition
(\ref{eqnp}), and  a point $z$ in the lattice~$\lattice_{n}$. We put $\delta_{n}=\exp(-r\subz
\epsilon n)$ as in the last subsection. 
 Let
$S=\{x_i\}_{i=0}^{q\subz}$ be an ordered subset in  the lattice
$\lattice(\delta_{n}/40)$.  We denote by 
$\ss_{1}(\bchi,\epsilon, n, p, z;S)$ the set of mappings $F$ in 
$\ss_{1}(\bchi,\epsilon, n, p, z)$ such that the subset $\{w_{i}\}_{i=0}^{q\subz}$ in the definition can be taken
so that  
\begin{itemize}
\item[(S4)] $d(w_i, x_i)< \delta_n/20$\quad  for $0\le i\le q\subz$.
\end{itemize}
The subset $\ss_{1}(\bchi, \epsilon, n,p, z)$ is contained in
the union of $\ss_{1}(\bchi,\epsilon, n, p, z;S)$ over all 
ordered subsets $S=\{x_i\}_{i=0}^{q\subz}$ of the lattice $
\lattice(\delta_{n}/40)$. 
And the number of such ordered sets $S$ is bounded by
$(40\delta_{n}^{-1}+1)^{2(q\subz+1)}$. Therefore, in order to prove the
inequality in  proposition \ref{ss12}, it is enough to show 
\begin{equation}\label{claim-eq12}
\mom_s(\base_{G}^{-1}(\ss_{1}(\bchi,\epsilon, n, p, z;S))\cap
\disk^{s-3}(d))<\exp((2(\chicl-\chiul)-2r\subz (q\subz+2)\epsilon)
n)
\end{equation}
for sufficiently large $n$.

Take a mapping $F$ in $\ss_{1}(\bchi,\epsilon, n, p,
z;S)$ arbitrarily and consider the family of mappings
$F_\t$ defined for the ordered subset $S$ in the last subsection.
Note that the conditions
(\ref{bdisj})  and (\ref{disjbi}) follows from the conditions (S3) and (S4).  
Let $\Y$ be the set of  continuous mappings
\[
\y:R\to M\times M\times \dots \times M, \quad
\y(\t)=(y_i(\t))_{i=0}^{q\subz}
\]
  that satisfy the conditions (Y1),
(Y2) and (Y3) in the last subsection.  
A family 
$\y(\t)$ in $\Y$ is uniquely determined once $\y(0)$ is given because of  the conditions
(Y1) and (Y2).  Thus we
have 
\begin{align*}
\#\mathcal{Y}&\le (\#(\Lambda(\bchi,\epsilon, (2h\subz+3)\epsilon
n,n;F)\cap F^{-n}(z)))^{q\subz+1}\\
 &\le \kappa\sube\exp((\chiuu+\chi_{c}^{++}+7\epsilon+6(2h\subz+3)\epsilon)(q\subz+1)
 n)\\
 &\le \exp((\chiuu+\chi_{c}^{++})(q\subz+1)n+C\subz \epsilon 
 n)
\end{align*}
for sufficiently large $n$, from corollary
\ref{prop-back-card} and the condition (Y2).
 
For a family $\y\in \Y$, we denote by $Z(\y)$ the set of parameters $\t\in
R$ such that 
\begin{align*}
\angle(DF_{\t}^{n}(\E^{u}(y_i(\t))), DF_{\t}^{n}(&\E^{u}(y_{0}(\t))))\\
&\le
\exp((\chicu-\chiul+6\epsilon+h\subz(\chicdif+\chiudif+4\epsilon))n)
\end{align*}
for all $1\le i\le q\subz$.
Lemma \ref{lem-jac} implies that we have
\[
\leb(Z(\y))\le \exp(
(\chicu-\chiul+h\subz(\chicdif+\chiudif))q\subz  n+C\subz   \epsilon n)
\]
provided that $n\ge N\sube$.

Suppose that $F_{\s}$ belongs to $\ss_{1}(\chi,\epsilon,n,p,z;S)$ for a parameter $\s\in R$. Then
there are points $w_i\in F^{-p}_{\s}(z)$, $0\le i\le q\subz$, and subsets $E_i\subset  F^{-(n-p)}_{\s}(w_i)$,
$0\le i\le q\subz$, which satisfy the conditions (S1)-(S4) with $F$ replaced by
$F_{\s}$. Consider a combination $(y_i)_{i=0}^{q\subz}$ of points such that
$y_i\in E_{i}$ for
$0\le i\le q\subz$. 
From  (\ref{magn}), we can check that
\[
d_{C^1}(F_{\t},F_{\s})<\rho\sube\exp((\chicl-5\epsilon)n-3\cdot 2(h\subz+1)\epsilon n)\qquad
\mbox{ for any $\t\in R$}
\]
provided $n$ is sufficiently large. 
Thus, the condition (S1) and 
lemma \ref{prop-rho}, we can check that there  exists a unique
element $\y(\t)=(y_i(\y))_{i=0}^{q\subz}$ in $\Y$ such that $y_{i}(\s)=y_i$ for $0\le i\le
q\subz$. The condition (S2) implies that $\s$ belongs to the subset $Z(\y)$. Therefore, if
$F_{\s}$ belongs to 
$\ss_{1}(\chi,\epsilon,n,p,z;S)$,  the parameter $\s$  belongs to the
subset
$Z(\y)$ for at least
\[
\prod_{i=0}^{q\subz} \#E_{i}\ge \exp(( \chicl+\chiul-\chicdif-\chiudif-r\subz \epsilon)(q\subz+1) n)
\]
elements $\y$ in $\Y$. Now we arrive at the estimate
\begin{align*}
&\leb(\{\t\in R \mid F_{\s}\in \ss_{1}(\chi,\epsilon,n,p,z;S)\})\le  \frac{\sum_{Y\in \mathcal
Y}
\leb(Z(\y))}{\prod_{i=0}^{q\subz} \#E_{i}}\\ &\le
\frac{\exp(
((\chicu-\chiul +h\subz(\chicdif+\chiudif))
q\subz+(\chiuu+\chi_{c}^{++})(q\subz+1))n+ C\subz
  \epsilon n)}
{\exp((
\chicl+\chiul-\chicdif-\chiudif-r\subz\epsilon)(q\subz+1)n) }.
\end{align*}
Note that we have this estimate uniformly for the mappings $F$ in $\ss_{1}(\chi,\epsilon,n,p,z;S)$.  Put
$m=q\subz$, $T_{i}=\exp(\chicl n)$ and $\psi_i=\psi_{n,x_i}$ for
$1\le i\le q\subz$ in  lemma
\ref{lem-meas-est}. Then the  assumption 
 (\ref{eqn-small-norm})  holds provided that $n$ is sufficiently large.
The conclusion is 
\begin{align*}
\mom_s(&\base_{G}^{-1}(\ss_{1}(\bchi,\epsilon, n, p, z;S))\cap
\disk^{s-3}(d))\notag\\
&<2^{\q\subz+1} \exp\left((\chi_{c}^{++}-\chicl-\chiul+(h\subz+2)(\chicdif+\chiudif))q\subz
n \right)\\
&\qquad \times \exp((\chi_{c}^{++}-\chicl+\chicdif+2\chiudif+ C\subz \epsilon)n).
\end{align*}
Using the condition in  the choice of $q\subz$, we  
 obtain (\ref{claim-eq12}) for sufficiently large $n$, provided that we take 
sufficiently small $\epsilon\subz$.


\section{Genericity of the no flat contact condition}\label{sec-nonflat}

In this section, we consider the situation where the images of admissible curves under 
an iterate of a mapping $F\in \U$ have flat contacts with  
the curves in the critical set $\criticalset(F)$, and investigate whether we can 
  resolve all of such flat contacts by perturbations. Our goal is the proof of 
theorem~\ref{th-ge-nonflat}, which will be carried out in the last subsection. The key idea in the proof  
is that  the non-flatness of contacts between curves is easier to establish if we assume higher
differentiability. The reader should notice that the content and the notation in this section is independent of
those in the last two sections. 

\subsection{The jet spaces of curves}

We begin  with
formulating a sufficient condition for the no flat contact condition  in terms of
{\em  jet}. For an integer $1\le q\le r$ and a point
$z\in \man$,  let
$\Curve_{z}^q$ be the set of germs of
$C^{q}$curves
$\curve:(\Real,0)\to (\man,z)$ at $z$. Recall that we always assume the curves to be parameterized by
length.   Two germs
$\curve_1$ and
$\curve_2$ in
$\Curve_{z}^q$ are said to have  contact of order $q$ if 
$d(\curve_1(t),\curve_2(t))/|t|^{q}\to 0$ as $t\to 0$.
This is an equivalence relation on the space $\Curve_{z}^q$. 
The equivalence classes are called {\em  $q$-jets of curve} and the  quotient space
is denoted by
$\Jet^q\Curve_{z}$. 
Suppose that  a $q$-jet $\j$ of curve at $z\in M$ is represented by  $\curve\in \Curve_{z}^q$. Then the 
 the tangent vector $\frac{d}{dt}\curve
(0)\in T_z M$ at $z$ does not depend on the choice of the representative $\curve$, and neither do the
differentials 
$d^i\curve(0)$,
$2\le i\le q$, which are defined in subsection \ref{ss-admc}. Thus we put
\[
\j^{(0)}=z,\quad \j^{(1)}=\frac{d}{dt}\curve
(0)\quad\mbox{and} \quad \j^{(i)}=d^i\curve(0)\qquad\mbox{ for $2\le i\le q$}.
\]

The jet space of 
curves of order $q$ is the disjoint union  $\Jet^q\Curve:=\amalg_{z\in \man} \Jet^q\Curve_{z}$, 
 which is equipped  with  the distance defined by
\[
d_{\Jet}(\j_1,\j_2)=\max\left\{d(\j^{(0)}_1,\j^{(0)}_1), \angle(\j^{(1)}_1,\j^{(1)}_2),\max\left\{
|\j^{(i)}_1-\j^{(i)}_2|; 2\le i\le q\right\}\right\}.
\]
Then the
following mapping is a homeomorphism:
\[
\j\in \Jet^q\Curve\mapsto \left(\j^{(1)},(\j^{(i)})_{i=2}^{q}\right)
\in T^{1}\man\times\Real^{q-1}
\]
where $ T^{1}\man$ is the unit tangent bundle of $\man$.  
Each mapping $F\in \U$ acts naturally on the space $\Jet^{q}\Curve$. 
 We denote this
action simply by
\[
F:\Jet^{q}\Curve\to\Jet^{q}\Curve,\qquad [\curve]\mapsto [F_* \curve].
\]
For $2\le q<r$, let $\Jet^q\admc\subset \Jet^q\Curve$ be the compact subset of $q$-jets that 
are represented by germs of admissible curves.  Lemma
\ref{lem-admc-def} tells that $F^{n}( \Jet^q\admc)\subset
\Jet^q\admc$ for $n\ge n\subg$. 

For a $C^{q}$curve $\curve:I\to \man$ defined on an interval $I$, its $q$-jet extension 
 is the mapping $\Jet^{q}\curve:I\to \Jet^{q}\Curve$ that carries a parameter $t\in I$ to 
the jet  in $\Jet^q\Curve_{\curve(t)}$ that is 
 represented by the germ of
$\curve$ at
$t$. Recall that the critical set $\criticalset(F)$ for any 
mapping
$F$ in
$\U$ consists of finitely many $C^{r-1}$curves. Let 
$\Crit(F)\subset \Jet^{r-2}\Curve$ be the union of the images of their $(r-2)$-jet
extensions:
\[
\Crit(F)=\{\Jet^{r-2}\curve(I)\mid \mbox{$\curve:I\to M$ is a $C^{r-1}$curve contained in
$\criticalset(F)$.}\}
\]

\begin{lemma}\label{lemma-jetnf}If a mapping $F\in \U$ satisfies 
\begin{equation}\label{cond:jet}
F^{n}(\Jet^{r-2}\admc)\cap 
\Crit(F)=\emptyset\qquad \mbox{for some $n\ge 1$,}
\end{equation}
then $F$ satisfies the no flat contact condition. 
\end{lemma}
\begin{proof} 
For each point in $\criticalset(F)$, we can find a small 
$C^{r-1}$coordinate neighborhood $(U,\psi:U\to \Real^2)$ such that $\psi(\criticalset(F)\cap U)$
is  an  interval in the $x$-axis
$\Real\times\{0\}$ and that it holds
either 
\begin{itemize}
\item[(a)] $D\psi(\cone^c(z))$ contains the $x$-axis $\Real\times\{0\}$ for every $z\in U$, or 
\item[(b)] $D\psi(\cone^c(z))$ contains the $y$-axis $\{0\}\times\Real$ for every $z\in U$. 
\end{itemize}
Since the critical set $\criticalset(F)$ is compact, we can cover it  by finitely many coordinate
neighborhood with these properties. So, for the purpose of proving the lemma,  it is enough to
show the following claim for each coordinate neighborhood $(U,\psi)$ as above: For  a
constant
$C>0$ and
$n_0>0$, it holds
\[
\leb_{\Real}\left(\left\{t\in [0,a]\mid 
F^n(\curve(t))\in U\mbox{ and
}d(F^{n}(\curve(t)),\criticalset(F))<\epsilon\right\}\right)<C\cdot 
\epsilon^{1/r-2}\max\{a,1\}
\]
for any $a>0$, 
$\curve\in \admc(a)$, $n\ge n_{0}$, $\epsilon>0$.
If  the condition (a) above holds, this claim is clear
because the images of the admissible curves in $U$ by the mapping $\psi$ are curves whose slope is
uniformly bounded away from
$0$.   Thus it remains to check the claim above in the case where the condition
  (b) holds.  To this end, it is enough to show the following lemma, because, in the case (b), the images of the admissible curves  by  $\psi$ are graphs of $C^{r-1}$functions whose slopes are bounded by some constant~$C\subg$.
\begin{claim}\label{claim-local-curve}
If  a $C^{r-1}$function
$\varphi$ on a compact interval $I\subset \Real$ satisfies, 
\begin{align*}
&\max_{x\in I}\;\max\left\{ |{d^q \varphi}/{dx^{q}}(x)|\; ;\;  1\le q\le r-1\right\} \le K\quad\mbox{and}\\
&\min_{x\in I}\;\max\left\{|{d^q \varphi}/{dx^{q}}(x)|\; ;\; 1\le q\le
r-2\right\}> \rho\
\end{align*}
for some positive constants $K$
and
$\rho$,  then we have
\begin{equation*}
\leb_{\Real}\{x\in \Real\mid |\varphi(x)|<\epsilon\}<C(r,\rho,K,I)\cdot \epsilon^{1/(r-2)}
\quad \mbox{for any $\epsilon>0$,}
\end{equation*}
where $C(r,\rho,K,I)$ is a constant that depends
 only on $r$, $\rho$, $K$ and the length of $I$. 
\end{claim} 
\noindent
We  show this claim by using the following lemma\cite[Lemma 5.3]{BST}. 
\begin{lemma}\label{morse}
If a $C^{q}$function 
$h$ on an interval $J$ satisfies
$|d^{q}h/dx^q(x)| \geq \rho>0$
for all $x\in J$. 
Then $
\leb_{\Real}(\{x\in J| \ | h(x)|\leq
\varepsilon\})\le 2^{q+1}\left(\varepsilon/\rho\right)^{1/q}$
for any $\epsilon>0$. 
\end{lemma}
\begin{proof}[Proof of claim \ref{claim-local-curve}]
Let $X\subset I$ be the set of points $x\in I$ such that $|\varphi(x)|\le
\rho/2$. For
each point $x\in X$,  there is an integer 
$1\le m\le r-2$ such that $|d^{m}\varphi/dx^{m}(x)|>\rho$ and hence that
$|d^{m}\varphi/dx^{m}|\ge\rho/2$ on the interval $J(x):=(x-\rho/(2K), x+\rho/(2K))$.
We can take points  $x_i\in X$, $i=1,2,\dots,
i_0$, so that the intervals $J(x_i)$ cover the subset $X$ and that the intersection multiplicity is
$2$, so  $i_0\le (2\leb_{\Real}(I)/(\rho/K))+1$. Applying lemma \ref{morse} to  each interval
$J(x_i)$, we can see that  $\leb_{\Real}\{x\in \Real\mid |\varphi(x)|<\epsilon\}$ is bounded by
$ i_0 \cdot
2^{r-1}\left(\epsilon/(\rho/2)\right)^{1/(r-2)}$ provided $\epsilon<\rho$. This implies claim 
\ref{claim-local-curve}.
\end{proof}
We have finished the proof of lemma \ref{lemma-jetnf}.
\end{proof}

\subsection{Lattices in the jet space}
In this subsection, we consider lattices in the space of admissible jets
$\Jet^{r-2}\admc$ and formulate  a sufficient condition for the no flat contact
condition by using them. Henceforth, we fix integers  
$2<
\nu<r\le s$ satisfying the condition   (\ref{rs}). Note that the condition (\ref{rs})
can be written in the form
\[
(r-2)\left(r-1-\frac{r-3}{2}\right)
<(r-\nu-2)\left(r-3-\frac{2s-r-\nu+1}{2\nu}\right).
\]
Thus we can and do cover the interval
$[\lambda\subg/2, 2\lmax]$ by finitely many intervals
$I(\ell)=(\lambda^{-}(\ell),\lambda^{+}(\ell))$, $1\le \ell\le \ell_{0}$, such that 
$\lambda^{-}(\ell)>\lambda\subg/4$ and that 
\[
(r-2)\left(r-1-\frac{r-3}{2}\frac{\lambda^{-}(\ell)}{\lambda^{+}(\ell)}\right)
<(r-\nu-2)\left(r-3-\frac{2s-r-\nu+1}{2\nu}\right).
\] 
For $n\ge 1$ and $1\le \ell\le \ell_{0}$, let $\spl(n,\ell)$
be the set of  jets $\j$  in $\Jet^{r-2}\admc$ that satisfy
\begin{itemize}
\item[(Q1)] the point $\j^{(0)}$ is contained in  the lattice
$\lattice(\exp(-\lambda^{+}(\ell)(r-2)n ))$,
\item[(Q2)] the angle $\angle(\j^{(1)},\eu(\j^{(0)}))$ is a multiple of $
\exp(-\lambda^{+}(\ell)(r-3)n)$, and
\item[(Q3)] $\j^{(q)}$ is a multiple of
$\exp((-\lambda^{+}(\ell)(r-3)+\lambda^{-}(\ell)(q-1))n)$ for
$2\le q\le r-2$.
\end{itemize}
We have
\begin{equation}\label{card-q}
\#\spl(n,\ell)\le C\subg
\exp\left((r-2)\left((r-1){\lambda^{+}(\ell)}-\frac{r-3}{2}\lambda^{-}(\ell)\right)
n\right).
\end{equation}
For integers $n\ge 1$, $1\le \ell\le \ell_{0}$, a mapping $F\in \U$ and $\sigma=0,1$, we define 
$V_{\sigma}(n,\ell;F)$ as the set of jets $\j$ in $\Jet^{r-2}\admc$ 
that satisfy
\[
\exp(\lambda^{-}(\ell)n-\sigma)\le |D_{*}F^{n}(\j^{(1)})|\le
\exp(\lambda^{+}(\ell)n+\sigma).
\]
Then, from the choice of the numbers $\lambda^{\pm}(\ell)$, the subsets $V_{0}(n,\ell;F)$ for 
$1\le \ell\le \ell_{0}$ cover $\Jet^{r-2}\admc$ provided that $n$ is larger than some
constant
$C\subg$. 
\begin{lemma}\label{lemma-aprx}
There is a constant $B\subg>1$ such that, for any jet
$\j$ in $V_{0}(n,\ell;F)$ with $n\ge B\subg$ and $1\le \ell\le
\ell_{0}$, there exists a jet $\i\in \spl(n,\ell)\cap V_{1}(n,\ell;F)$ such that
\begin{equation}\label{claim-dist}
d_{\Jet}(F^{n}(\j),F^{n}(\i))<B\subg \exp(-\lambda^{+}(\ell) (r-3)n).
\end{equation}
\end{lemma}
\begin{proof}
Let us take a jet 
$\j\in V_{0}(n,\ell;F)$ arbitrarily.  Let  $w$ be the point in the lattice
$\lattice(\exp(-\lambda^{+}(\ell)(r-2)n ))$ that is closest to $\j^{(0)}$. As
$\j^{(1)}$ belongs to $ \cone^u(\j^{(0)})$,  the minimum angle between 
$\j^{(1)}$  and the cone
$\cone^{u}(w)$ is bounded by $C\subg \cdot d(\j^{(0)},w)$. Hence
 we can choose a jet $\i\in
\spl(n,\ell)$ such that
\begin{itemize}
\item[(I1)] $\i^{(0)}=w$ and hence
$d(\j^{(0)},\i^{(0)})<\exp(-\lambda^{+}(\ell)(r-2)n)$,
\item[(I2)] $\angle(\j^{(1)},\i^{(1)})<\exp(-\lambda^{+}(\ell)(r-3)n)+C\subg
\exp(-\lambda^{+}(\ell)(r-2)n) 
$ and 
\item[(I3)]
$|\j^{(q)}-\i^{(q)}|<\exp((-\lambda^{+}(\ell)(r-3)+\lambda^{-}(\ell)(q-1))n)$ for
$2\le q\le r-2$.
\end{itemize}
For $0\le m\le n$, we put
$z(m)=F^{m}(\j)^{(0)}=F^{m}(\j^{(0)})$, $w(m)=F^{m}(\i)^{(0)}=F^{m}(\i^{(0)})$  and 
\[
\Delta^{q}_{m}=
\begin{cases}
d(F^{m}(\j)^{(0)}, F^{m}(\i)^{(0)})=d(z(m),w(m)),&\mbox{for $q=0$};\\
\angle(F^{m}(\j)^{(1)},F^{m}(\i)^{(1)})=\angle(DF^{m}(\j^{(1)}), DF^{m}(\i^{(1)})),&\mbox{for $q=1$};\\
|F^{m}(\j)^{(q)}-F^{m}(\i)^{(q)}|,&\mbox{for $2\le q\le r-2$.}
\end{cases}
\]
In order to prove the inequality (\ref{claim-dist}), it is enough
to show 
\[
\Delta_{n}^{q}\le C\subg \exp(-\lambda^{+}(\ell) (r-3)n)\qquad \mbox{for $0\le q\le
r-2$.}
\]

First we prove 
\begin{equation}\label{delta0}
\Delta^{0}_{m}\le 2\|DF^{m}_{z(0)}\|\cdot \Delta^{0}_{0}\le C\subg
\exp(-\lambda^{+}(\ell) (r-3)n)
\end{equation}
for
$1\le m<n$. 
As $\j\in V_{0}(n,\ell;F)$, 
we have
\begin{align}\label{jetzero}
\|DF^{m-k}_{z(k)}\|&\le C\subg \cdot D_{*}F^{m-k}(DF^{k}(\j^{(1)}))\le \frac{C\subg\cdot D_{*}F^{n}(\j^{(1)})}
{D_{*}F^{n-m}(DF^{m}(\j^{(1)}))\cdot D_{*}F^{k}(\j^{(1)})}\\
&\le C\subg
\exp(\lambda^{+}(\ell)n-\lambda\subg(n-m+k))\notag 
\end{align}
for $0\le k\le  m\le n$.  So the second inequality in (\ref{delta0}) follows from the condition~(I1). 
We prove the first inequality in (\ref{delta0}) by induction on $1\le m<n$. Using the
simple estimate
\[
\left\|
\exp_{z(m)}^{-1}(w(m))-DF_{z(m-1)}(\exp_{z(m-1)}^{-1}(w(m-1)))
\right\|
\le C\subg (\Delta_{m-1}^{0})^2
\]
repeatedly, we can get the following inequality for
$\Delta_{m}^{0}=\|\exp_{z(m)}^{-1}(w(m))\|$:
\begin{align}\label{inejetdel}
\Delta_{m}^{0}\le 
\|DF_{z(0)}^{m}\|\Delta_{0}^{0}+C\subg\sum_{k=0}^{m-1}
\|DF_{z(k+1)}^{m-k-1}\|(\Delta_{k}^{0})^2\qquad \mbox{for $0\le m\le n$.}
\end{align}
Note that we have, from (\ref{eqn-c2}),
\begin{align}\label{jetzero1}
\|DF_{z(k+1)}^{m-k-1}\|\|DF^k_{z(0)}\|&\le C\subg D_{*}F^{m-k-1}(DF^{k+1}(\eu(z_0)))\cdot 
D_{*}F^{k}(\eu(z_0))\\
&\le C\subg  \frac{D_{*}F^{m}(\eu(z_0))}
{D_{*}F(DF^{k}(\eu(z_0)))}
\le C\subg
\|DF^{m}_{z(0)}\|\notag
\end{align}
for $0\le k\le m-1$. 
Consider an integer $0\le m_0\le n$ and suppose that the left inequality in (\ref{delta0}) holds
for $0\le m< m_0$. Then, using them and the estimates (\ref{jetzero}) and (\ref{jetzero1}) in
(\ref{inejetdel}),  we can obtain 
\begin{align*}
\Delta_{m_0}^{0}&\le 
\|DF_{z(0)}^{m_0}\|\Delta_{0}^{0}+C\subg\sum_{k=0}^{m_0-1}
\|DF_{z(k+1)}^{m_0-k-1}\|\cdot 2\|DF_{z(0)}^{k}\|\Delta_{0}^{0}\cdot \Delta_{k}^{0}\\
&\le \|DF_{z(0)}^{m_0}\| \Delta_{0}^{0}
\left(
1+C\subg\cdot  n \cdot \exp(-\lambda^{+}(\ell)(r-3)n)
\right). 
\end{align*}
This implies the first inequality in (\ref{delta0}) for $m=m_0$, provided that  $n$ is larger
than some constant~$C\subg$. Thus we can obtain (\ref{delta0}) for $1\le m\le n$ by induction.

Next we estimate  $\Delta_{m}^{1}$ for $0\le m\le n$. We have
\begin{align*}
\Delta_{m}^{1}&\le \angle(DF^m_{z(0)}(\j^{(1)}), DF^m_{z(0)}(\i^{(1)}))\\
&\qquad\qquad +
\sum_{k=0}^{m-1}
\angle(DF_{z(k)}^{m-k}(DF^k_{w(0)}(\i^{(1)})),
DF_{z(k+1)}^{m-k-1}(DF^{k+1}_{w(0)}(\i^{(1)})))
\end{align*}
where we omit the operations of parallel translation. (See the remark given in the proof of
lemma
\ref{prop-rho}.)
For $0\le k<n$, we have $DF^{k}_{w(0)}(\i^{(1)})\in \cone^u(w(k))$ and 
$
d(z(k),w(k))=\Delta_{k}^{0}\le C\subg
\exp(-\lambda^{+}(\ell) (r-3)n)$. Hence  the parallel
translation of $DF^{k}_{w(0)}(\i^{(1)})$  to
$z(k)$ does not belong to the central cone
$\cone^c(z(k))$ provided that  
$n$ is larger than some constant~$C\subg$. Using (\ref{anglerel}), we can obtain
\begin{align*}
&\angle(DF_{z(k)}^{m-k}(DF^k_{w(0)}(\i^{(1)})),
DF_{z(k+1)}^{m-k-1}(DF^{k+1}_{w(0)}(\i^{(1)})))\\
& \le A\subg
\frac{|D^{*}F^{m-k-1}(DF^{k+1}(\j^{(1)}))|}{D_{*}F^{m-k-1}(DF^{k+1}(\j^{(1)}))}
\angle(DF_{z(k)}(DF^k_{w(0)}(\i^{(1)})),
DF_{w(k)}(DF^{k}_{w(0)}(\i^{(1)})))\\
& \le C\subg \exp(-\lambda\subg(m-k-1))\Delta_{k}^{0}<
C\subg
\exp(-\lambda\subg(m-k-1)-\lambda^{+}(\ell) (r-3)n).
\end{align*}
Likewise we can obtain $
\angle(DF_{z(0)}^{m}(\j^{(1)}), DF_{z(0)}^{m}(\i^{(1)}))\le C\subg
\exp(-\lambda\subg m)\Delta_{0}^{1}$. 
Therefore, by condition (I2), we conclude 
\begin{align*}
\Delta_{m}^{1}&\le  C\subg
\exp(-\lambda\subg m)\Delta_{0}^{1}+\sum_{k=0}^{m-1}C\subg \exp(-\lambda\subg(m-k-1)-\lambda^{+}(\ell)
(r-3)n)\\ 
&\le C\subg
\exp(-\lambda^{+}(\ell) (r-3)n).
\end{align*}

Finally, we estimate $\Delta_{n}^{q}$ for $2\le q\le r$.  From the formula
(\ref{eqn-inductive}), we can see 
\begin{equation}\label{eqn-rec}
\Delta^{q}_{m}\le
\frac{|D^{*}F(DF^{m-1}(\j^{(1)}))|}{D_{*}F(DF^{m-1}(\j^{(1)}))^{q}}\Delta^{q}_{m-1}+C\subg\sum_{0\le
d<q}\Delta^{d}_{m-1}.
\end{equation}
Consider this inequality (\ref{eqn-rec}) for $m=n$ and estimate the right hand side   by using
(\ref{eqn-rec}) recurrently as long as there exist  terms
$\Delta^{q}_{m}$  with $q>1$ or $m>0$ on the right hand side. Then we see that 
$\Delta^{q}_{n}$ is bounded by 
\begin{equation}\label{bigea}
\begin{aligned}
&\frac{|D^{*}F^{n}(\j^{(1)})|}{D_{*}F^{n}(\j^{(1)})^{q}}\Delta^{q}_{0}\\
&+
C\subg \sum_{1< d<q}\;
\sum_{
\substack{0= n_{d}\le n_{d+1}\le \cdots\\\qquad\cdots\le n_q<  n_{q+1}=n+1}}\;
\prod_{\ell=d}^{q}\;\prod_{n_{\ell}\le j< n_{\ell+1}-1}\;
\frac{|D^{*}F(F^{j}(\j)^{(1)})|}{D_{*}F(F^{j}(\j)^{(1)})^\ell}\Delta^{d}_{0}\\
 &
+C\subg\sum_{\substack{d=0,1\\0\le m<
n}}\;
\sum_{\substack{m= n_{d}\le n_{d+1}\le \cdots\\\qquad\cdots\le n_q<  n_{q+1}=n+1}}\;
\prod_{\ell=d}^{q}\;\prod_{n_{\ell}\le j< n_{\ell+1}-1}\;
\frac{|D^{*}F(F^{j}(\j)^{(1)})|}{D_{*}F(F^{j}(\j)^{(1)})^\ell}\Delta^{d}_{m}.
\end{aligned}
\end{equation}
Note that, for any sequence $m= n_{d}\le n_{d+1}\le \cdots\cdots\le n_q<  n_{q+1}=n+1$ with $q\le r$, we have
\begin{align*}
\prod_{\ell=d}^{q}\;\prod_{n_{\ell}\le j< n_{\ell+1}-1}\;
\frac{|D^{*}F(F^{j}(\j)^{(1)})|}{D_{*}F(F^{j}(\j)^{(1)})^\ell}
&\le \frac{\exp(-\lambda\subg(n-m-q)+c\subg q)}{D_{*}F^{n-m}(F^{m}(\j)^{(1)})^{d-1}\exp(- q^2 \Lambda\subg)}\\
&\le C\subg \frac{\exp(-\lambda\subg(n-m))}{D_{*}F^{n-m}(F^{m}(\j)^{(1)})^{d-1}}.
\end{align*}
Hence  it follows from (\ref{bigea}) 
\begin{align*}
\Delta^{q}_{n}\le&\frac{\exp(-\lambda\subg n )}{D_{*}F^{n}(\j^{(1)})^{q-1}}\Delta^{q}_{0}+
C\subg n^q \sum_{1< d<q}\frac{\exp(-\lambda\subg n)}
{D_{*}F^{n}(\j^{(1)})^{d-1}}\Delta^{d}_{0}\\
&\qquad +C\subg\sum_{0\le m< n}(n-m)^q \exp(-\lambda\subg (n-m))
(\Delta^{0}_{m}+\Delta^{1}_{m})\\
\le & C\subg\max_{0\le m<n}(\Delta^{0}_{m}+\Delta^{1}_{m})+C\subg \sum_{1< d\le
q}\exp(-(d-1)\lambda^{-}(\ell)n)\Delta_{0}^{d},
\end{align*}
where the second inequality follows from the fact that the jet $\j$ belongs to
$V_{0}(n,\ell;F)$. 
Using the estimates on
$\Delta_m^0$ and $\Delta_m^1$ and the condition (I3) in the inequality above, we can conclude
\[
\Delta^{q}_{n}<C\subg \exp(-\lambda^{+}(\ell)(r-3)n)
\quad\mbox{for $2\le q\le r-2$.}
\]
We have proved the inequality (\ref{claim-dist}). The jet  $\i$ belongs to $V_{1}(n,\ell;F)$
because 
\[
\log \frac{D_{*}F^{n}(\e^{u}(\i^{(0)}))}{D_{*}F^{n}(\e^{u}(\j^{(0)}))}
\le C\subg \!\!\sum_{m=0}^{n-1}
(\Delta_{m}^{0}+\Delta_{m}^{1})
\le C\subg n \exp(-\lambda^{+}(\ell)(r-3)n)<1,
\]
 provided that $n$ is larger than some constant
$C\subg$. 
\end{proof}

For integers $n\ge 1$, $1\le \ell\le \ell_{0}$ and a 
jet
$\j\in \spl(n, \ell)$, let $\ss_2(n,\ell,\j)$ be 
  the set of mappings $F\in \U$ such that $\j\in V_{1}
(n,\ell; F)$ and that
\[
d_{\Jet}(F^{n}(\j), \Crit(F))\le 2 B\subg \exp(-\lambda^{+}(\ell)n (r-3)).
\]
Then the last lemma implies
\begin{corollary}\label{cirqv}
If there exists $n\ge B\subg$ such that $F\notin \ss_{2}(n,\ell,\j)$
for all $1\le \ell\le \ell_{0}$ and  $\j\in \spl(n,\ell)$, then 
 $F$ satisfies the no flat contact condition.
\end{corollary}
In the remaining part of this section, we shall estimate the measure of the subsets $\ss_{2}(n,\ell,\j)$ for
$\j\in
\spl(n,\ell)$ by using lemma \ref{lem-meas-est}.

\subsection{Perturbations}\label{just}
In this subsection,  we introduce some families of mappings  and
give a few estimates on the variation of the images of jets under the iterates of  mappings in the families.
 In the argument below, we fix
$1\le
\ell\le
\ell_0$ and put
\[
\delta_{n}=\exp(-\lambda^{+}(\ell)n/\nu)\quad \mbox{ for $n\ge 1$.}
\]
For $1\le q\le
r-2$, we take and fix a $C^{\infty}$function 
$\psi_{q}:\Real^{2}\to \Real$ such that 
\[
\psi_{q}(x,y)=
\begin{cases}
x^{q}/q!, &\mbox{for  $(x,y)\in \ball(0,1/10)$};\\
0,&\mbox{for  $(x,y)\notin \ball(0,1)$.}
\end{cases}
\]
\begin{remark} We can and do take the functions $\psi_{q}$ so that their $C^r$norm is bounded by some constant $C\subg$. 
\end{remark}
 For each point $\zeta\in M$,
we consider an isometric embedding 
\[
\varphi_{\zeta}:\{w\in \Real^2\mid
\|w\|<1\}\to
\torus
\] 
that carries the origin to
the point $\zeta$ and the $x$-axis $\Real\times\{0\}$ to $\E^{u}(\zeta)$.

Recall that we took the
subset $\U$ of mappings  as a neighborhood of a
$C^r$mapping $\bF$ in subsection \ref{cnbl}. 
For positive integers $n$, $1\le q\le r-2$ and a point $\zeta$  in $M$, we define a $C^{\infty}$
mapping
$\psi_{q,n,\zeta}:\man\to
\Real^{2}$ by
\[
\psi_{q,n, \zeta}(z)=
\begin{cases}
\delta_{n}^{s}\cdot\psi_{q}(\varphi_{\zeta}^{-1}(z)/\delta_{n})\cdot 
\e^{c}(\bF(\zeta))&\mbox{if $d(z,\zeta)<\delta_n$};\\
0,&\mbox{otherwise,}
\end{cases}
\]
where  $\e^{c}(\cdot)$ is either of the two unit tangent vectors in  the central subspace $\E^c(\cdot)$. Note
that, for any mapping
$F\in
\U$, the parallel translation of the vector
$\e^{c}(\bF(z))$ to $F(z)$ is contained in 
$\cone^c(F(z))$  from the choice of 
 $\U$  in subsection \ref{cnbl}.

For a positive  integer~$n$, a mapping
$F\in
\U$ and a point
$\zeta$ in $M$, we define
\begin{equation}\label{deffam2}
F_{\t}(z)=F(z)+\sum_{q=\nu+1}^{r-2} t_{q}\psi_{q,n,\zeta}(z):M\to \torus
\end{equation}
where $\t=(t_{\nu+1},t_{\nu+2},\cdots,t_{r-2})$ is the  parameter that
ranges over  $R=[-1,1]^{r-2-\nu}$. 
This the family of mappings that we are going to consider. 
\begin{remark}
The purpose of considering the family above is  to move the the images $F_{\t}^n(\j)$ of the jets $\j\in \spl(n,\ell)$ by choosing the point $\zeta$ appropriately. As it will turn out, we can keep control of the coordinates $F_{\t}^n(\j)^{(q)}$ with  $ q\ge \nu+1$ but {\em not} of those with $0\le q\le \nu$. 
This is the reason why we restricted the range of $q$ in between $\nu+1$ and $r-2$ in  (\ref{deffam2}). Note that, if we take smaller  $\nu$, we can keep control of more coordinates but the magnitude of the perturbation  becomes smaller. Thus, we have to choose a good value for $\nu$. The inequality (\ref{rs}) is related to this choice. 
\end{remark}

 Obviously we have
\begin{equation}\label{diff-F2}
d_{C^{q}}(F_{\t},F)\le C\subg\delta_{n}^{s-q}
\qquad\mbox{and}\qquad
\|\partial_{\t}F_{\t}\|_{C^{q}}\le C\subg\delta_{n}^{s-q}
\end{equation}
 for $0\le q\le r$ and $\t\in R$. Especially, $F_{\t}(M)\subset M$ if $n$ is sufficiently large.

We consider a jet 
$\j\in \spl(n,\ell)\cap V_{1}(n,\ell;F)$ and give some estimates on the
variation of  the image~$F_{\t}^{n}(\j)$.
 We begin with the estimate on the position
$F_{\t}^{n}(\j)^{(0)}$. 
\begin{lemma}\label{lem:fdif} We have, for $0\le m\le n$ and $\t\in R$, 
\begin{equation*}\label{diff-F}
d(F_{\t}^{m}(\j^{(0)}),F^{m}(\j^{(0)}))<C\subg \|DF^{m}_{\j^{(0)}}\|\delta_{n}^{s}\le 
C\subg \delta_{n}^{s-\nu} 
\end{equation*}
and
\begin{equation*}\label{est-del-F}
\|\partial_{\t}F_{\t}^{m}(\j^{(0)})\|<C\subg \|DF^{m}_{\j^{(0)}}\|\delta_{n}^{s} \le C\subg
\delta_{n}^{s-\nu}
\end{equation*}
provided that $n$ is larger than some constant $C\subg$. 
\end{lemma}
\begin{proof} The following argument  is a  modification of  that
in the former part of the  proof of lemma~\ref{lemma-aprx}. We denote 
$z(m)=F^{m}(\j^{(0)})$,
$w(m)=F_{\t}^{m}(\j^{(0)})$ and
$\Delta_{m}=d(z(m),w(m))$ for $0\le m\le n$, so $\Delta_{0}=0$.  
 Using the simple estimate
\[
\left\|
\exp_{z(m)}^{-1}(w(m))-(DF)_{z(m-1)}(\exp_{z(m-1)}^{-1}(w(m-1)))
\right\|
\le C\subg(\delta_{n}^{s}+ (\Delta_{m-1})^2)
\]
repeatedly, we  obtain
\begin{equation}\label{eqn-jco}
\Delta_{m}\le 
\sum_{k=0}^{m-1}
\|(DF^{m-k-1})_{z(k+1)}\|\cdot C\subg(\delta_{n}^{s}+ (\Delta_{k})^2)
\end{equation}
for $0\le m\le n$.
Consider an integer  $0\le m_0\le n$  and a positive number $K$. And suppose that we have
\begin{equation}\label{inedelt}
\Delta_{m}<K \|(DF^{m})_{z(0)}\|\cdot\delta_{n}^{s} 
\end{equation}
for $0\le m<m_0$.
 Then,  using this, the inequality (\ref{jetzero1})   and the simple estimate
 \[
C\subg^{-1} \exp(\lambda\subg k)\le \|DF^{k}_{z(0)}\|\le C\subg \|DF^{n}_{z(0)}\|\le C\subg\delta_{n}^{-\nu}\quad\mbox{ for $0\le k\le m\le n$}
 \]
 in the right hand side of the inequality  (\ref{eqn-jco}) for $m=m_0$, we obtain
\begin{align*}
\Delta_{m_0}&\le  C\subg
\|(DF^{m_0})_{z(0)}\|\cdot  \sum_{k=0}^{m-1}
\left(\delta_{n}^{s}\|DF^{k}_{z(0)}\|^{-1}+K^2 \delta_{n}^{2s}\|DF^{k}_{z(0)}\|\right)
 \\
 &\le  C\subg
\|(DF^{m_0})_{z(0)}\|\cdot \delta_{n}^{s}\cdot  \sum_{k=0}^{m-1}
\left(
\exp(-\lambda\subg k)+K^2 \delta_{n}^{s-\nu}\right).
\end{align*}
This implies  (\ref{inedelt}) for $m=m_0$, provided 
$K$ and $n$ are  larger than some constant~$C\subg$. Thus we can obtain the first
claim of the lemma by induction on $m$.

Put $\Delta'_{m}=\partial_{\t}F_{\t}^{m}(\j^{(0)})$ for $0\le m\le n$. Using the simple estimate
\[
\|\Delta'_{m}-(DF)_{z(m-1)}\Delta'_{m-1}\|\le C\subg (\delta_{n}^{s}+\Delta_{m-1}\|\Delta'_{m-1}\|)
\]
repeatedly, we can obtain
\begin{align*}
\Delta'_{m}\le 
\sum_{k=0}^{m-1}
\|(DF^{m-k-1})_{z(k+1)}\|\cdot C\subg(\delta_{n}^{s}+\Delta_{k}\|\Delta'_{k}\|
).
\end{align*}
From this and the estimates on $\Delta_m$ we have proved above, we can obtain the second
claim of the lemma by induction on $m$, in a similar manner as above.
\end{proof}

Next we give the estimates on $\partial_{\t}F^{n}_{\t}(\j)^{(q)}$ for $1\le q\le r-2$. 
We denote by  $\partial_{p}$ the
differentiation by the parameter 
$t_{p}$. For integers $p$ and $q$ satisfying $\nu+1\le p\le r-2$ and $1\le q\le r-2$, and for  a
jet
$\i\in
\Jet^{r-2}\admc$ and $\t\in R$, we define 
\[
\beta_{p}^{(q)}(\i,\t)=\pm
\frac{\sin(\angle(\ec(\bF(z)),F_{\t}(\i^{(1)}))\cdot 
\|\partial_{p}((D^{q}F_{\t})_{\i^{(0)}}(\i^{(1)},\i^{(1)},\cdots,\i^{(1)}))\|}
{D_{*}F_{\t}(\i^{(1)})^{q}}
\]
where  $(D^{q}F_{\t})_{z}: \otimes^q
T_{z}M
\to T_{F(z)}M$ is the $q$-th differential of $F_{\t}$ at
$z$ and the sign on the right hand side will be chosen appropriately in the argument below.
We have
\begin{equation}\label{estbeta} |\beta_{p}^{(q)}(\i,\t)|\le C\subg \delta_{n}^{s-q}.
\end{equation}
\begin{lemma} \label{est-ent}
There exists a positive constant $C\subg$ such that, if $n\ge C\subg$, it holds 
\begin{align*}
&\left|\partial_{p}(F_{\t}^{m}(\j)^{(q)})-\sum_{k=0}^{m-1}
\frac{D^{*}F_{\t}^{m-k-1}(F_{\t}^{k+1}(\j)^{(1)})}{D_{*}F_{\t}^{m-k-1}(F_{\t}^{k+1}(\j)^{(1)})^{q}}
\cdot \beta^{(q)}_{p}(F_{\t}^{k}(\j),\t)\right|
<C\subg\delta_{n}^{s-q+1}
\end{align*}
for any $\nu+1 \le q\le r-2$,  $\nu+1\le p\le r-2$, $\t\in R$ and  $0\le
m\le n$, provided that we choose the sign in the definition of $\beta^{(q)}_{p}(F_{\t}^{k}(\j),\t)$ appropriately.
\end{lemma}
\begin{proof} Take and fix   $\nu+1\le p\le r-2$ arbitrarily. For $0\le q\le r-2$ and  $0\le
m\le n$, we put
\[
\Delta_{m}^{(q)}=
\begin{cases}
\|\partial_{p} F^{m}_{\t}(\j)^{(0)}\|,&\mbox{if $q=0$};\\
\partial_{p}(\angle(F^{m}_{\t}(\j)^{(1)}, v_0)),&\mbox{if $q=1$};\\
\partial_{p}(F^{m}_{\t}(\j)^{(q)}),&\mbox{if $q\ge 2$}
\end{cases}
\]
where $v_0$ is some fixed vector. 
For $0< m\le n$, it holds
\[
\left|\Delta_{m}^{(1)}-\frac{D^{*}F_{\t}(F_{\t}^{m-1}(\j)^{(1)})
}
{D_{*}F_{\t}(F_{\t}^{m-1}(\j)^{(1)})}\Delta_{m-1}^{(1)}
-\beta_{p}^{(1)}(F_{\t}^{m-1}(\j),\t)\right|
\le C\subg \Delta_{m-1}^{(0)}\le C\subg\delta_{n}^{s-\nu}
\]
where the second inequality follows from lemma \ref{est-del-F}.
From this inequality and the estimate (\ref{estbeta}) for $q=1$, we can see
\begin{align*}
|\Delta_{m}^{(1)}|&\le C\subg 
\sum_{k=0}^{m-1}\frac{|D^{*}F_{\t}^{m-k}(F_{\t}^{k}(\j)^{(1)})|
}{D_{*}F_{\t}^{m-k}(F_{\t}^{k}(\j)^{(1)})}\left(\beta_{p}^{(1)}
(F_{\t}^{k}(\j),\t)+\delta_{n}^{s-\nu}\right)< C\subg
\delta_{n}^{s-\nu}
\end{align*}
for  $0\le
m\le n$. Recall the formula (\ref{eqn-inductive}) and the remark after it. By differentiating
the both sides of  (\ref{eqn-inductive}) with $F$ replaced by $F_\t$ and using (\ref{diff-F2}), 
we can obtain
\begin{equation}\label{eqnimp}
\left|\Delta_{m}^{(q)}-\frac{D^{*}F_{\t}(F_{\t}^{m-1}(\j)^{(1)})
}
{D_{*}F_{\t}(F_{\t}^{m-1}(\j)^{(1)})^{q}}\Delta_{m-1}^{(q)}
-\beta_{p}^{(q)}(F_{\t}^{m}(\j),\t))\right|
\le C\subg \delta_{n}^{s-q+1}+C\subg \sum_{0\le d< q}\Delta_{m-1}^{(d)}
\end{equation}
for $2\le q\le r-2$ and  $0\le
m\le n$.  
 Especially we have, from  (\ref{estbeta}),
\begin{equation*}
\left|\Delta_{m}^{(q)}-\frac{D^{*}F_{\t}(F_{\t}^{m-1}(\j)^{(1)})}
{D_{*}F_{\t}(F_{\t}^{m-1}(\j)^{(1)})^{q}}\Delta_{m-1}^{(q)}
\right|
\le C\subg \delta_{n}^{s-q+1}+C\subg \sum_{0\le d< q}\Delta_{m-1}^{(d)}
\end{equation*}
for $2\le q\le r-2$ and  $0\le
m\le n$. 
 Using this inequality repeatedly, we reach 
 \begin{align*}
|\Delta_{m}^{(q)}|&\le C\subg 
\sum_{k=1}^{m}\left(\frac{|D^{*}F_{\t}^{m-k}(F_{\t}^{k}(\j)^{(1)})|
}{D_{*}F_{\t}^{m-k}(F_{\t}^{k}(\j)^{(1)})^{q}}\left(
\delta_{n}^{s-q}+\sum_{0\le d<
q}\Delta_{k-1}^{(d)}\right)\right)\\
&\le C\subg\left(
\delta_{n}^{s-q}+\max_{0\le d<
q}\max_{0<k<m}\Delta_{k}^{(d)}\right) .
\end{align*}
Hence, we can show the estimate $|\Delta_{m}^{(q)}|\le C\subg\delta_{n}^{s-\nu}$ for $2\le q\le \nu$ and  $0\le
m\le n$  by induction
on
$q$,  by lemma \ref{lem:fdif} and the estimate on $|\Delta_{m}^{(1)}|$ above.
 
Next, using the inequality (\ref{eqnimp}) repeatedly, we can see that the left hand side of the inequality in the lemma is bounded by 
\[
 C\subg
\sum_{k=1}^{m}\frac{|D^{*}F_{\t}^{m-k}(F_{\t}^{k}(\j)^{(1)})|}{D_{*}
F_{\t}^{m-k}(F_{\t}^{k}(\j)^{(1)})^q}
\left(\delta_{n}^{s-q+1}+\sum_{0\le d<q}\Delta_{k-1}^{(d)}\right).
\]
By induction on $\nu+1\le q\le r-2$, we obtain the inequality  in the
lemma. 
\end{proof}

Note that, for any $C^r$ mapping $G:M\to M$ such that $d_{C^r}(G, \bF)<2\rho\subg$, the level curves of the
function
$\det G: z\mapsto
\det DG_{z}$ are regular in the
neighborhood $\ball(\criticalset(G),\rho\subg)$ of the critical set
$\criticalset(G)$, from the  choice of the constant $\rho\subg$ in subsection \ref{cnbl}.
For a point $w\in\ball(\criticalset(G),\rho\subg)$, we denote by $\c(w;G)$ 
the $(r-2)$-jet at $w$ that is given by the level curve passing through $w$. 

Suppose that a jet $\j\in \Jet^{r-2}\admc$ satisfies, for all
$\t\in R$,
\begin{itemize}
\item[{\rm (V1)}] $d(F_{\t}^{n-1}(\j)^{(0)},\zeta)<  \delta_{n}/10$,
\item[{\rm (V2)}] $d(F_{\t}^{n-1}(\j)^{(0)},\criticalset(F_{\t}))>3 \delta_{n}$, and
\item[{\rm (V3)}] $d(F_{\t}^{n}(\j)^{(0)},\criticalset(F_{\t}))< \delta_{n}$.
\end{itemize}
From the condition (V3), we can define the mapping $\Psi:R\to \Real^{r-\nu-2}$
  by
\[
\Psi(\t)= 
\left(\frac{F_{\t}^{n}(\j)^{(q)}-\c(F^{n}_{\t}(\j)^{(0)};F_{\t})^{(q)}}
{\delta_{n}^{s-q}}\right)_{q=\nu+1}^{r-2}
\]
provided that $n$ is so large that $\delta_n<\rho\subg$. The following is the goal of this subsection.

\begin{lemma}\label{lem-psiest} 
If the conditions (V1),(V2) and (V3) hold for all $\t\in R$, the mapping
$\Psi$ is a diffeomorphism and 
$|\det D\Psi(\t)|$ is bounded from below by a constant $C\subg^{-1}$, provided that $n$ is larger
than some constant $C\subg$. 
\end{lemma}
\begin{proof} From the condition (V1) and the definition of the family $F_\t$, we
have
\begin{alignat*}{2}
\beta^{(q)}_{p}(F_{\t}^{n-1}(\j),\t)&=0 & &\quad \mbox{ for $q>p$, and }\\
|\beta^{(q)}_{p}(F_{\t}^{n-1}(\j),\t)|&\ge C\subg ^{-1}\delta_{n}^{s-q} &
 &\quad \mbox{ for $q=p$,
}
\end{alignat*}
in addition to (\ref{estbeta}).
We show 
\begin{equation}\label{est-bq}
\left| \sum_{k=0}^{n-2}
\frac{D^{*}F_{\t}^{n-k-1}(F_{\t}^{k+1}(\j)^{(1)})}{D_{*}F_{\t}^{n-k-1}(F_{\t}^{k+1}(\j)^{(1)})^{q}}
\beta^{(q)}_{p}(F_{\t}^{k}(\j),\t)\right|<C\subg \delta_{n}^{s-q+1}.
\end{equation}
Suppose that $\beta^{(q)}_{p}(F_{\t}^{k}(\j),\t)\neq 0$ for some integer $0\le k\le m-2$.
Then we have 
$d(F_{\t}^{k}(\j)^{(0)}, \zeta)<\delta_n$ and 
\begin{align*}
&d(F_{\t}^{k+1}(\j)^{(0)}, \criticalset(F_{\t}))\\
&\qquad \le d(F_{\t}^{k+1}(\j)^{(0)},
F_{\t}(\zeta))+d(F_{\t}(\zeta),
F_{\t}^{n}(\j)^{(0)})+
d(F_{\t}^{n}(\j)^{(0)},\criticalset(F_{\t}))<C\subg
\delta_n
\end{align*}
from (V1) and (V3).  This and (\ref{eqn-c1}) imply $|D^{*}F_{\t}(F_{\t}^{k+1}(\j)^{(1)})|<C\subg
\delta_n$,  and  hence
\[
\left|\frac{D^{*}F_{\t}^{m-k-1}(F_{\t}^{k+1}(\j)^{(1)})}{D_{*}F_{\t}^{m-k-1}(F_{\t}^{k+1}(\j)^{(1)})^{q}}
\beta^{(q)}_{p}(F_{\t}^{k}(\j),\t)\right|\le
C\subg \delta_n^{s-q+1}
\exp(-\lambda\subg(m-k-1)+2c\subg).
\]
Therefore we obtain  (\ref{est-bq}). 

The jet $\c(w;F_{\t})$ for $w\in \ball(\criticalset(F),\delta_n)$ 
does not depend on the parameter
$\t\in R$ because 
 $\ball(\zeta,\delta_n)\cap
\ball(\criticalset(F),\delta_n)=\emptyset$ from (V1) and 
(V2).
So we have 
\[
\|\partial_{p}(\c(F_{\t}^{n}(\j)^{(0)};F_{\t})^{(q)})\|<C\subg
\|\partial_{p}(F_{\t}^{n}(\j)^{(0)})\|< C\subg
\delta_{n}^{s-\nu}\quad\mbox{for $\nu+1\le p,q\le r-2$}
\]
  from  lemma \ref{lem:fdif}. From 
(\ref{est-bq}) and lemma \ref{est-ent}, it follows
\[
\left|\partial_{p}\left(F_{\t}^{m}(\j)^{(q)}-\c(F_{\t}^{n}(\j)^{(0)};F_{\t})^{(q)}\right)-
\beta^{(q)}_{p}(F_{\t}^{m-1}(\j),\t)\right|
<C\subg\delta_{n}^{s-q+1}.
\]
Denote by  $D\Psi(\t)_{q,p}$  the $(q,p)$-entree of the representation matrix of
$D\Psi(\t)$ with respect to the standard basis of $\Real^{r-2-\nu}$.  Then, from 
the estimates above, we have
 \begin{alignat}{2}
|D\Psi(\t)_{q,p}|&<C\subg \delta_n& \quad&\mbox{ if $q>p$,}\notag\\
|D\Psi(\t)_{q,p}|&<C\subg& &\mbox{ if $q\le p$, and }\notag\\
|D\Psi(\t)_{q,p}|&>C\subg^{-1}& &\mbox{ if $q=p$.}\notag
\end{alignat}
Now we can conclude  the lemma by an elementary argument.
\end{proof}

\subsection{Resolution of the flat contacts}
In this subsection, we prove theorem~\ref{th-ge-nonflat}.  Until the
last part of the proof, we fix $1\le
\ell\le \ell_0$ and put $\delta_{n}=\exp(-\lambda^{+}(\ell)n/\nu)$ for $n\ge 1$ as in the last subsection. Let $n$
be a large integer, $\zeta$ a point in the lattice $\lattice(\delta_{n}/20)$ and 
 $\j$ a jet in 
$\spl(n,\ell)$.
We denote, by
$Y_{0}(n,\ell,\j, \zeta)$ (resp. $Y_{1}(n,\ell,\j, \zeta)$),  the set of mappings
$F\in C^{r}(M,M)$ that satisfy
\begin{equation}\label{eqv}
F^{n-1}(\j)^{(0)}\in
\ball\left(\zeta,
 \delta_{n}/20\right)\qquad(\mbox{ resp. }
F^{n-1}(\j)^{(0)}\in
\ball\left(\zeta,
 \delta_{n}/5\right)  ).
\end{equation}
Below we  estimate 
\[
\mom_s\left(\base^{-1}_{G}(\ss_2(n,\ell,\j)\cap Y_{0}(n,\ell,\j,\zeta))\cap
\disk^{r}(d)
\right)\qquad \mbox{for  $G\in C^{r}(M,\torus)$ and $d>0$,}
\]
 where $\Phi_G$ and 
$\disk^{r}(d)$ are those defined by (\ref{PhiG}) and (\ref{defD}) respectively. 

Take a mapping $F$ in  
$\ss_2(n,\ell,\j)\cap Y_{0}(n,\ell,\j,\zeta)$ arbitrarily and consider the family
$F_{\t}$ defined by (\ref{deffam2}) in the last subsection. Note that the jet $\j$ belongs to $V_1(n,\ell;F)$
from the definition of $\ss_2(n,\ell,\j)$.  We  check that the conditions
(V1), (V2) and (V3) hold for  $\t\in R$ provided that 
$n$ is larger than some constant~$C\subg$.  
Since 
$F$ belongs to $\ss_{2}(n,\ell,\j)$, there exists a point ${w_0}\in \criticalset(F)$ such that
\begin{equation}\label{eqn-eqbv}
d_{\Jet}(F^{n}(\j),\c({w_0};F))\le 2B\subg
\delta_{n}^{(r-3)\nu}.
\end{equation}
Especially we have $d(F^{n}(\j)^{(0)},{w_0})<\rho\subg$ and 
$
\angle(F^{n}(\j)^{(1)},\c({w_0};F)^{(1)})< \rho\subg$, provided that $n$ is larger than some constant $C\subg$. 
It follows from the condition (C5) in 
the choice of  the constant $\rho\subg$
in subsection \ref{cnbl} that 
\begin{equation}\label{eqn-jdel}
d(F^{n-1}(\j)^{(0)},\criticalset(F))>\rho\subg.
\end{equation}
 Using (\ref{eqv}), we can see  
\[
d(\zeta,\criticalset(F))\ge
d(F^{n-1}(\j)^{(0)},\criticalset(F))-d(F^{n-1}(\j)^{(0)},\zeta)>\rho\subg-2B\subg
\delta_{n}^{(r-3)\nu}
>4\delta_n
\]
provided that 
$n$ is larger than some constant~$C\subg$. 
This implies that the critical set $\criticalset(F_\t)$ does not depend on
$\t\in R$.  Hence (V1), (V2) and (V3) follow from (\ref{eqv}), (\ref{eqn-eqbv}), (\ref{eqn-jdel})
and lemma
\ref{lem:fdif}, provided that 
$n$ is larger than some constant~$C\subg$.

Let $\Psi:R\to \Real^{r-\nu-2}$ be the mapping that we defined
in the last subsection. Note that the
conclusion of lemma \ref{lem-psiest} holds for this $\Psi$. Suppose that 
$F_{\s}$ belongs to
$\ss_2(n,\ell,\j)\cap Y_{0}(n,\ell,\j,\zeta)$ for a parameter $\s\in R$. 
Then, by definition, there exists a point $w_1\in \criticalset(F)$ such that 
$d_{\Jet}(F_{\s}^{n}(\j),\c({w_1};F_{\s}))<2 B\subg \exp(-\lambda^{+}(\ell)n (r-3))$. Since
$\c(\cdot;F_\s)=\c(\cdot;F):\ball(\criticalset(F),\rho\subg)\to \Jet^{r-2}\Curve$ is a
$C^1$mapping whose first-order differentials are bounded by some constant~$C\subg$, it follows
\begin{align*}
d_{\Jet}(F_{\s}^{n}(\j),\c(F^{n}_{\s}\j^{(0)};F_{\s}))&
\le d_{\Jet}(F_{\s}^{n}(\j),\c({w_1};F_{\s}))
+d_{\Jet}(\c({w_1};F_{\s}), \c(F_{\s}^{n}(\j)^{(0)};F_{\s}))\\
&<(1+C\subg) d_{\Jet}(F_{\s}^{n}(\j),\c({w_1};F_{\s})) 
<C\subg 
\delta_{n}^{\nu(r-3)}.
\end{align*}
Hence the image  $\Psi(\s)$ is contained in 
\[
\prod_{q=\nu+1}^{r-2}\left[-C\subg\delta_{n}^{\nu(r-3)-(s-q)},C\subg\delta_{n}^{\nu(r-3)-(s-q)}\right]
\subset \Real^{r-\nu-2}.
\]
We arrive at the estimate 
\[
\leb_{\Real^{r-\nu-2}}\{\t\in R\mid  F_{\t}\in \ss_2(n,\ell,\j)\cap Y_{0}(n,\ell,\j,\zeta)\}\le
C\subg \prod_{q=\nu+1}^{r-2}\delta_{n}^{\nu(r-3)-(s-q)},
\]
which holds uniformly for $F\in \ss_2(n,\ell,\j)\cap Y_{0}(n,\ell,\j,\zeta)$, provided that $n$
is larger than some constant $C\subg$.  

Now we apply lemma
\ref{lem-meas-est}. Fix a small number $0<T<1$ such that
\[
\max_{|t_q|\le T}\left\|\sum_{q=\nu+1}^{r-2}t_{q}\psi_{q,n,\zeta}\right\|_{C^s}\le
r\cdot \max_{\nu\le q\le r-2}\|\psi_{q}\|_{C^s}\cdot  T<\bound_{s}(d)
\]
where $\bound_{s}(d)$ is that in  lemma \ref{lem-mom}. Note that we can take $T$ independently of
$n$.  Put $X=\ss_2(n,\ell,\j)\cap Y_{0}(n,\ell,\j,\zeta)$ and $T_{i}=T$ in
lemma
\ref{lem-meas-est}. Then the assumption (\ref{eqn-small-norm}) holds from the choice of $T$, and the subset $Y$ in
the statement of lemma
\ref{lem-meas-est} is contained in 
$Y_{1}(n,\ell,\j,\zeta)$ from  the condition (V1) which we have proved above. Therefore we can obtain, as the
conclusion, 
\[
\frac{\mom_s  (\base^{-1}_{G}(\ss_{2}(n,\ell,\j)\cap Y_{0}(n,\ell,\j,\zeta))\cap
\disk^{r}(d))) }
{\mom_s(\base^{-1}_{G}(Y_{1}(n,\ell,\j,\zeta)))}\le 
C\subg T^{-r+\nu+2}\prod_{q=\nu+1}^{r-2}\delta_{n}^{\nu(r-3)-(s-q)}
\]
provided that $n$ is larger than some constant $C\subg$. 
The the subsets  $Y_{0}(n,\ell,\j,\zeta)$ for $\zeta\in
\lattice(\delta_{n}/20)$  cover 
$C^{r}(M,M)$ while the intersection multiplicity of the subsets 
$Y_{1}(n,\ell,\j,\zeta)$ for $\zeta\in
\lattice(\delta_{n}/20)$ is bounded by some absolute constant. Hence we can conclude that the measure 
 $\mom_s  (\base^{-1}_{G}(\ss_{2}(n,\ell,\j))\cap
\disk^{r}(d))$ is bounded by 
\begin{align*}
C\subg &T^{-r+\nu+2}\prod_{q=\nu+1}^{r-2}\delta_{n}^{\nu(r-3)-(s-q)}\\
&= C\subg T^{-r+\nu+2}
\exp\left((r-\nu-2)\left(-(r-3)+\frac{2s-r-\nu+1}{2\nu}\right)\lambda^{+}(\ell)n\right).
\end{align*}
The subset  $\ss_2$ is contained in the closed subset
\[
\ss'_{2}:=\bigcap_{n\ge B\subg}\bigcup_{\ell=1}^{\ell_0}
\bigcup_{\j\in \spl(n,\ell)}\ss_{2}(n,\ell,\j)
\]
from corollary \ref{cirqv}.
Hence the measure $\mom_s (\base^{-1}_{G}(
\ss'_{2})\cap \disk^{r}(d))$ is bounded by 
\[
 C\subg T^{-r+\nu+2}\sum_{\ell=1}^{\ell_0}\#\spl(n,\ell)
\cdot 
\exp\left((r-\nu-2)\left(-(r-3)+\frac{2s-r-\nu+1}{2\nu}\right)
\lambda^{+}(\ell)n\right)
\]
for sufficiently large $n$. From the estimate (\ref{card-q}) on the cardinality of 
$\spl(n,\ell)$ and the condition in the choice of $\lambda^{\pm}(\ell)$, this
converges to $0$ exponentially fast as $n\to \infty$. Thus we conclude $\mom_s 
(\base^{-1}_{G}(
\ss_{2})\cap \disk^{r}(d))=\mom_s
(\base^{-1}_{G}(
\ss'_{2})\cap \disk^{r}(d))=0$. As $d$ is an arbitrary positive number,  $\mom_s 
(\base^{-1}_{G}(
\ss_{2}))=0$ or $\ss_{2}$ is shy with respect to $\mom_s$. 

Suppose that $r\ge 19$. Then the inequality 
(\ref{rs}) holds for
$s=r+3$ and
$\nu=3$ and so
$\mom_{r+3}(\base_{G}^{-1}(\ss'_2))=0$ for any
$G\in C^{r}(M,\torus)$. This implies that  $\U\setminus \ss'_{2}$
is  dense.  Therefore  $\ss_2$ is contained in  the closed nowhere dense  subset 
$\ss'_2$. 

\appendix
\section{Proof of corollary \ref{lm-trans-meas}}
To see that  corollary \ref{lm-trans-meas} follows
from theorem
\ref{main_theorem2},  it is enough to show 
\begin{lemma}\label{lm-trans-meas-2}If
$X$ is a Borel subset in
$C^{r}(\man,\torus^2)$ that is timid for the class  $\mathcal Q_s^r$ of measures for some $s>r$, then the  subset 
\[
Y=\left\{F(z,t)\in C^{r}(\man\times [-1,1]^{k},\torus)\mid \leb_{\Real^k}(\{t\in [-1,1]^k\mid F(\cdot,t)\in
X\})>0\right\},
\]
 is timid for the class of 
 Borel finite measures on  $C^{r}(\man\times [-1,1]^{k},\Real^2)$ that are quasi-invariant
along
$C^{s}(\man\times [-1,1]^{k},\Real^2)$. 
\end{lemma}
\begin{proof}
Take a mapping $G\in C^{r}(\man\times [-1,1]^{k},\torus)$ and put $G_0(z)=G(z,0)$. We define the
mapping
\[
P(f, \t):= G(\cdot, \t)-G_{0}(\cdot)+f(\cdot,\t): C^{r}(\man\times [-1,1]^{k},\Real^2)\times [-1,1]^{k}\to
C^{r}(\man,\Real^2),
\]
so that $\base_{G_0}\circ P(f,\t)=G(\cdot,\t)+f(\cdot,\t)$. 
Let  $\mon$ be a Borel  finite  measure on  $C^{r}(\man\times [-1,1]^{k},\Real^2)$ that is
quasi-invariant along
$C^{s}(\man\times [-1,1]^{k},\Real^2)$. Then the measure 
$(\mon\times \leb_{\Real^k}|_{[-1,1]^k})\circ P^{-1}$ on $C^{r}(\man,\Real^2)$ belongs to 
$\mathcal Q_{s}^{r}$ and so we have $(\mon\times \leb_{\Real^k})( (\base_{G_0}\circ P)^{-1}(X))=0$ from the
assumption. This  and Fubini theorem imply $\mon\circ \Phi_{G}^{-1}(Y)=0$ and hence the claim of the lemma. 
\end{proof}

\section{Proof of lemma \ref{lem-mom}}
We use the definitions and results in the book\cite{Skorohod} by Skorohod. We consider the
functions
$e_{nm}(x,y)=\exp\left(2\pi\sqrt{-1} (n x+m y)\right)$ for  $n,m\in \Z$,  as a complete orthonormal basis of the
space
$L^{2}(\torus,\vol)$. Let $A: L^{2}(\torus,\vol)\to L^{2}(\torus,\vol)$ be the operator defined
 by
\[
A\left( \sum_{(n,m)\in \Z^{2}}a_{nm}e_{nm}\right)=\sum_{(n,m)\in \Z^{2}}(n^{2}+m^{2}+1)^{-1/2}a_{nm}e_{nm}.
\] 
Let $\mon$ be  the Gaussian measure\cite[\S 5 ]{Skorohod}
on
$L^{2}(\torus, \vol)$ whose characteristic function is
$\Theta(\psi)=\exp(-(1/2)(A^{2s-3}\psi,\psi)_{L^{2}})$. Then $\mon$ is supported
on the Sobolev space
$W^{s-3}:=A^{s-3}(L^{2}(\torus,\vol))$. We can see, from
\cite[\S 16 theorem 2]{Skorohod}, that  $\mon$ is quasi-invariant along 
$W^{s-(3/2)}\supset C^{s-1}(\torus,\Real)$ and it holds
\begin{align*}
\frac{d(\mon\circ \tau_{\psi}^{-1})}{d\mon}(\varphi)&=\exp
\left((A^{-s}\psi,A^{-s+3}\varphi)_{L^{2}}-(1/2)\left\|A^{-s+(3/2)}\psi\right\|^2_{L^{2}}\right)\\
&\le \exp(\|\psi\|_{W^{s}}\cdot \|\varphi\|_{W^{s-3}})
\end{align*}
for $\psi\in
W^{s}$ and $\mon$-a.e. $\varphi\in W^{s-3}$.

We show that the measure $\mon$ is
actually supported on
$C^{s-3}(\torus,\Real)$. We  follow the argument in the proof of the fact that the measure corresponding to 
Brownian motion is supported on the space of continuous paths\cite{IW}.  Let
$\varphi^{(s-3)}$ be one of the 
$(s-3)$-th partial differentials of
$\varphi$. Denoting the expectation with respect to the measure
$\mon$ by $E(\cdot)$, we have
\[
E(|\varphi^{(s-3)}(z)-\varphi^{(s-3)}(w)|^{5})
\le const \cdot d(w,z)^{5/2}
\]
because the distribution of
$\varphi^{(s-3)}(z)-\varphi^{(s-3)}(w)$ is  a Gaussian distribution with average 0 and variance 
bounded  by
\[
\sum_{(n,m)\in \Z^2}\left(\min\{2,(n^{2}+m^{2}+1)^{1/2}d(z,w)\}
(n^{2}+m^{2}+1)^{-3/4}\right)^{2}
\le const. \cdot d(z,w).
\]
By Borel-Cantelli lemma,   
there is a constant
$i_0>0$ for $\mon$-almost every $\varphi$ such that 
\[
\left|\varphi^{(s-3)}(2^{-i}p, 2^{-i}q)-\varphi^{(s-3)}(2^{-i}p', 2^{-i}q')\right|^{5}\le
2^{-i/3}
\]
for $i>i_0$ and  $p,q,p',q'\in \Z$ such that $|p-p'|\le 1$ and
$|q-q'|\le 1$. This implies that $\varphi^{(s-3)}$ is continuous for  $\mon$-almost every $\varphi$
and hence $\mon$ is supported on $C^{s-3}(\torus, \Real)$. 

As   $C^{s-3}(\torus,\Real^2)$ is naturally identified with  $C^{s-3}(\torus,\Real)\times C^{s-3}(\torus,\Real)$,
we regard  the product 
$\mon\times\mon$ as a measure on $C^{s-3}(\torus,\Real^2)$. Put
$\mom_s=(\mon\times\mon)\circ
\pi^{-1}$ where  $\pi:C^{s-3}(\torus,\Real^2)\to
C^{s-3}(M,\Real^2)$ is the mapping that corresponds to the restriction to $M$.  Then $\mom_s$ satisfies the
conditions in the lemma.

\end{document}